\documentclass[12pt,dvips]{article}

\usepackage{rotating}
\usepackage{graphicx,latexsym,amsmath,amsfonts,amssymb,mathrsfs}

\usepackage{tikz}
\usetikzlibrary{shapes,snakes}

\newtheorem{thm}{Theorem}

\newtheorem{cnj}[thm]{Conjecture}
\newtheorem{cor}[thm]{Corollary}

\newtheorem{lem}[thm]{Lemma}

\newtheorem{qst}[thm]{Question}

\def\a{{\alpha}}

\def\k{{\kappa}}

\def\r{{\rho}}
\def\s{{\sigma}}
\def\t{{\tau}}

\def\x{{\chi}}

\def\bv{{\mathbf v}}

\def\bx{{\mathbf x}}

\def\bP{{\mathbf P}}
\def\bT{{\mathbf T}}

\def\zN{{\mathbb N}}

\def\cG{{\cal G}}

\def\cP{{\cal P}}

\def\cT{{\cal T}}

\def\sNP{{\sf NP}}
\def\sP{{\sf P}}

\def\deg{{\sf deg}}
\def\diam{{\sf diam}}
\def\dist{{\sf dist}}

\def\mt{{\emptyset}}

\def\pr{{\prime}}
\def\rar{{\rightarrow}}

\def\sqr#1#2{{\vcenter{\hrule height.#2pt
        \hbox{\vrule width.#2pt height#1pt \kern#1pt
                \vrule width.#2pt}
        \hrule height.#2pt}}}
\def\gbox{{\mathchoice\sqr34\sqr34\sqr{2.1}3\sqr{1.5}3}}
\def\L2{{L\gbox L}}

\def\Mx{{\rm Max.}}
\def\st{{\rm s.t.}}

\def\proof{{\noindent {\it Proof}.\ \ }}
\def\pf{{\hfill$\Box$\\ \medskip}}

\providecommand{\binom}[2]{{#1\choose#2}}

\begin{document}

% ##########################################################################
% ##########################################################################
%
%       TITLE PAGE
%
\title{A Linear Optimization Technique for Graph Pebbling\\
\ \\}

\author{
Glenn Hurlbert\thanks{Much of this work was carried out while vising the
Centre de Recerca Matematica during the research program on Variational
Analysis and Optimization, Fall 2010.}\\
School of Mathematics and Statistics\\
Arizona State University\\
Tempe, AZ 85287 USA\\
\texttt{hurlbert@asu.edu}
}
\maketitle

%\vspace{1.0 in}

% ##########################################################################
% ##########################################################################
%
%       ABSTRACT
%
\begin{abstract}
Graph pebbling is a network model for studying whether or not a given
supply of discrete pebbles can satisfy a given demand via pebbling moves.
A pebbling move across an edge of a graph takes two pebbles from one
endpoint and places one pebble at the other endpoint; the other pebble
is lost in transit as a toll.
It has been shown that deciding whether a supply can meet a demand on
a graph is \sNP-complete.
The pebbling number of a graph is the smallest $t$ such that every supply
of $t$ pebbles can satisfy every demand of one pebble.
Deciding if the pebbling number is at most $k$ is $\Pi_2^\sP$-complete.

In this paper we develop a tool, called the Weight Function Lemma, for
computing upper bounds and sometimes exact values for pebbling numbers 
with the assistance of linear optimization.
With this tool we are able to calculate the pebbling numbers of much
larger graphs than in previous algorithms, and much more quickly as well.
We also obtain results for many families of graphs, in many cases by hand, 
with much simpler and remarkably shorter proofs than given in previously
existing arguments (certificates typically of size at most the number of
vertices times the maximum degree), especially for highly symmetric graphs.

Here we apply the Weight Function Lemma to several specific
graphs, including the Petersen, Lemke, $4^{\rm th}$ weak Bruhat, Lemke 
squared, and two random graphs, as well as to a number of infinite 
families of graphs, such as trees, cycles, graph powers of cycles, cubes,
%Kneser graphs, 
and some generalized Petersen and Coxeter graphs.
This partly answers a question of Pachter, et al., by computing the pebbling
exponent of cycles to within an asymptotically small range.
It is conceivable that this method yields an approximation
algorithm for graph pebbling.

\end{abstract}

\newpage

% ##########################################################################
% ##########################################################################
%
%       INTRODUCTION
%
\section{Introduction}\label{Intro}

Graph pebbling is like a number of network models, including network flow,
transportation, and supply chain, in that one must move some commodity
from a set of sources to a set of sinks optimally according to certain
constraints.
Network flow constraints restrict flow along edges and conserve flow through
vertices, and the goal is to maximize the amount of commodity reaching the sinks.
The transportation model includes per unit costs along edges and aims to 
minimize the total cost of shipments that satisfy the source supplies and 
sink demands.
At its simplest, the supply chain model ignores transportation costs while
seeking to satisfy demands with minimum inventory.
The graph pebbling model introduced by Chung \cite{Chun} also tries to 
meet demands with minimum inventory, but constrains movement across an edge
by the loss of the commodity itself, much like an oil tanker using up the
fuel it transports, not unlike heat or other energy dissipating during transfer.

Specifically, a {\em configuration} $C$ of pebbles on the vertices of a 
connected graph $G$ is a function $C:V(G)\rar\zN$ (the nonnegative integers), 
so that $C(v)$ counts the number of pebbles placed on the vertex $v$.
We write $|C|$ for the {\em size} $\sum_vC(v)$ of $C$; i.e. the number of
pebbles in the configuration.
A {\em pebbling step} from a vertex $u$ to one of its neighbors $v$ reduces
$C(u)$ by two and increases $C(v)$ by one (so that one can think of it as
moving one pebble at the {\em cost} of another as toll).
Given two configurations $C$ and $D$ we say that $C$ is $D$-{\em solvable}
if some sequence of pebbling steps converts $C$ to $D$.
In this paper we study the traditional case in which the target distribution
consists of a single pebble at some {\em root} vertex $r$ (one can peruse 
\cite{Hurl1,Hurl2,GPP} for a wide array of variations on this theme).
We are concerned with determining $\pi(G,r)$, the minimum number $t$ of
pebbles so that every configuration of size $t$ is $r$-solvable.
Then the {\em pebbling number} of $G$ equals $\pi(G)=\max_r\pi(G,r)$.
Alternatively, $\pi(G)$ is one more than the maximum $s$ such that
there is some root $r$ and some size $s$ configuration $C$ so that 
$C$ does not solve $r$.
The primary focus of this paper is to exploit this duality with newly
discovered algebraic constraints.

\medskip

% ##########################################################################
%
%       CALCULATING
%
\subsection{Calculating Pebbling Numbers}\label{Calc}

Given a graph $G$, configuration $C$, and root $r$, one can ask how
difficult it is to determine if $C$ solves $r$.
In \cite{HurKie} it was determined that this problem is \sNP-hard.
Subsequently, \cite{MilCla,Wats} proved that the problem is \sNP-complete,
with \cite{MilCla} showing further that answering the question
``is $\pi(G)\le k$?'' is $\Pi_2^\sP$-complete (and hence both \sNP-hard
and co\sNP-hard, and therefore in neither \sNP\ nor co\sNP\ unless
\sNP\ $=$ co\sNP).
Finding classes of graphs on which we can answer more quickly is
therefore relevent, and there is some evidence that one can be
successful in this direction.
Besides what we share in this introduction, we show later that many
graphs can have very short certificates that $\pi(G)\le k$.

The $r$-unsolvable configuration with one pebble on every vertex other than 
the root $r$ shows that $\pi(G)\ge n$, where $n=n(G)$ denotes the number
of vertices of $G$.
In \cite{PaSnVo} it is proved that graphs of diameter two satisfy
$\pi(G)\le n+1$, with a characterization separating the two classes
({\it Class 0} means $\pi(G)=n$ and {\it Class 1} means $\pi(G)=n+1$)
given in \cite{BlaSch,ClHoHu}.
One of the consequences of this is that 3-connected diameter two graphs
are Class 0.
As an extension it is proved in \cite{CzHuKiTr} that $2^{2d+3}$-connected 
diameter $d$ graphs are also Class 0, and they use this result to show
that almost every graph with significantly more than $n(n\lg n)^{1/d}$ 
(for any fixed $d$) edges is Class 0.
Consequently, it is a very (asymptotically) small 
collection of graphs that cause all the problems.

Knowing the pebbling number of a graph and actually solving a particular
configuration are two different things, as even a configuration that is
known to be solvable (say, one of size equal to the pebbling number) can
be difficult to solve.
Evidence that most configurations are not so difficult, though, comes in the
following form.
The work of \cite{BeBrCzHu} shows that every infinite graph sequence
$\cG=(G_1,G_2,\ldots,G_N,\ldots)$ has a {\em pebbling threshold}
$\t_\cG:\zN\rar\zN$, which yields the property that almost every 
configuration $C_N$ on $G_N$ of size $|C_N|\gg\t(N)$ is solvable 
(and almost every configuration of size $|C_N|\ll\t(N)$ is not).
In papers such as \cite{BekHur,CzyHur1,CzyHur2} we find that $\t_\cG(N)$
is significantly smaller than $\pi(G_N)$ --- for example, $\sqrt{N}$ as
opposed to $N$ for the complete graph $K_N$, and roughly $N2^{\sqrt{\lg N}}$
as opposed to $2^{N-1}$ for the path $P_N$.
Moreover, the proof techniques reveal that almost all of these solvable
configurations can be solved {\em greedily}, meaning that every pebbling
step reduces the distance of the pebble to the root.
So the hardness of the problem stems from a rare collection of configurations.

With these results as backdrop, \cite{BekCus} presents a polynomial
algorithm for determining the solvability of a configuration on diameter two
graphs of connectivity some fixed $\k$.
Furthermore, \cite{Sieb} contains an algorithm that calculates pebbling
numbers, and is able to complete the task for every graph on at most 9 vertices.
Also, the proof in \cite{Chun} that the $d$-dimensional cube is of Class 0
is a polynomial algorithm (actually bounded by its number of edges $n\lg n$).
Along these lines, our main objective is to develop algorithmic tools 
that will in a reasonable amount of time yield good upper bounds 
on $\pi(G)$ for much larger graphs, and in particular decide in some cases
whether or not a graph is of Class 0.

This latter determination is motivated most by the following conjecture
of Graham in \cite{Chun}.
For graphs $G$ and $H$, let $G\gbox H$ denote the {\em Cartesian product}
whose vertices are $V(G\gbox H)=V(G)\times V(H)$, with edges
$(u,x)\sim (v,x)$ whenever $u\sim v$ in $G$ and
$(u,x)\sim (u,y)$ whenever $x\sim y$ in $H$.

\begin{cnj}\label{graham}
{\bf (Graham)}
Every pair of graphs $G$ and $H$ satisfy $\pi(G\gbox H)\le\pi(G)\pi(H)$.
\end{cnj}

\begin{figure}
\centerline{\includegraphics[height=2.0in]{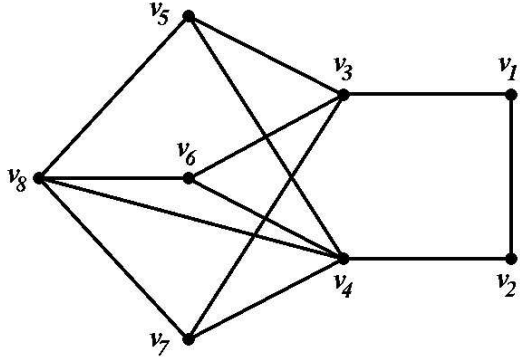}}
\caption{The Lemke graph}\label{lemke}
\end{figure}

The conjecture has been verified for many graphs; see \cite{Hers} for
the most recent work.
However, as noted in \cite{Hurl3}, there is good reason to suspect that 
$\L2$ might be a counterexample to this conjecture, if one exists, 
where $L$ is the {\em Lemke} graph of Figure \ref{lemke}.
Since $L$ is Class 0, Graham's conjecture requires that $\L2$ is also,
but it is a formiddable challenge to compute the pebbling number of a
graph on 64 vertices.
One hopes that graph structure and symmetry will be of use, but purely
graphical methods have failed to date.
The methods of this paper represent the first strides toward the
computational resolution of the $\$64$ question\footnote{Yes, I'll pay
if you beat me to it!}, ``Is $\pi(\L2)\le 64$?''.
Certainly, these methods alone will not suffice%
\footnote{We obtain evidence that $\pi(\L2)\le 108$ in Theorem \ref{Lem2} ---
in fact, for one root $r$ we show $\pi(\L2,r)\le 68$.}% 
, but if they produce a decent upper bound then the methods of \cite{Sieb} 
might be able to finish the job.

\medskip

% ##########################################################################
%
%       RESULTS
%
\subsection{Results}\label{Results}

The main tool we develop is the Weight Function Lemma \ref{wfl}.
This lemma allows us to define a (very large) integer linear 
optimization problem that yields an upper bound on the pebbling number.
This has several important consequences, including the following.

\begin{enumerate}
\item
The pebbling numbers of reasonably small graphs often can be computed easily.
Moreover, it is frequently the case that the fractional relaxation suffices 
for the task, allowing the computation for somewhat larger graphs.
\item
It is also common that only a small portion of the constraints are
required, expanding the pool of computable graphs even more.%
\footnote{We present some findings along these lines in Section
\ref{Random}, with graphs on $15$ and $20$ vertices.}
One can restrict the types of constraints to greedy, bounded depth, 
and so on, with great success, seemingly because of the comments above.
Potentially, this allows one to begin to catalog special classes of graphs 
such as Class 0, (semi-)greedy, and tree-solvable.
\item
The dual solutions often yield very short certificates of the results,
in most cases quadratic in the number of vertices, and usually at most 
the number of vertices times the degree of the root.
These certificates are remarkably simple compared to the usual
solvability arguments that chase pebbles all over the graph in a
barrage of cases.
One can sometimes find such certificates for infinite families of graphs
by hand, without resorting to machine for more than the smallest
one or two of its members.  
This was our approach in Section \ref{Classes}, for example.
\item
Our method gives trivial proofs of 
\begin{enumerate}
\item
$\pi(C_{2k})=2^k$ and $\pi(C_{2k+1})=2\lfloor 2^{k+1}/3\rfloor+1$,
which we write as $\lceil (2^{k+2}-1)/3\rceil$, and
\item
$C_n^{(k)}$ is Class 0 for $k\ge n/2(\lg n-\lg\lg n)$, where $G^{(k)}$ 
denotes the $k^{\rm th}$ {\it graph power} of $G$ (as opposed to the 
{\it Cartesian} power $G^k$).
This answers a question of \cite{PaSnVo}, who defined the {\it pebbling
exponent} of $G$ minimum such $e=e_\pi(G)$ for which $\pi(G^{(e)})=n(G)$.
Thus $e_\pi(G)\le n/2(\lg n-\lg\lg n)$ (see Theorem \ref{PebbExpo}),
which is fairly close to the obvious lower bound of $n/2\lg n$.
\end{enumerate}
\end{enumerate}

In this paper we apply the Weight Function Lemma to several specific
graphs, including the Petersen, Lemke, $4^{\rm th}$ weak Bruhat, Lemke 
squared, and two random graphs, as well as to a number of infinite 
families of graphs, such as trees, cycles, graph powers of cycles, cubes,
%Kneser graphs, 
and some generalized Petersen and Coxeter graphs.

%\bigskip

% ##########################################################################
% ##########################################################################
%
%       WEIGHT FUNCTION LEMMA
%
\section{The Weight Function Lemma}\label{WFL}

Let $P_n$ be the path $v_1v_2\cdots v_n$ on $n$ vertices.
Then $\pi(P_n)=2^{n-1}$ is easily proved by induction.
In particular, any configuration of at least $2^{n-1}$ pebbles solves $v_1$.
But one can say more about smaller $v_1$-solvable configurations as well,
with the use of a weight function $w$.
Define $w$ on $V(G)$ by $w(v_{n-i})=2^i$, and extend the weight function 
to configurations by $w(C)=\sum_{v\in V}w(v)C(v)$.
Then a pebbling step can only preserve or decrease the weight of a 
configuration.
Since the weight of a configuration with a pebble on $v_1$ is at least
$2^{n-1}$, we see that $2^{n-1}$ is a lower bound on every $v_1$-solvable
configuration.
In fact, induction shows that every $v_1$-unsolvable configuration has 
weight at most $2^{n-1}-1$, which equals $\sum_{i=2}^nw(v_i)$.
That is, this inequality characterizes $v_1$-unsolvable 
configurations on $P_n$.
The Weight Function Lemma \ref{wfl} generalizes this result on trees,
and we explore the applications of the lemma in the following sections.

\medskip

% ##########################################################################
%
%       LINEAR OPTIMIZATION
%
\subsection{Linear Optimization}\label{LinOpt}

Let $G$ be a graph and $T$ be a subtree of $G$ rooted at vertex $r$, with
at least two vertices.
For a vertex $v\in V(T)$ let $v^+$ denote the {\em parent} of $v$; i.e.
the $T$-neighbor of $v$ that is one step closer to $r$ (we also say that
$v$ is a {\em child} of $v^+$).
We call $T$ a {\em strategy} when we associate with it a nonnegative, nonzero
weight function $w$ with the property that $w(r)=0$ and $w(v^+)=2w(v)$
for every other vertex that is not a neighbor of $r$ (and $w(v)=0$ for
vertices not in $T$).
Let $\bT$ be the configuration with $\bT(r)=0$, $\bT(v)=1$ for all 
$v\in V(T)$, and $\bT(v)=0$ everywhere else.
With this notation note that the path result above can be restated:
$C$ is $v_1$-unsolvable if and only if $C(v_1)=0$ and $w(C)\le w(\bT)$,
where $T$ is the strategy $T=P_n$ with associated weight function $w$.

\begin{lem}{\sf [Weight Function Lemma]}\label{wfl}
Let $T$ be a strategy of $G$ rooted at $r$, with associated weight function $w$.
Suppose that $C$ is an $r$-unsolvable configuration of pebbles on $V(G)$.
Then $w(C)\le w(\bT)$.
\end{lem}

\proof
By contrapositive and induction.
The base case is when $T$ is a path, which is proved above.
Suppose $w(C)>w(\bT)$, let $y$ be a leaf of $T$, and define $P$ to be the 
path from $y$ to $r$ in $T$, with $P_y$ being the subpath from $y$ to its 
closest vertex $x$ on $P$ of degree at least 3 in $T$ (or $r$ if none exists).
Denote by $T^\pr$ the tree $T-P_y+x$, and among all such $r$-unsolvable 
configurations, choose $C$ to be the one having largest weight on $T^\pr$.
The restriction $w^\pr$ of $w$ to $T^\pr$ witnesses that $T^\pr$ is a strategy 
for the root $r$, so induction requires that $w^\pr(C)\le w^\pr(\bT^\pr)$.
Likewise, the restriction $w_y$ of $w$ to $P_y$ witnesses that $P_y$
is a strategy for the root $x$.
Because $w(C)=w^\pr(C)+w_y(C)$ and $w(\bT)=w^\pr(\bT^\pr)+w_y(P_y)$,
we must have $w_y(C)>w_y(P_y)$ which by induction means that the
restriction of $C$ to $P_y-x$ solves $x$ (and so $x\not= r$).
Let $C_x$ be the resulting configuration after moving a pebble from 
$P_y-x$ to $x$.
Since $C$ is $r$-unsolvable, so is $C_x$.
Now $w(C_x)=w(C)$, but $w^\pr(C_x)>w^\pr(C)$, which contradicts the initial
choice of $C$.
\pf

For a graph $G$ and root vertex $r$, let $\cT$ be the set of all $r$-strategies
in $G$, and denote by $z_{G,r}$ the optimal value of the integer linear 
optimization problem $\bP_{G,r}$:
\begin{equation}\label{ILOP}
\Mx\ \sum_{v\not= r}C(v)\ \st\ w(C)\le w(\bT),\ {\rm and\ }\bT\in\cT\ 
{\rm with\ witnessing\ weight\ function\ }w\ .
\end{equation}
We also let $\hat{z}_{G,r}$ be the optimum of the relaxation, which allows
configurations to be rational.
We will find the relation $z_{G,r}\le \lfloor \hat{z}_{G,r}\rfloor$ useful
at times.
The following corollary is straightforward.

\begin{cor}\label{Upper}
Every graph $G$ and root $r$ satisfies $\pi(G,r)\le z_{G,r}+1$.
\end{cor}

\proof
By definition, the pebbling number is one more than the size of the
largest unsolvable configuration.
\pf

Until now, one could only use trees in an individual manner:
$\pi(G,r)\le\pi(T,r)$ for every spanning tree $T$ rooted at $r$.
The Weight Function Lemma allows one to consider all subtrees rooted at $r$
(not only spanning trees) simultaneously, which we will see is significantly
more powerful.
One strength of the method is that the relaxation frequently has an integer
optimum.
This means that the dual solution will point out which tree constraints
certify the result, and because the dual problem has only $n(G)-1$
constraints there are at most that many such trees in the certificate.
Experience has shown, however, that usually one can find a certificate with
only $\deg(r)$ trees (or sometimes a few extra).
We will see this behavior starting in Section \ref{GenApps}.

\medskip

% ##########################################################################
%
%       BASIC APPLICATIONS
%
\subsection{Basic Applications}\label{BasicApps}

We begin with the pebbling number of trees, whose formula was first
discovered and proved in \cite{Chun}.
View a tree $T$ with root $r$ as a directed graph with every edge directed
toward $r$.
Then a {\em path partition} $\cP$ of $T$ is a set of edge-disjoint directed
paths whose union is $T$.
One path partition {\em majorizes} another if its nonincreasing sequence
of path lengths majorizes that of the other.
A path partition is {\em maximum} if it majorizes all others.
We can use Corollary \ref{Upper} to give a new proof of the following 
result of \cite{Chun}.

\begin{thm}\label{tree}
For a tree $T$ and root $r$ we have $\pi(T,r)=\sum_{P\in\cP}2^{e_P}-|\cP|+1$,
where $e_P$ denotes the length (number of edges) of $P$.
\end{thm}

\proof
We begin by showing that a maximum size $r$-unsolvable configuration has
pebbles on leaves only, and in fact on all leaves.
Indeed, if $C$ has a pebble on the nonleaf $x$, then we define a
{\em pushback} of $C$ at $x$ to be any configuration obtained by removing
the $C(x)$ pebbles from $x$, adding $2C(x)+1$ pebbles to one of the
children of $x$, and adding $1$ pebble to all other children of $x$.
Certainly, if $C$ is $r$-unsolvable and has no pebbles past $x$ (on the
subtree of $T$ rooted at $x$, minus $x$ itself), then the pushback will 
also be $r$-unsolvable, and thus satisfy the constraints of $\bP_{T,r}$.
It will also be larger than $C$.
The configuration $C^*$ that places $2^{e_P}-1$ pebbles on the leaf of the
path $P\in\cP$ is one possible result of pushing back the empty configuration,
and so satisfies the constraints of $\bP_{T,r}$.
Hence $\pi(T,r)\ge \sum_{P\in\cP}2^{e_P}-|\cP|+1$.

For the upper bound we prove that $C^*$ is optimal by using induction
to show that the optimal configuration has $2^{e_P}-1$ pebbles on the leaf
$y_P$ of $P$ for every $P\in\cP$.
This is true if $|\cP|=1$, so suppose $|\cP|\ge 2$ and let $Q$ denote
one of the paths in $\cP$ whose leaf has the highest weight in $w$, with 
a tie going to one of the shortest length.
This is to guarantee that the graph $T-Q+x$, where $x$ is the root of $Q$,
is a tree (i.e. is connected).
In order to maximize the number of pebbles that satisfy 
$\sum_{P\in\cP}w(y_P)C(y_P)\le\sum_{P\in\cP}w(P)$ we would transfer
as many pebbles from $y_Q$ to other leaves as possible because their
weights are at most $w(y_P)$ and we could add extra pebbles when the
weight is smaller.
But by induction on $T-Q+x$ we know from the constraint 
$\sum_{P\in\cP-\{Q\}}w(y_P)C(y_P)\le\sum_{P\in\cP-\{Q\}}w(P)$
that each $C(y_P)\le 2^{e_P}-1$, with equality for all $P\not=Q$ if and
only if $C$ is maximum.
Therefore, since $w(P)=w(y_P)(2^{e_P}-1)$ for all $P\in\cP$, we have
\begin{eqnarray*}
w(y_P)C(y_P)&\le&\sum_{P\in\cP}w(P)-\sum_{P\in\cP-\{Q\}}w(y_P)(2^{e_P}-1)\\
&=&w(Q)\\
&=&w(y_Q)(2^{e_Q}-1)\ ,
\end{eqnarray*}
which implies that $C(y_P)\le 2^{e_Q}-1$.
\pf

\begin{figure}
\centerline{\includegraphics[height=2.25in]{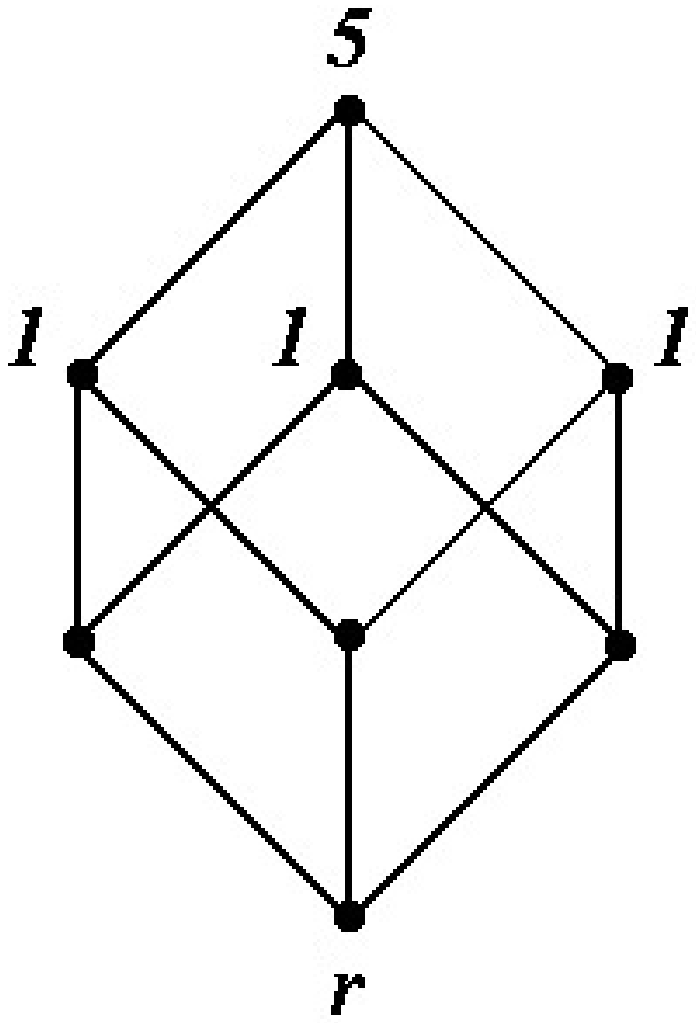}
$\qquad$\includegraphics[height=2.25in]{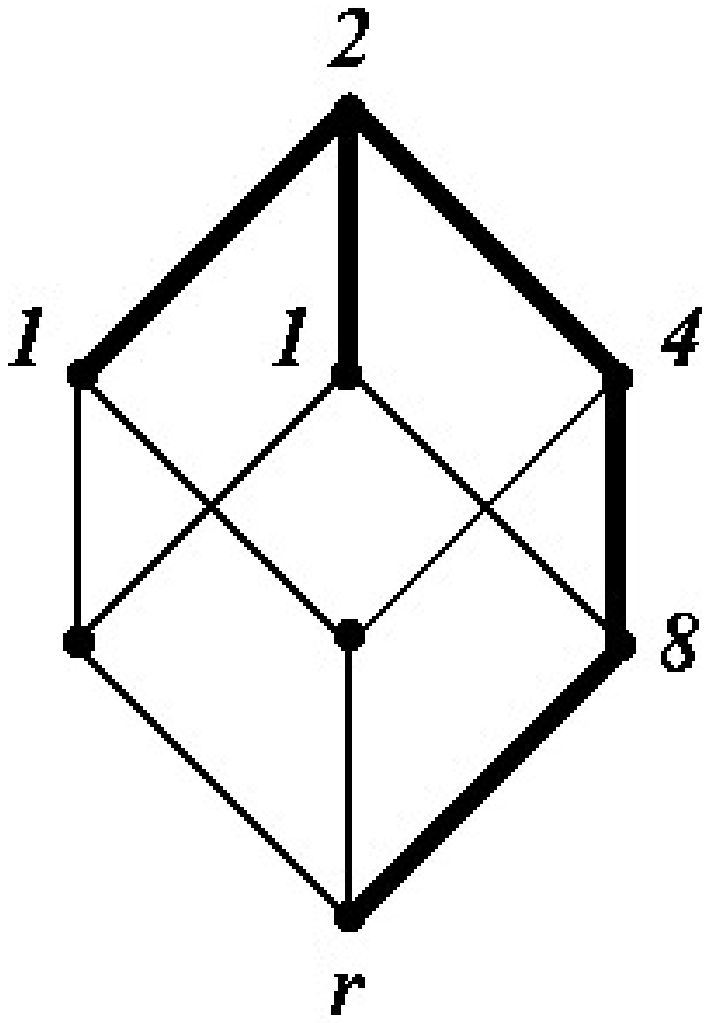}
$\qquad$\includegraphics[height=2.25in]{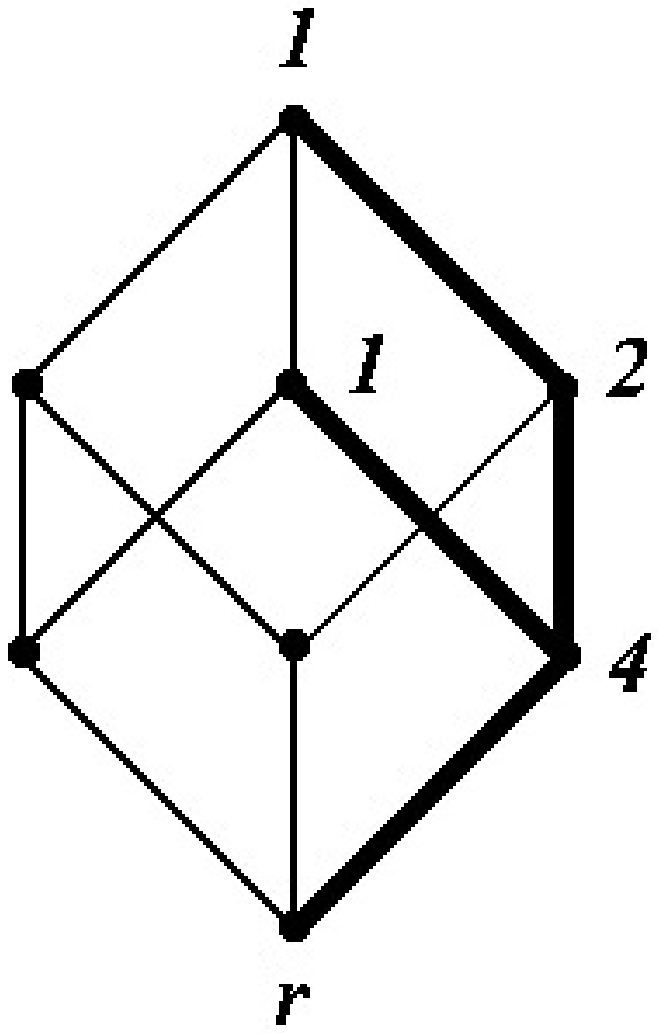}}
\caption{A solvable configuration (left) not recognized by any tree strategy;
a canonical strategy (center) used for certifying $\pi(Q^3)\le 9$;
a simplified nonbasic strategy (right) that does the same}
\label{3cube}
\end{figure}

A slight weakness of these tree constraints is that they do not classify
unsolvable configurations on trees the way that they do on paths.
This is because they let in a few solvable configurations.
For example, consider the star on four vertices with one of its leaves as
root $r$.
Then the configuration with $2$ pebbles on each of the other two leaves
is $r$-solvable and satisfies all tree contraints.
Since it is the average of the two $r$-unsolvable configurations that 
place either $1$ and $3$ or $3$ and $1$ on those other leaves, it cannot
be cut out by the tree constraints that don't cut out at least one of 
these two other constraints.
In this case it doesn't hurt us, since the strategy bound yields $\pi\le 5$ 
and the actual pebbling number is 5, but it can cause trouble on graphs in general.
For example, we know that the 3-cube $Q^3$ in Figure \ref{3cube} has
pebbling number 8, so that the shown configuration $C$ is solvable (pebbles
from to top must be split in two directions in its solution).
However, no strategy recognizes its solution, and Corollary \ref{Upper}
yields only $\pi(Q^3)\le 9$ (the three rotations of the strategy in the 
center in Figure \ref{3cube} certify this).
One can see where the aforementioned star appears in the Figure \ref{3cube}
configuration on $Q^3$ 
and is exploited accordingly: moving pebbles from the $5$ along one edge
yields a $(3,1)$ configuration, while splitting the moves along two edges
yields a $(2,2)$ configuration.

% ##########################################################################
% ##########################################################################
%
%       GENERAL APPLICATIONS
%
\section{General Applications}\label{GenApps}

In this section we illustrate the method more fully by presenting 
short proofs of both known and new results.
We begin by relaxing strategies in the following way.
We now use the term {\it basic} to describe the strategies as
currently defined.
A {\it nonbasic strategy} will satisfy the inequality $w(v^+)\ge 2w(v)$
in place of the equality used in a basic strategy (see Figure \ref{3cube}).
The following lemma shows that we can use nonbasic strategies in an
upper bound certificate since they are conic combinations of a nested
family of basic strategies.
Thus the use of nonbasic strategies can simplify and shorten certificates
significantly.

\begin{lem}\label{nonbasic}
If $T$ is a nonbasic strategy for the rooted graph $(G,r)$, then there
exists basic strategies $T_1,\ldots,T_k$ for $(G,r)$ and nonnegative
constants $c_1,\ldots,c_k$ so that $T=\sum_{i=1}^kc_iT_i$.
\end{lem}

\proof
We use induction, as the result is true when $T$ has two vertices
since $T$ is basic then.
Given $T$, let $S$ be a basic strategy on the edge set of $T$, 
define $c$ to be the largest constant for which $cS\le T$, and
denote $T^\pr=T-cS$.
Then some vertex $v$ of $G$ satisfies $cS(v)=T(v)$, so $T^\pr$ has 
fewer vertices than $T$.
Also, because $S$ is basic, any vertex $u$ whose unique $ur$-path
contains $v$ also satisfies $cS(u)=T(u)$, which means that $T^\pr$
is connected, and hence a strategy.
Moreover, $T^\pr$ is nonbasic since every nonneighbor $x$ of $r$ has
$T^\pr(x^+)=T(x^+)-cS(x^+)\ge 2T(x)-2cS(x)=2T^\pr(x)$.
By induction, $T^\pr$ is a conic combination of basic strategies,
and so therefore is $T$.
\pf

We use conic combinations of strategies to derive, for some $\a$, the inequality 
$|C|=\sum_{v\not=r}C(v)\le \a$ for $r$-unsolvable configurations $C$.
From this we surmise that $\pi(G)\le\lfloor \a\rfloor+1$.
Instead of writing our strategies algebraically, it will be somewhat
easier to show them graphically.
We will display them so as to derive 
$m\sum_{v\not=r}C(v)\le\sum_{v\not=r}m_vC(v)\le m\a$ for some sequence
$\{m_v\}_v$ with $m=\min_vm_v$, and let the reader divide by $m$.
In fact, in many instances we will derive $m_v=m$ for all $v\not=r$,
which makes for the following observation.

\begin{lem}{\sf [Uniform Covering Lemma]}\label{ucl} 
Let $\cT$ be a set of strategies for the root $r$ of the graph $G$.
If there is some $m$ such that, for each vertex $v\not=r$, we have
$\sum_{T\in\cT}T(v)=m$, then $\pi(G,r)=n(G)$.
\pf
\end{lem}

% ##########################################################################
%
%       SPECIFIC GRAPHS
%
\subsection{Specific Graphs}\label{Specific}

It has been said in jest that every graph theory paper should contain the
Petersen graph, so we get it out of the way first.

\begin{thm}\label{Pet}
Let $P$ denote the Petersen graph.
Then $\pi(P)\le 10$.
\end{thm}

Of course, the vertex lower bound implies $\pi(P)=10$, but since the focus
of this paper regards upper bounds, we prove them only.

\begin{figure}
\centerline{\includegraphics[height=1.8in]{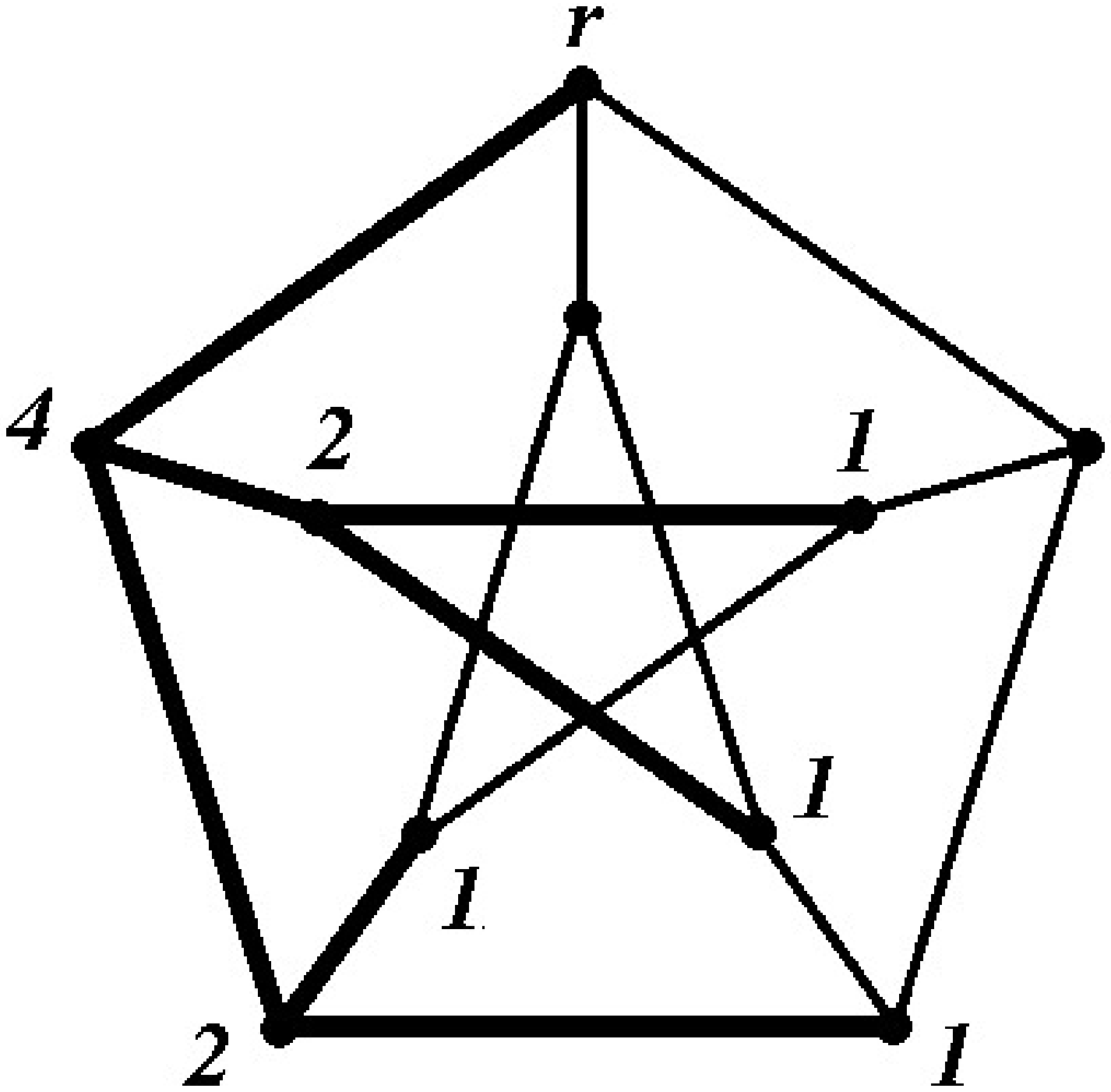}
\includegraphics[height=1.8in]{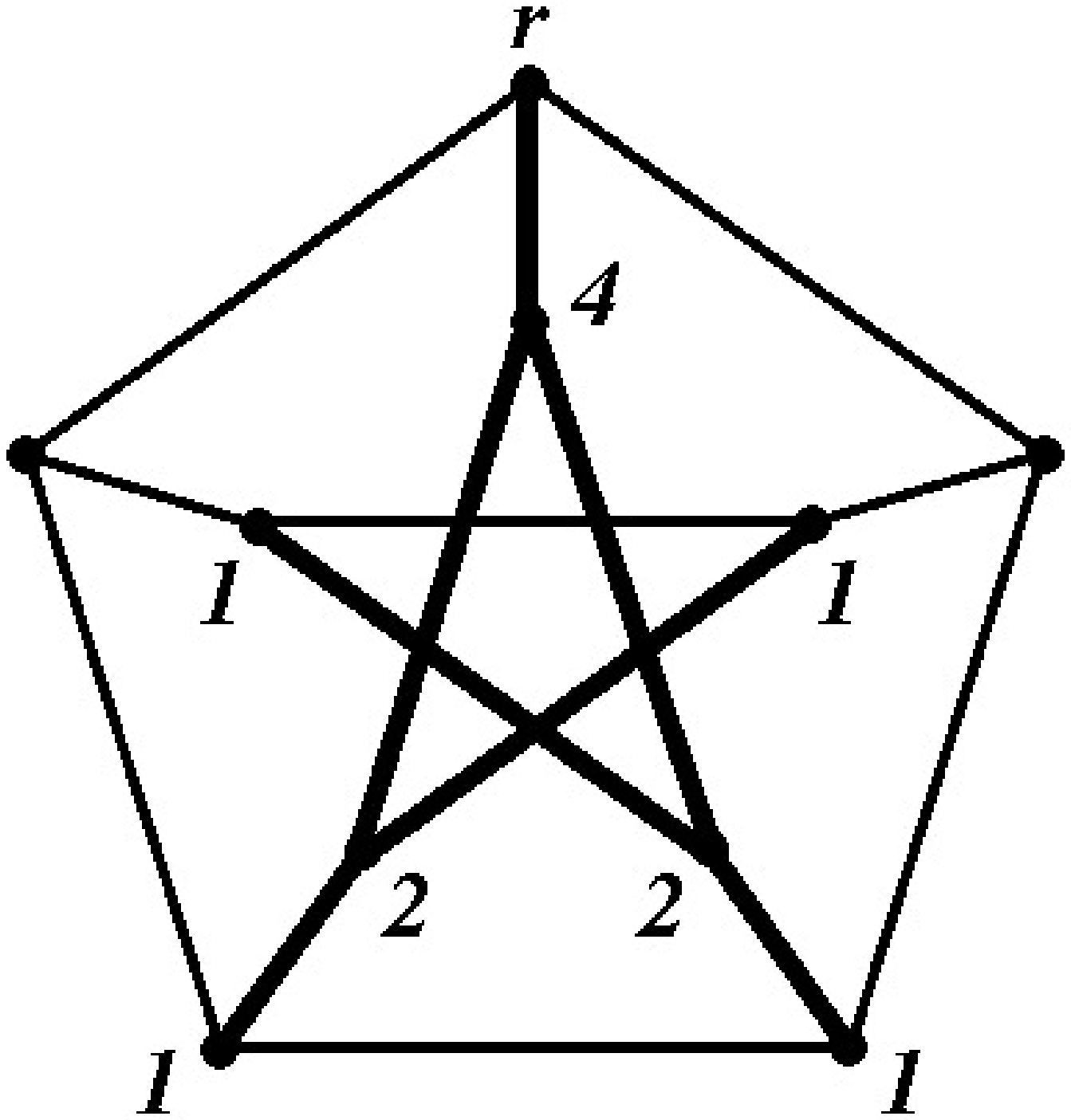}
\includegraphics[height=1.8in]{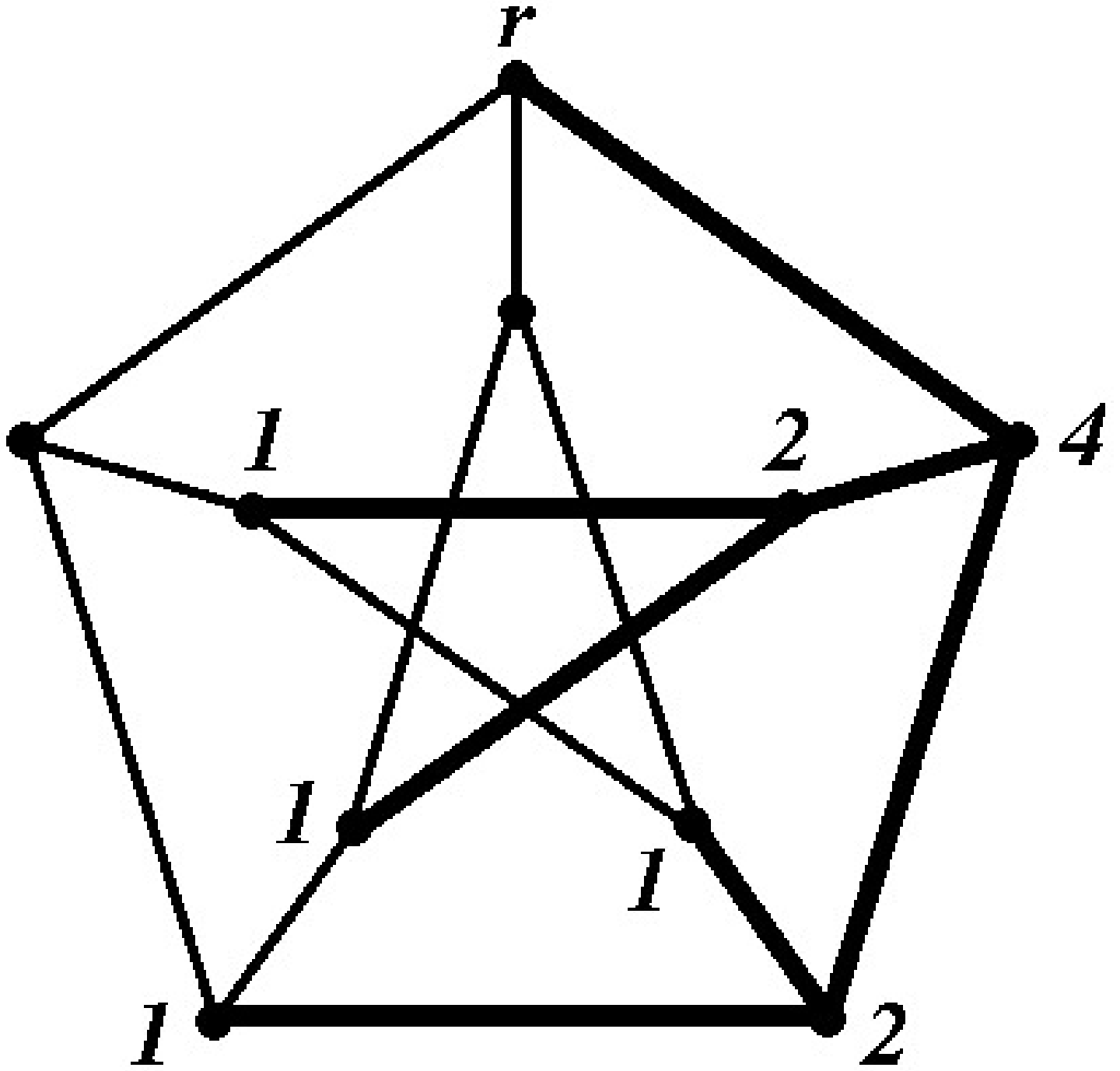}}
\caption{Petersen Class 0 certificate}
\label{PetStrat}
\end{figure}

\proof
The $3$ strategies shown in Figure \ref{PetStrat} certify the result.
\pf

Without such nice symmetry, the Lemke graph requires a different 
certificate for each possible root.

\begin{thm}\label{Lem}
Let $L$ denote the Lemke graph and suppose $r\not=v_1$.
Then $\pi(L,r)\le 8$.
\end{thm}

\proof
We show the strategies for each root vertex in turn, below.

For $r=v_2$:
\medskip

\centerline{\includegraphics[height=1.4in]{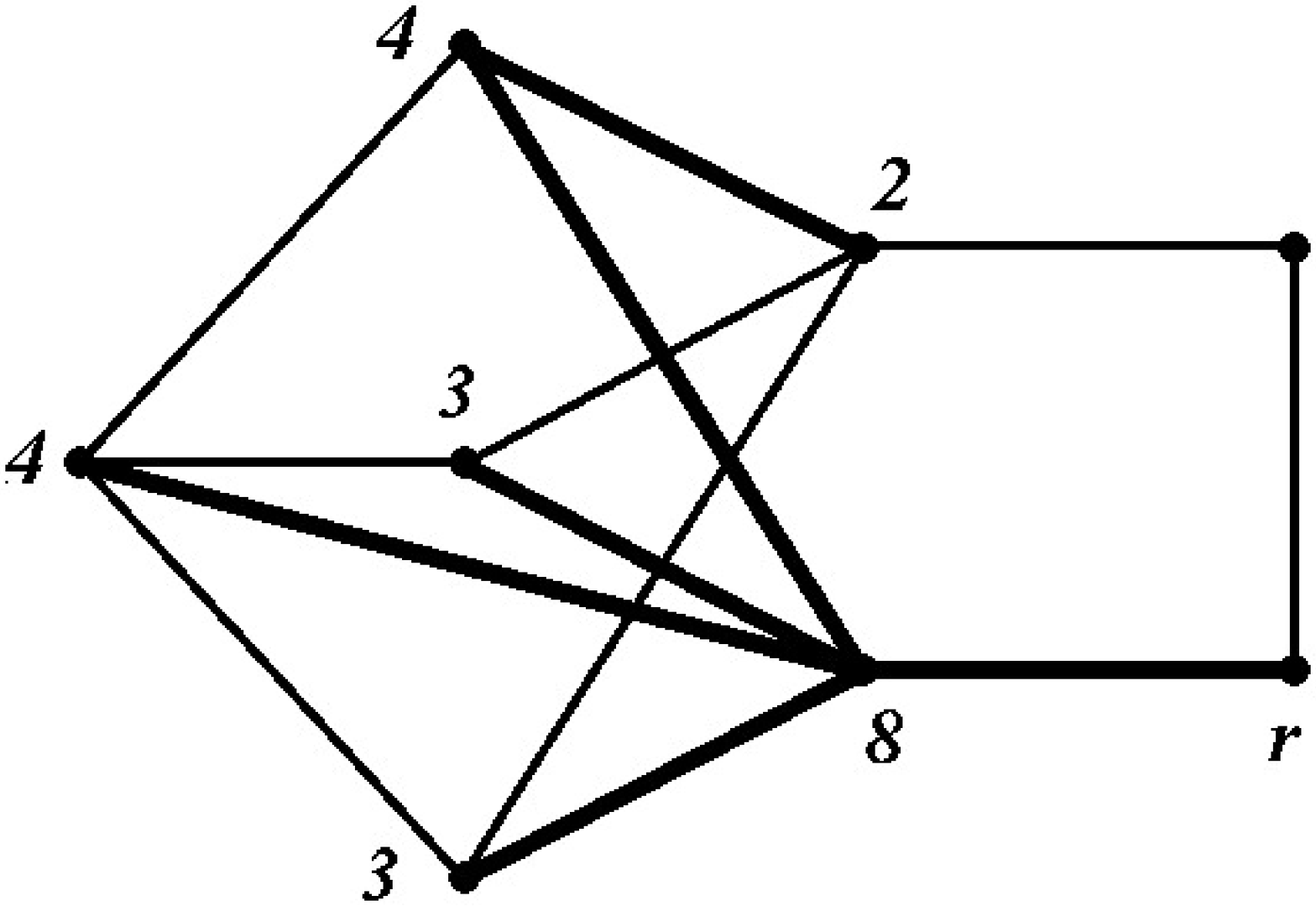}
\includegraphics[height=1.4in]{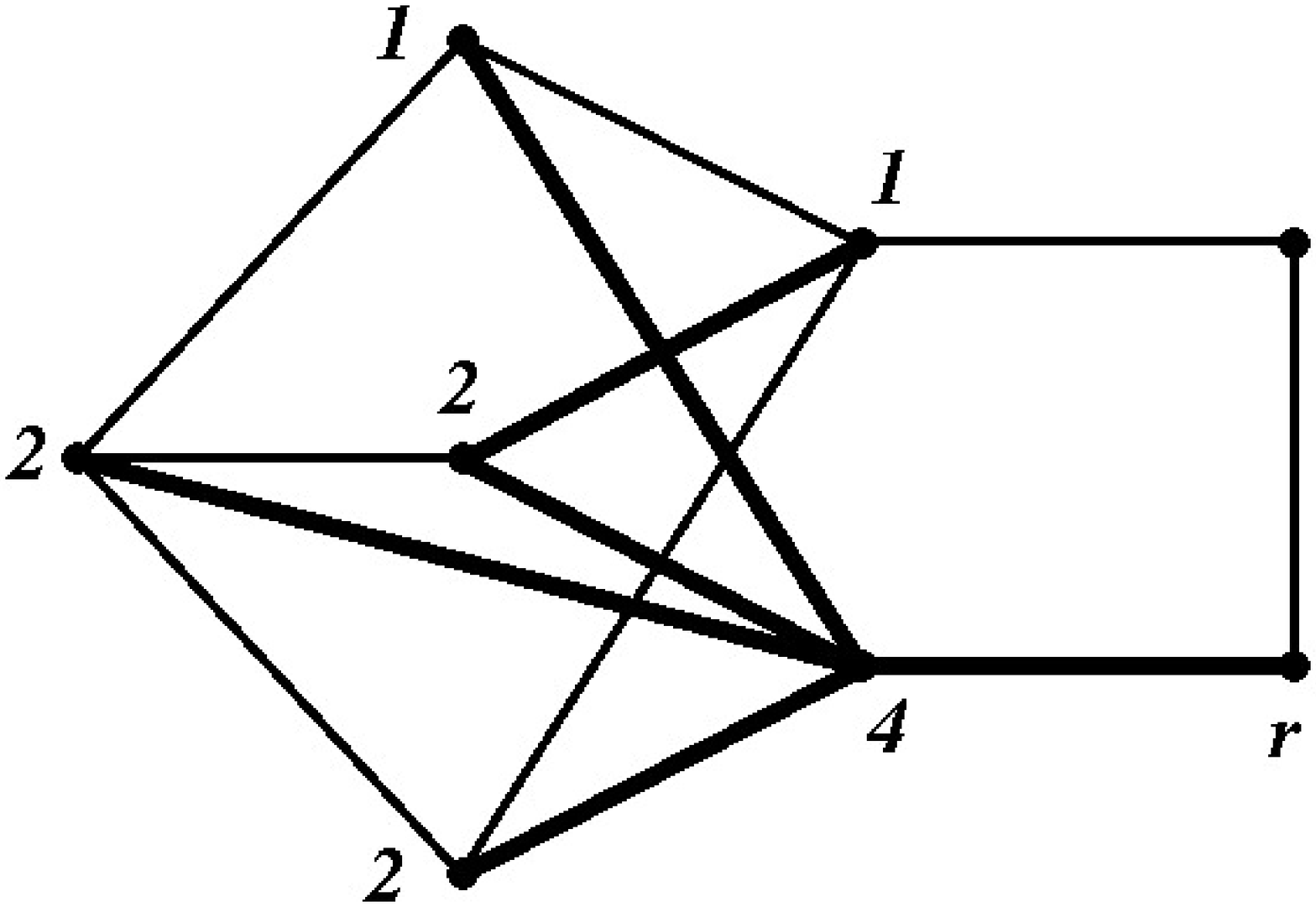}
\includegraphics[height=1.4in]{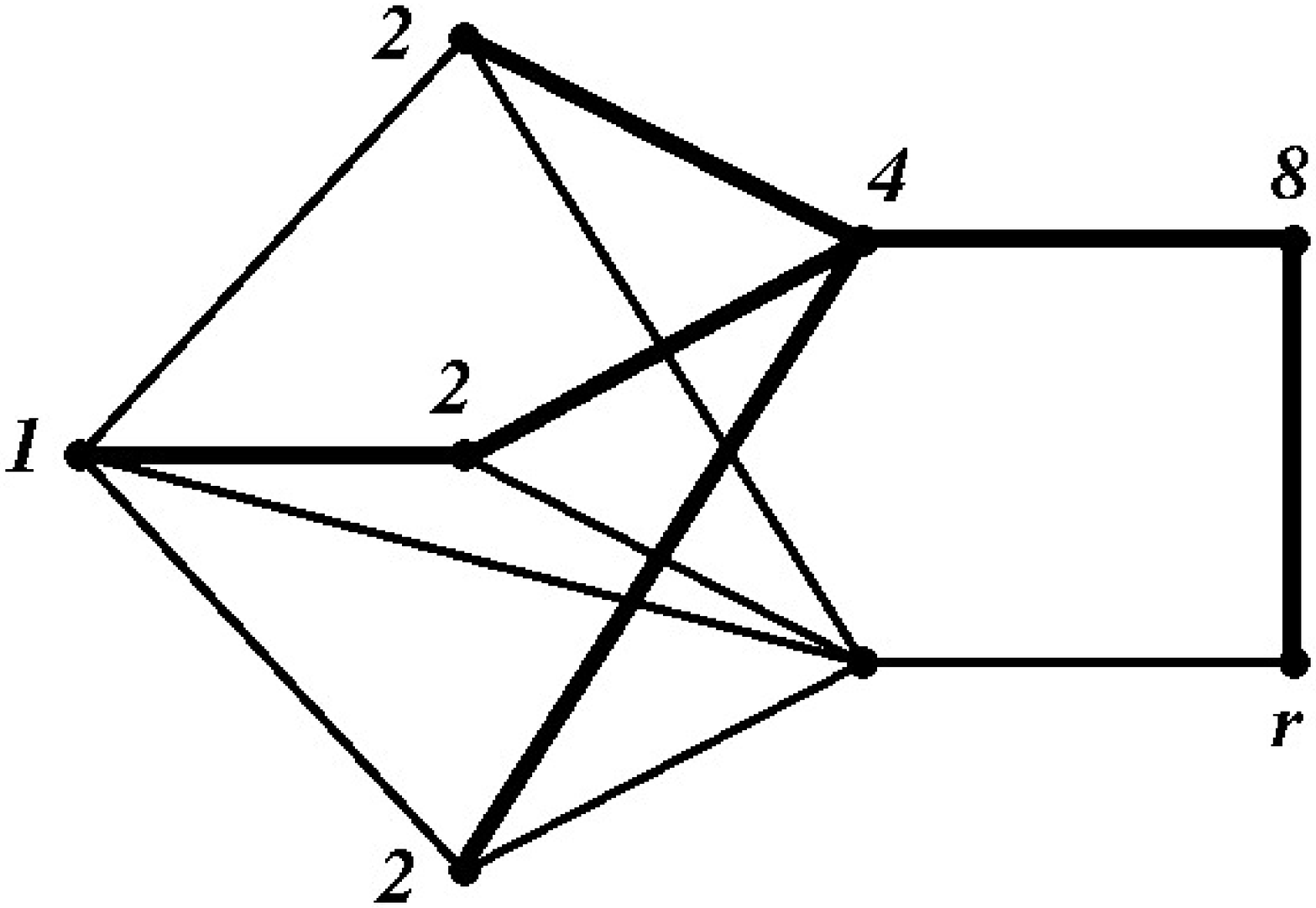}}
\medskip

For $r=v_3$:
\medskip

\centerline{\includegraphics[height=1.4in]{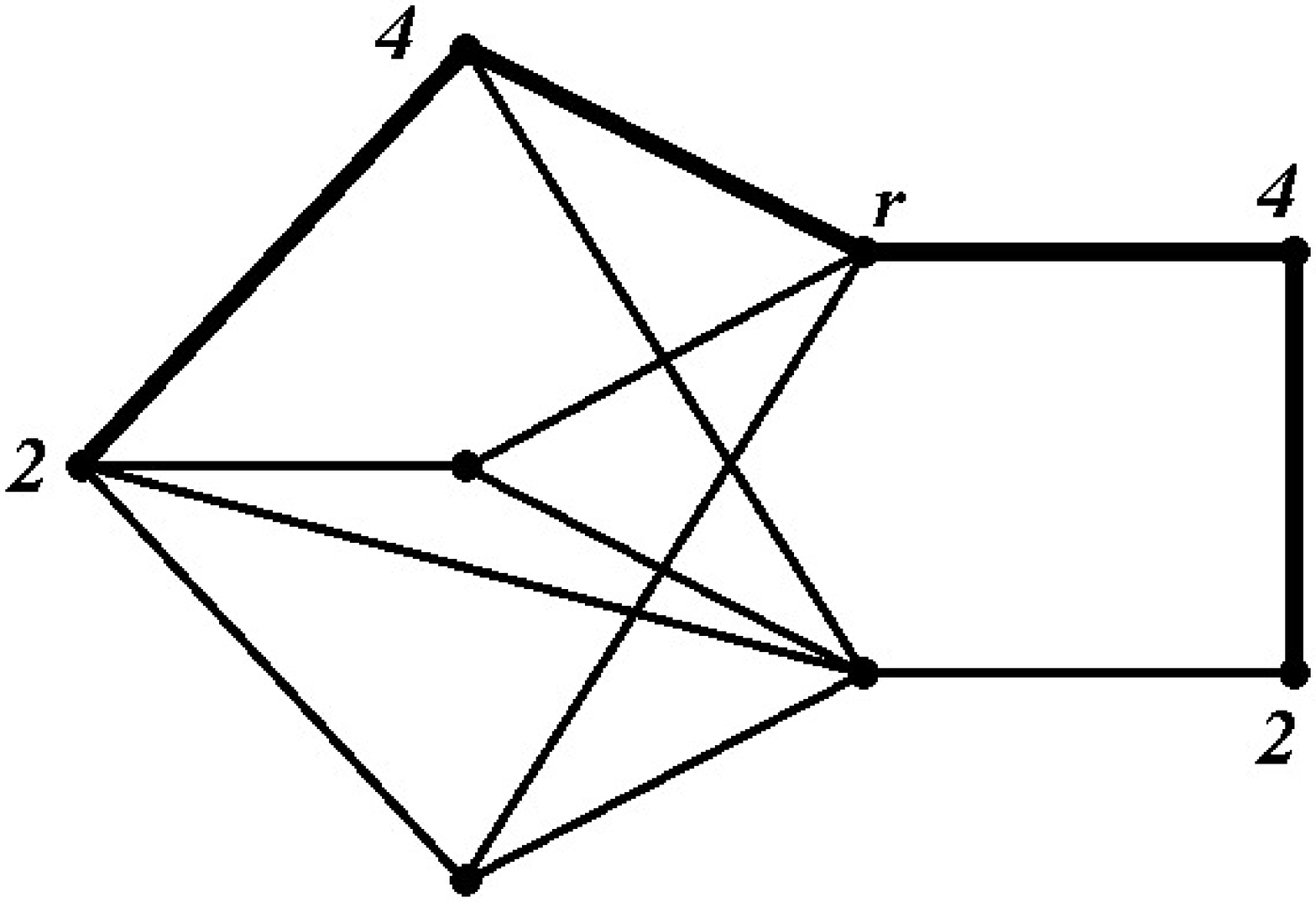}
\includegraphics[height=1.4in]{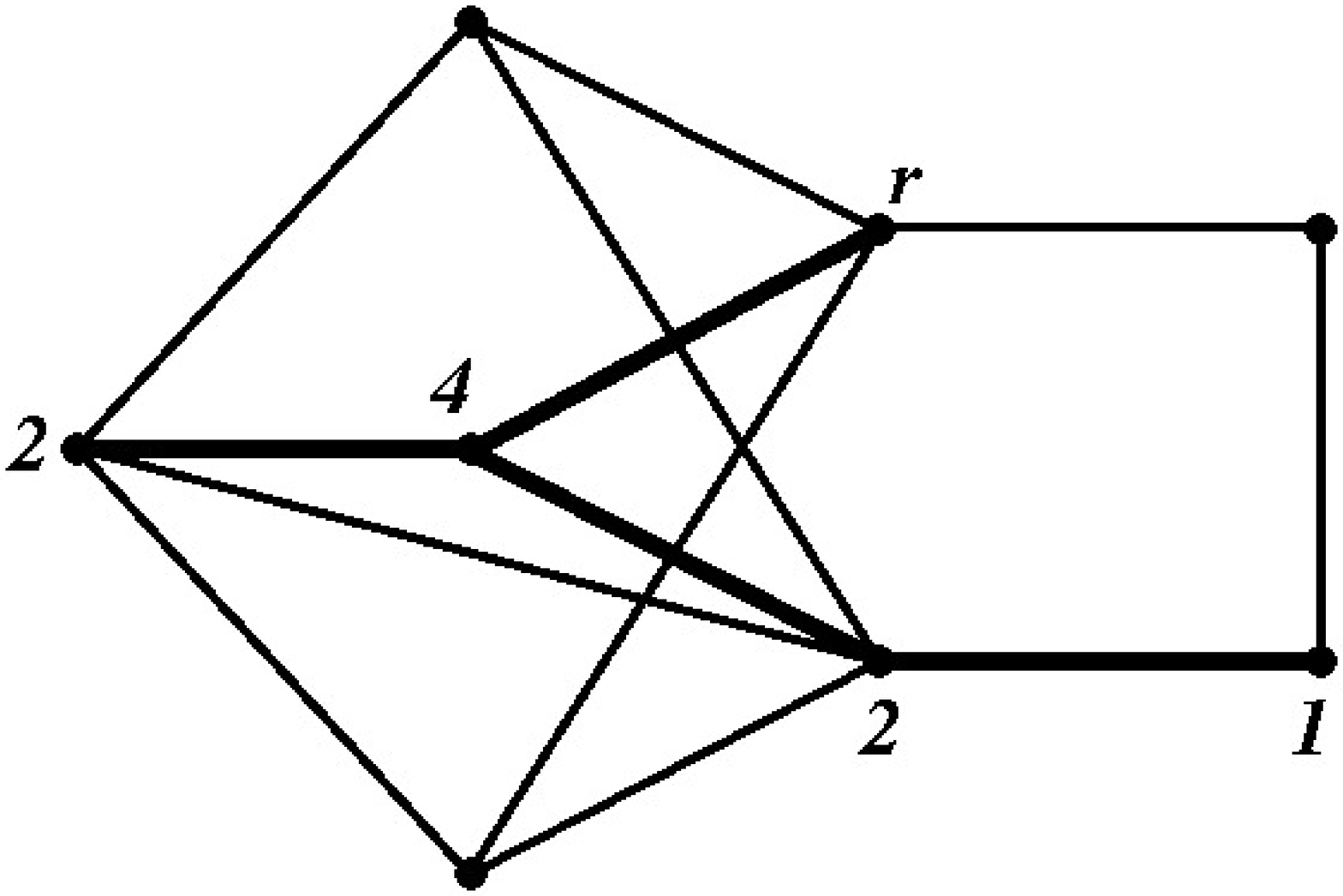}
\includegraphics[height=1.4in]{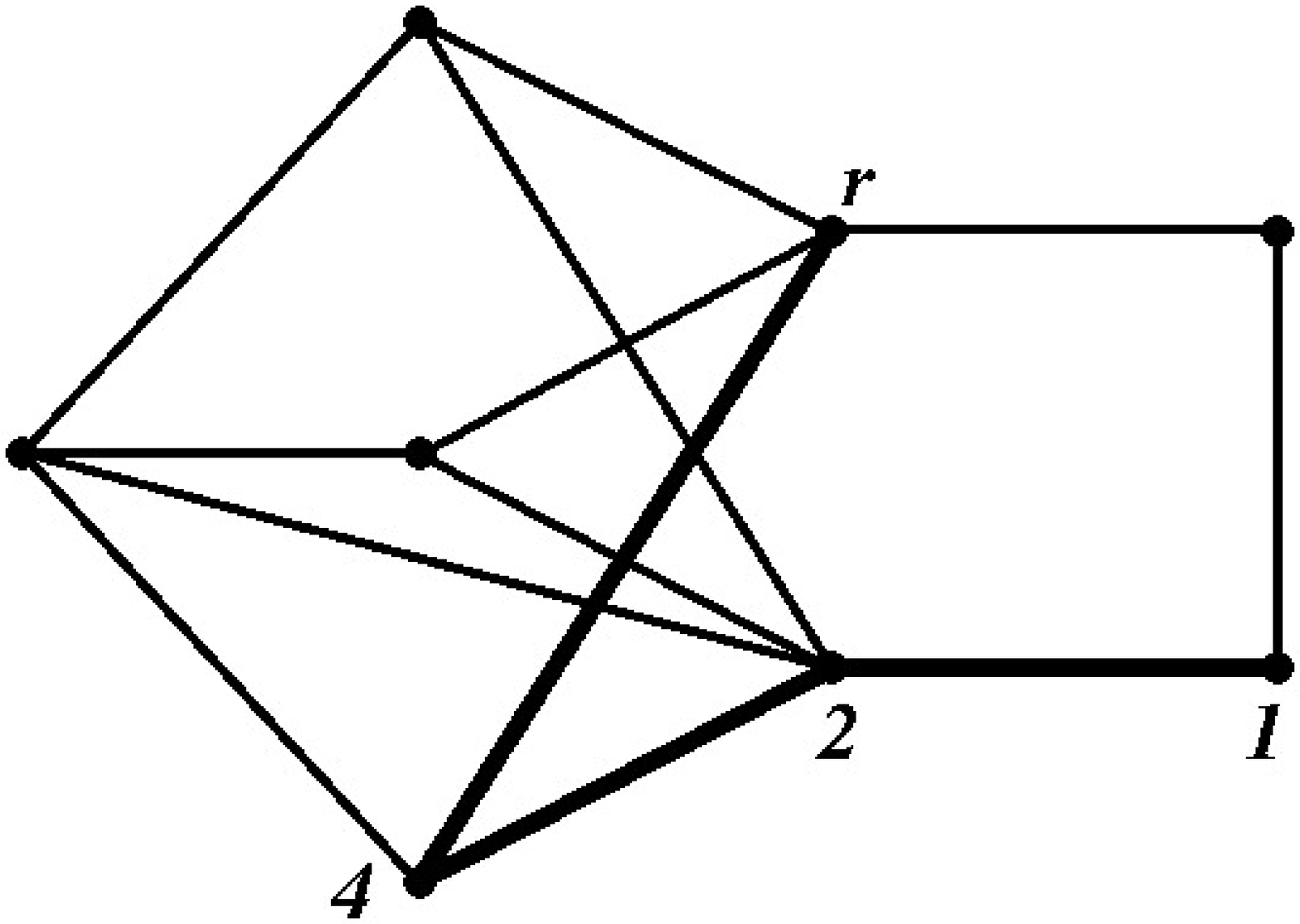}}
\medskip

Without the edge $v_4v_8$, the resulting graph would view $v_3$ and $v_4$
symmetrically.
Using that symmetry, the solutions for $r=v_3$ become solutions for $v_4$.

For $r=v_5$:
\medskip

\centerline{\includegraphics[height=1.4in]{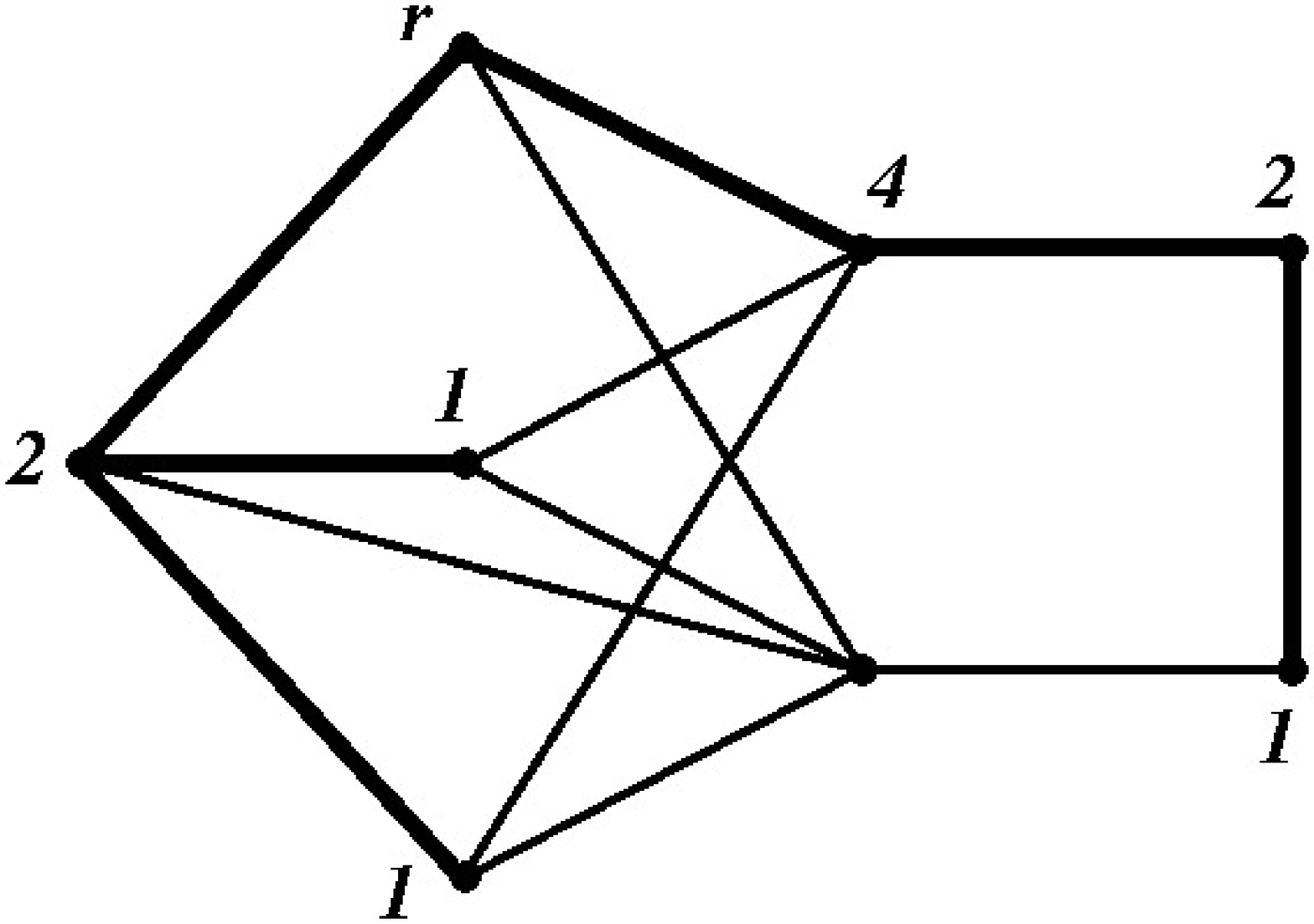}
$\qquad$
\includegraphics[height=1.4in]{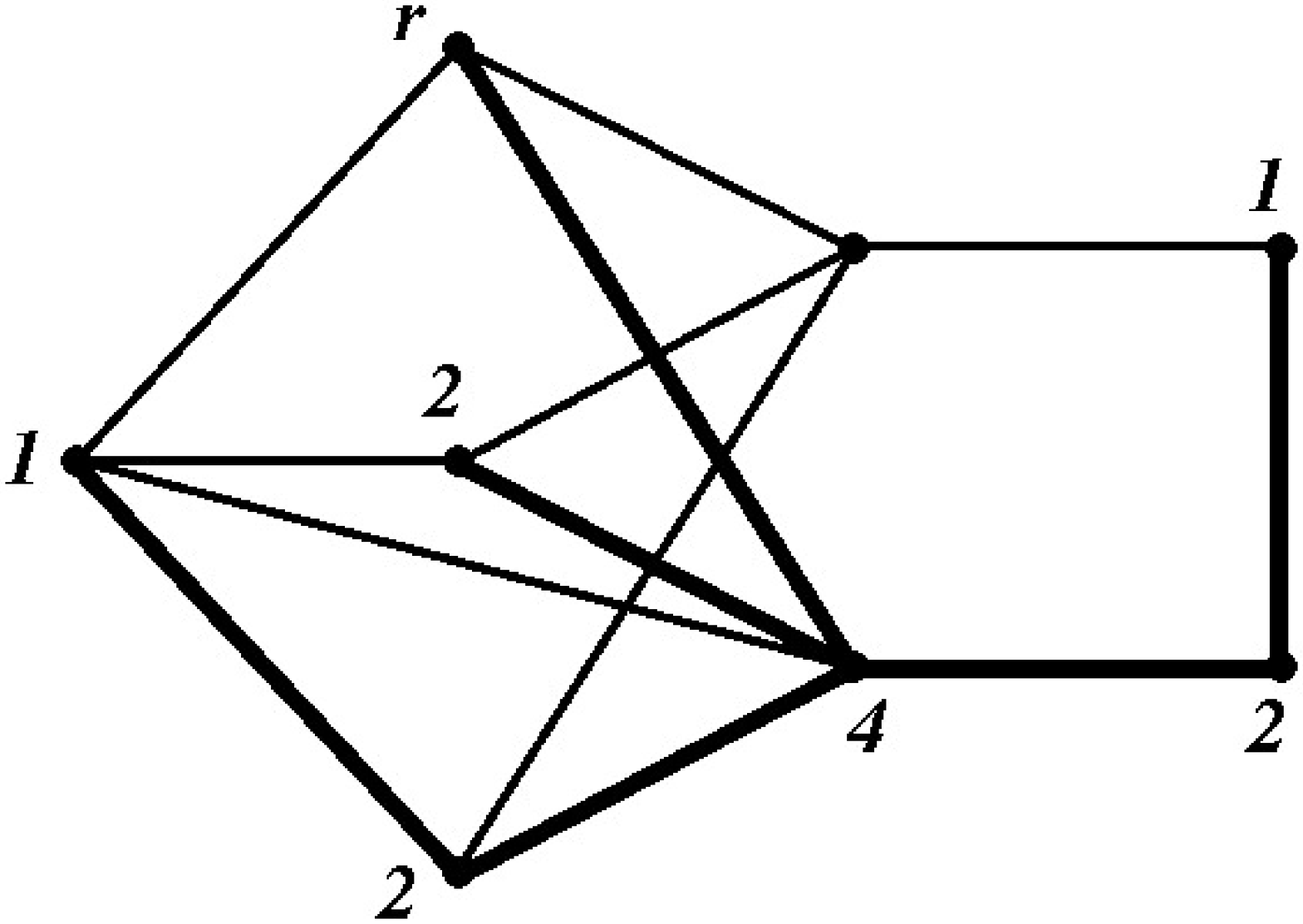}}
\medskip

For $r\in\{v_6,v_7\}$ we use the solutions from the case $r=v_5$ given by the
appropriate symmetry.
\medskip

For $r=v_8$:
\medskip

\centerline{\includegraphics[height=1.4in]{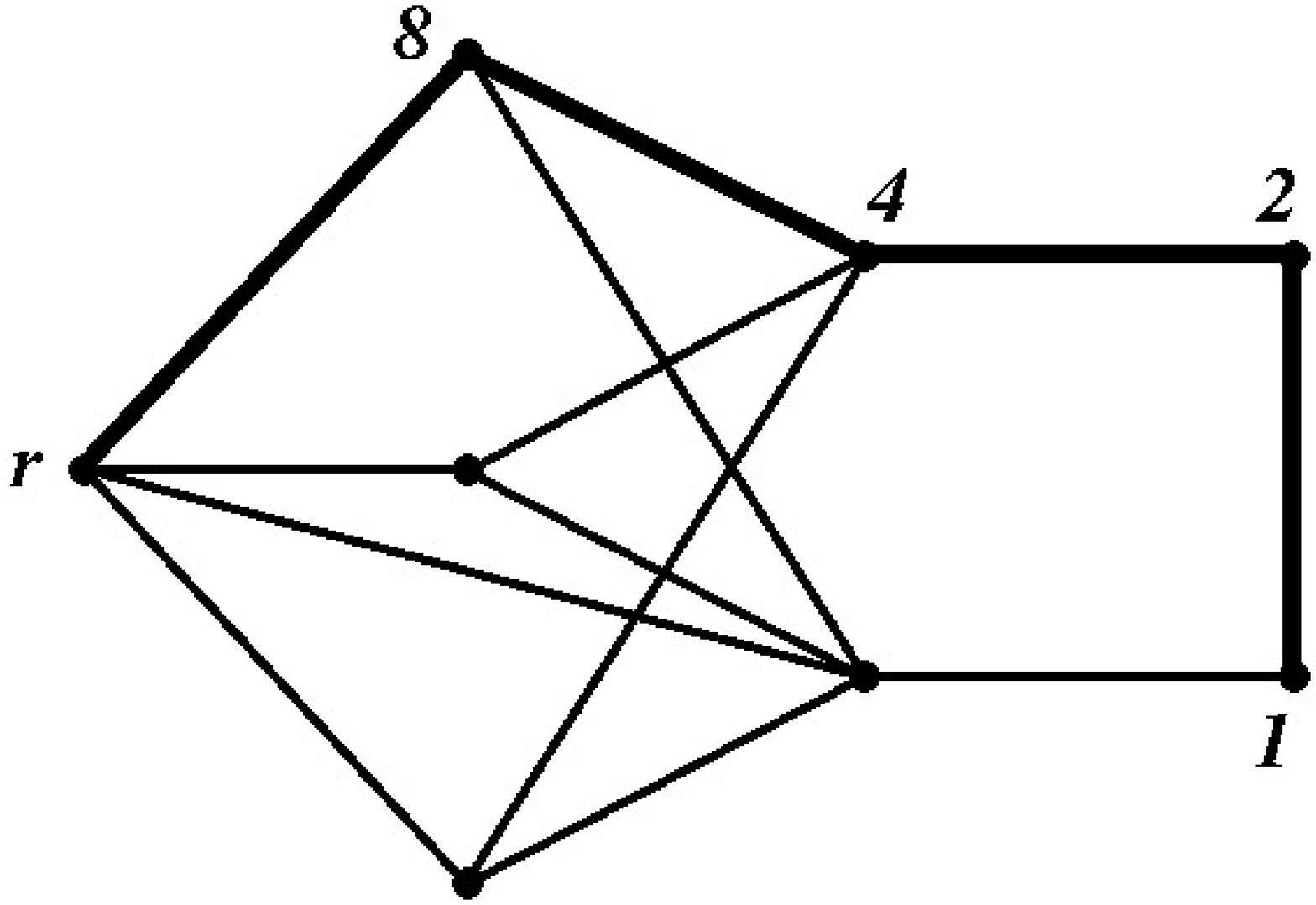}
$\qquad$
\includegraphics[height=1.4in]{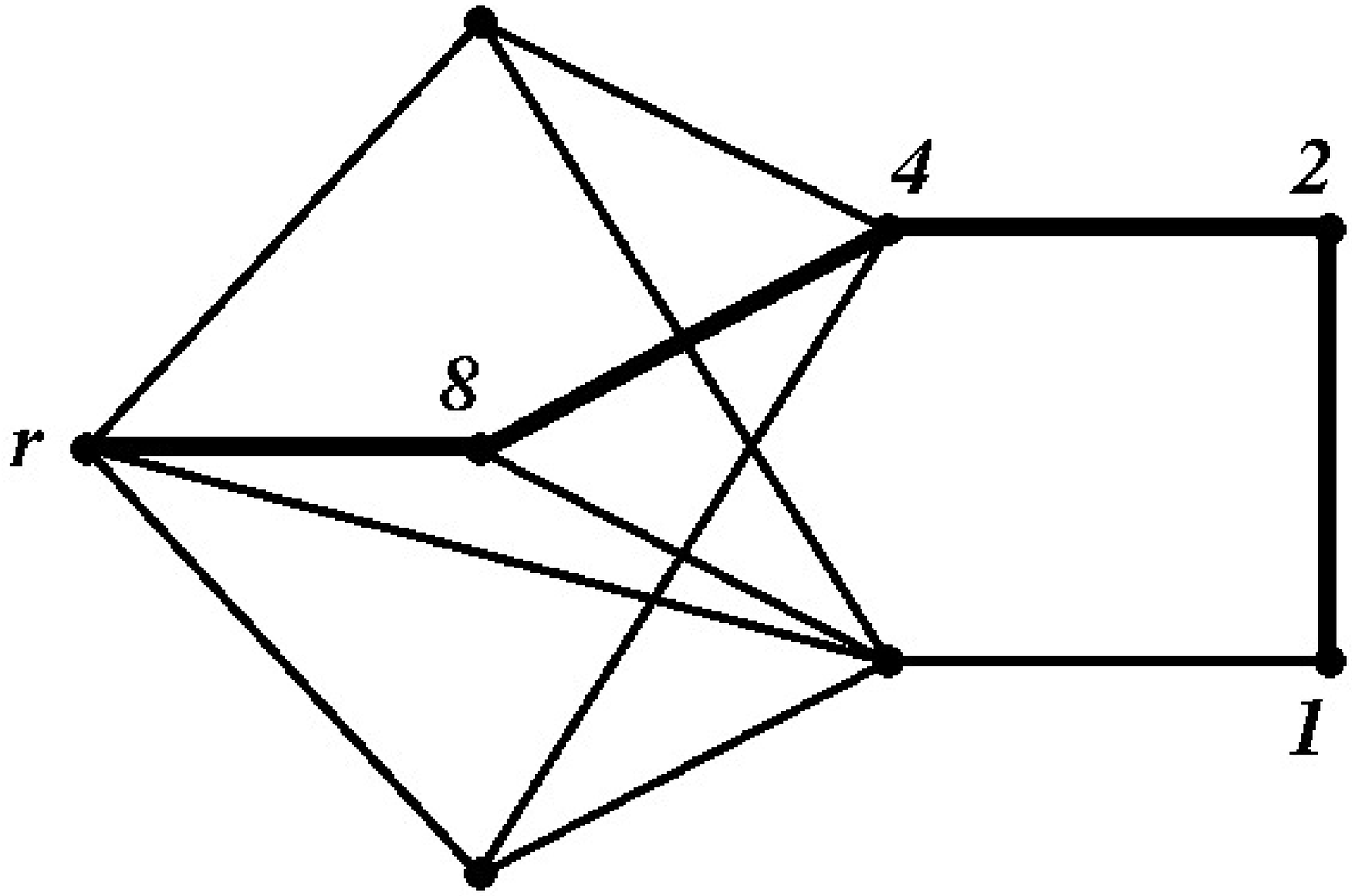}$\qquad\qquad\qquad$}

\centerline{$\qquad\qquad\qquad$\includegraphics[height=1.4in]{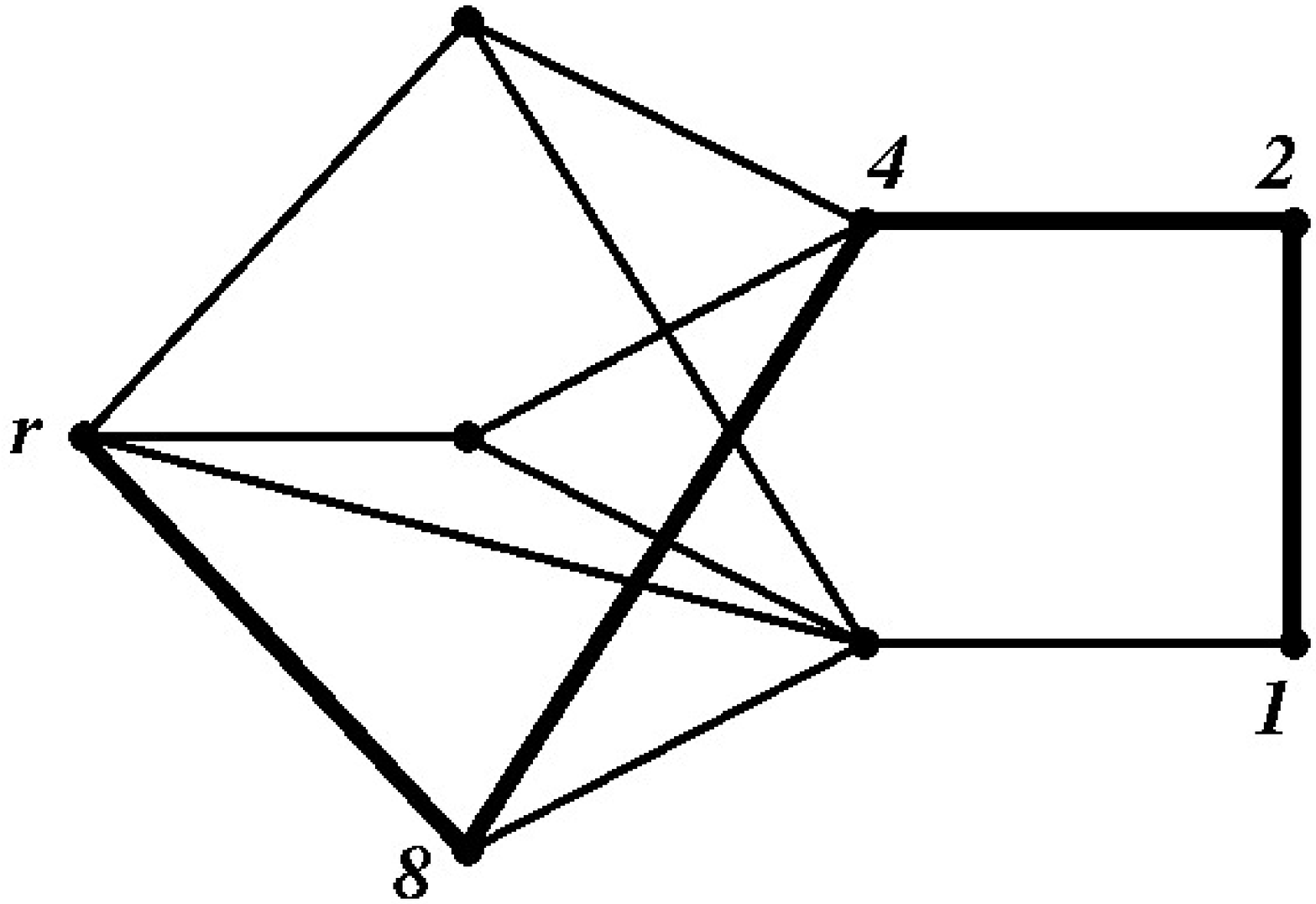}
$\qquad$
\includegraphics[height=1.4in]{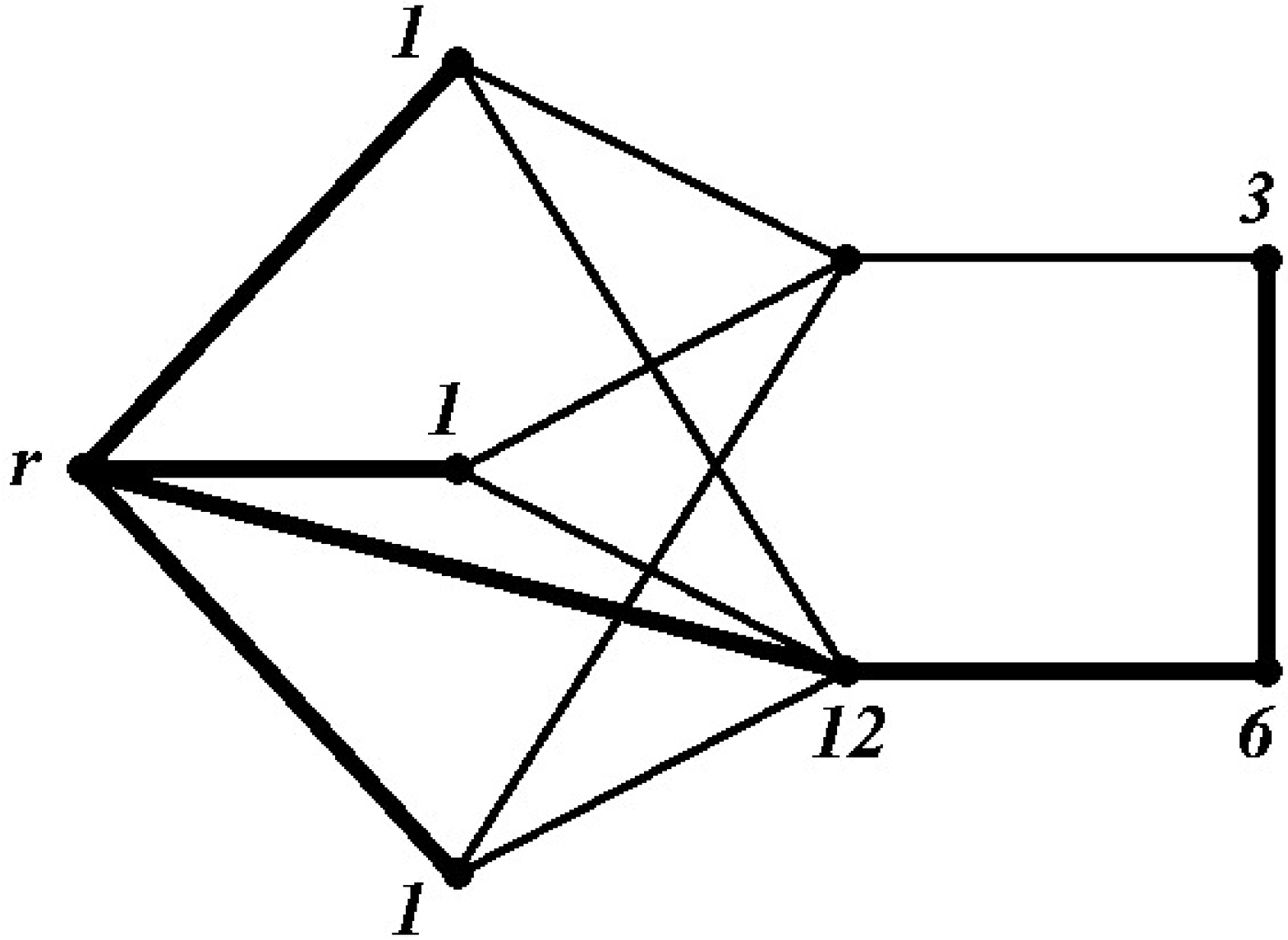}}

\pf

Because of the $3$-cube-like configuration with $(1,1,1,5)$ on 
$(v_5,v_6,v_7,v_8)$, the best that our tree strategies can muster
is $\pi(L,v_1)\le 9$.
Thus, to show that $L$ is Class 0 one must handle $v_1$ by more
traditional methods.

\begin{figure}
\centerline{\includegraphics[height=3.5in]{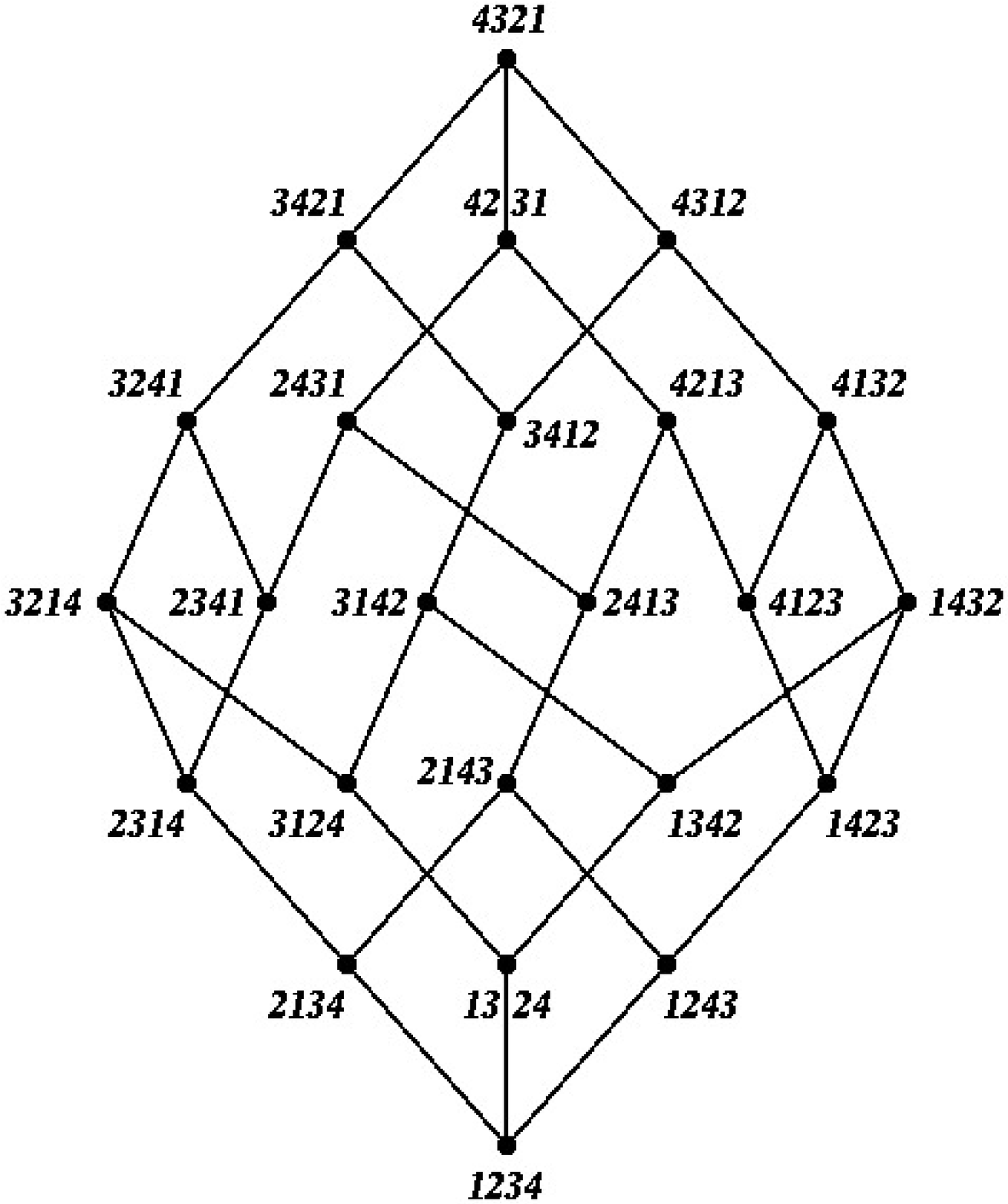}
$\qquad$
\includegraphics[height=3.5in]{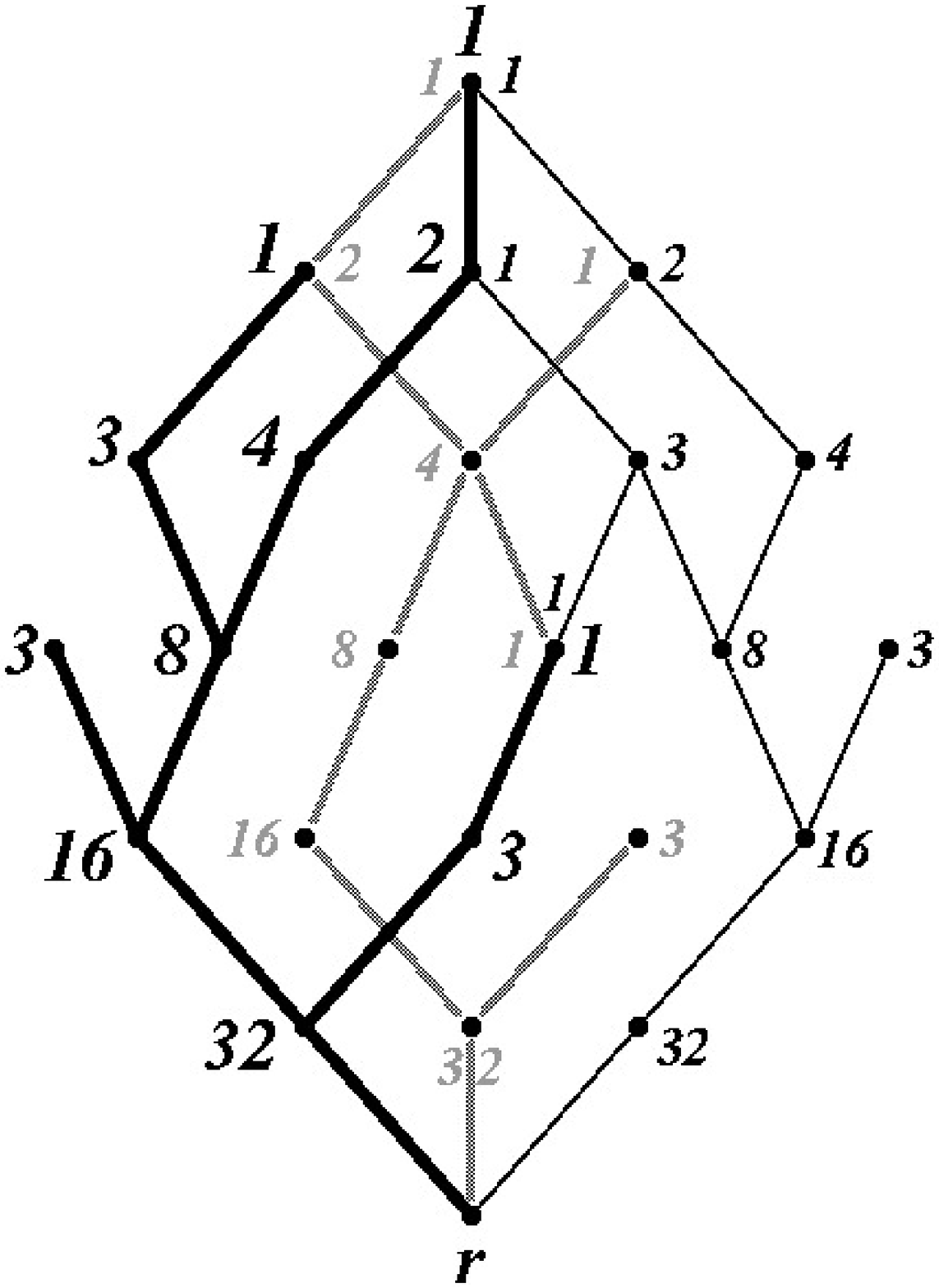}}
\caption{Bruhat graph $B_4$ (left) and strategies (right)}
\label{BruStrat}
\end{figure}

To illustrate that larger graphs can be tackled,
we move on to one of order 24 that is not Class 0.
Define the (weak) {\it Bruhat} graph of order $m$ (see Figure \ref{BruStrat})
to have all permutations of $[m]$ as vertices, with an edge between pairs of
permutations that differ by an adjacent transposition.
One can recognize it as the Cayley graph of $S_m$, generated by adjacent 
transpositions, and also note that $B_4$ is the cubic Ramanujan (expander) 
graph of \cite{Chiu}.
Intuitively, expander graphs would seem to have low pebbling numbers,
but because $B_4$ has diameter 6 we have the lower bound $\pi(B_4)\ge 64$.
We give here a fairly tight bound.

\begin{thm}\label{Bru}
Let $B_4$ be the Bruhat graph of order 4.
Then $\pi(B_4)\le 72$.
\end{thm}

\proof
Because the graph is vertex transitive, only one root must be checked.
The $3$ strategies shown in Figure \ref{BruStrat} certify the result.
We combined them into one figure, separated by edge styles.
\pf

Next we consider $\L2$, the square of the Lemke graph $L$.
Because $L$ has diameter 3, $\L2$ has diameter 6, and so
$\pi(\L2)\ge 64$.
Strategies deliver the following upper bounds.
Since the bounds are not that tight, we do not pursue bounds for
all $64$ vertices, although it is likely (since $v_1$ is the most
problematic root in $L$) that the upper bound of $108$
works for all roots $(r_1,r_2)$.

\begin{thm}\label{Lem2}
Let $\L2$ be the square of the Lemke graph $L$.
Then
$$\pi(\L2,(r,r))\le\Bigg\{
\begin{array}{rcl}
108&&{\rm if\ }r=v_1\ ,\\
96&&{\rm if\ }r=v_8{\rm\ ,\ and}\\
68&&{\rm if\ }r=v_4\ .
\end{array}$$
\end{thm}

\proof
For $r=v_1$ one can verify that a quarter of the sum of the four
strategies in Figure \ref{lemke2trees1} yields the weights 

$$\qquad\begin{array}{c|c|c|c|c|c|c|c|c|}
T_3&v_1&v_2&v_3&v_4&v_5&v_6&v_7&v_8\\
\hline
v_1&0&8&8&1&4&1&1&2\\
\hline
v_2&8&4&1&2&1&1&1&1\\
\hline
v_3&8&4&1&1&1&1&1&1\\
\hline
v_4&4&1&1&1&1&1&1&1\\
\hline
v_5&1&2&1&1&1&1&1&1\\
\hline
v_6&1&1&1&1&1&1&1&1\\
\hline
v_7&1&1&1&1&1&1&1&1\\
\hline
v_8&2&1&1&1&1&1&1&1\\
\hline
\end{array}\ ,$$
giving the bound of $63+4(7+3+1)+1=108$.
One verifies that each $T_i$ is a strategy by making sure that
each nonzero entry has a corresponding entry in its column or row
with at least twice its weight that is joined to it by an edge in
the appropriate copy of $L$.

\begin{figure}

{\footnotesize
$$\begin{array}{c|c|c|c|c|c|c|c|c|}
T_1&v_1&v_2&v_3&v_4&v_5&v_6&v_7&v_8\\
\hline
v_1&0&&32&&16&4&4&8\\
\hline
v_2&&&4&&&&&\\
\hline
v_3&&&4&&2&2&2&4\\
\hline
v_4&&&&&&&&\\
\hline
v_5&&&2&&1&1&1&2\\
\hline
v_6&&&2&&1&1&1&\\
\hline
v_7&&&2&&1&1&1&\\
\hline
v_8&&&&&&&&1\\
\hline
\end{array}
\qquad
\begin{array}{c|c|c|c|c|c|c|c|c|}
T_2&v_1&v_2&v_3&v_4&v_5&v_6&v_7&v_8\\
\hline
v_1&0&&&&&&&\\
\hline
v_2&32&&&&&&&\\
\hline
v_3&&&&&&&&\\
\hline
v_4&16&4&4&&2&2&2&2\\
\hline
v_5&4&&2&&1&1&1&\\
\hline
v_6&4&&2&&1&1&1&2\\
\hline
v_7&4&&2&&1&1&1&2\\
\hline
v_8&8&&4&&2&2&2&1\\
\hline
\end{array}
$$
\medskip

$$\begin{array}{c|c|c|c|c|c|c|c|c|}
T_3&v_1&v_2&v_3&v_4&v_5&v_6&v_7&v_8\\
\hline
v_1&0&&&&&&&\\
\hline
v_2&&&&&&&&\\
\hline
v_3&32&16&&4&2&2&2&\\
\hline
v_4&&&&&&&&\\
\hline
v_5&&8&&2&1&1&1&1\\
\hline
v_6&&4&&2&1&1&1&1\\
\hline
v_7&&4&&2&1&1&1&1\\
\hline
v_8&&4&&2&1&1&1&1\\
\hline
\end{array}
\qquad
\begin{array}{c|c|c|c|c|c|c|c|c|}
T_4&v_1&v_2&v_3&v_4&v_5&v_6&v_7&v_8\\
\hline
v_1&0&32&&4&&&&\\
\hline
v_2&&16&&8&4&4&4&4\\
\hline
v_3&&&&&&&&\\
\hline
v_4&&&&4&2&2&2&2\\
\hline
v_5&&&&2&1&1&1&1\\
\hline
v_6&&&&2&1&1&1&1\\
\hline
v_7&&&&2&1&1&1&1\\
\hline
v_8&&&&2&1&1&1&1\\
\hline
\end{array}$$
}
\smallskip

\caption{Strategies for $\L2$ at $(v_1,v_1)$.}\label{lemke2trees1}
\end{figure}

We use a similar argument when $r=v_8$, using the strategies from
Figure \ref{lemke2trees8}, along with their transposes, and divide by 8
to obtain the bound $63+8(3+1)+1=96$.

\begin{figure}

{\footnotesize
$$\begin{array}{c|c|c|c|c|c|c|c|c|}
T_1&v_1&v_2&v_3&v_4&v_5&v_6&v_7&v_8\\
\hline
v_1&1&2&2&4&4&4&4&8\\
\hline
v_2&&&&&&&&\\
\hline
v_3&&&&&&8&8&16\\
\hline
v_4&&&&&&&&\\
\hline
v_5&&&&&&&&\\
\hline
v_6&&&&&&&&\\
\hline
v_7&&&&8&8&8&4&32\\
\hline
v_8&&&&&&&&0\\
\hline
\end{array}
\qquad
\begin{array}{c|c|c|c|c|c|c|c|c|}
T_2&v_1&v_2&v_3&v_4&v_5&v_6&v_7&v_8\\
\hline
v_1&1&2&2&4&4&4&4&\\
\hline
v_2&2&4&4&8&8&8&8&16\\
\hline
v_3&&&&&&&&\\
\hline
v_4&&&&&&&&32\\
\hline
v_5&&&&&&&&\\
\hline
v_6&&&&&&&&\\
\hline
v_7&&&&&&&&\\
\hline
v_8&&&&&&&&0\\
\hline
\end{array}
$$
\smallskip

$$\begin{array}{c|c|c|c|c|c|c|c|c|}
T_3&v_1&v_2&v_3&v_4&v_5&v_6&v_7&v_8\\
\hline
v_1&1&2&&&&&&\\
\hline
v_2&&&&&&&&\\
\hline
v_3&&4&&8&&&&\\
\hline
v_4&&&&4&&&&\\
\hline
v_5&&&&16&4&&8&32\\
\hline
v_6&&&&&&&&\\
\hline
v_7&&&&&&&&\\
\hline
v_8&&&&&&&&0\\
\hline
\end{array}
\qquad
\begin{array}{c|c|c|c|c|c|c|c|c|}
T_4&v_1&v_2&v_3&v_4&v_5&v_6&v_7&v_8\\
\hline
v_1&1&&2&&&&&\\
\hline
v_2&&&&&&&&\\
\hline
v_3&2&&4&&8&&&\\
\hline
v_4&&&&&&&&\\
\hline
v_5&&&&&&&&\\
\hline
v_6&&&&8&16&4&&32\\
\hline
v_7&&&&&&&&\\
\hline
v_8&&&&&&&&0\\
\hline
\end{array}
$$
}
\smallskip

\caption{Strategies for $\L2$ at $(v_8,v_8)$.}\label{lemke2trees8}
\end{figure}

Likewise, for $r=v_4$, the strategies and their transposes from Figure
\ref{lemke2trees4} yield an upper bound of $63+\lfloor 8(.6)\rfloor+1=68$.
\pf

\begin{figure}

$$\begin{array}{c|c|c|c|c|c|c|c|c|}
T_1&v_4&v_8&v_7&v_6&v_5&v_2&v_3&v_1\\
\hline
v_4&0&&&&&&&\\
\hline
v_8&10&5&4&2&4&4&2&\\
\hline
v_7&&&&&&&&\\
\hline
v_6&&&&&&&&\\
\hline
v_5&&&&&&&&\\
\hline
v_2&&&&&&&&\\
\hline
v_3&&&&&&&&\\
\hline
v_1&&&&&&&&\\
\hline
\end{array}
\qquad
\begin{array}{c|c|c|c|c|c|c|c|c|}
T_2&v_4&v_8&v_7&v_6&v_5&v_2&v_3&v_1\\
\hline
v_4&0&&&&&&&\\
\hline
v_8&&&&&&&&\\
\hline
v_7&16&6&5&2&8&8&2&4\\
\hline
v_6&&&&&&&&\\
\hline
v_5&&&&&&&&\\
\hline
v_2&&&&&&&&\\
\hline
v_3&8&4&4&4&4&4&2&1\\
\hline
v_1&2&2&&2&2&2&1&1\\
\hline
\end{array}
$$
\smallskip

$$\begin{array}{c|c|c|c|c|c|c|c|c|}
T_3&v_4&v_8&v_7&v_6&v_5&v_2&v_3&v_1\\
\hline
v_4&0&&&&&&&\\
\hline
v_8&&&&&&&&\\
\hline
v_7&&&&&&&&\\
\hline
v_6&16&8&8&5&8&8&2&4\\
\hline
v_5&&&&&&&&\\
\hline
v_2&&&&&&&&\\
\hline
v_3&2&4&4&2&4&4&2&2\\
\hline
v_1&&2&2&&2&2&1&1\\
\hline
\end{array}
\qquad
\begin{array}{c|c|c|c|c|c|c|c|c|}
T_4&v_4&v_8&v_7&v_6&v_5&v_2&v_3&v_1\\
\hline
v_4&0&&&&&&&\\
\hline
v_8&&&&&&&&\\
\hline
v_7&&&&&&&&\\
\hline
v_6&&&&&&&&\\
\hline
v_5&16&6&2&2&5&8&2&4\\
\hline
v_2&&&&&&&&\\
\hline
v_3&&&&&&2&1&2\\
\hline
v_1&&&&&&&&1\\
\hline
\end{array}
$$
\smallskip

$$\begin{array}{c|c|c|c|c|c|c|c|c|}
T_5&v_4&v_8&v_7&v_6&v_5&v_2&v_3&v_1\\
\hline
v_4&0&&&&&&&\\
\hline
v_8&&&&&&&&\\
\hline
v_7&&&&&&&&\\
\hline
v_6&&&&&&&&\\
\hline
v_5&&&&&&&&\\
\hline
v_2&16&6&2&2&2&5&&2\\
\hline
v_3&&&&2&&&&\\
\hline
v_1&8&4&4&4&2&4&2&2\\
\hline
\end{array}
$$
\smallskip

\caption{Strategies for $\L2$ at $(v_4,v_4)$.}\label{lemke2trees4}
\end{figure}

One can see where the improvement from 108 to 96 comes from.
We obtain a bound from a sum of strategies by dividing by the minimum
weight, and since the graph has diameter six, the maximum weight is
at least 32.
So we need to increase the minimum weight as much as possible by 
including more strategies, but there are diminishing returns.
In the case of root $(v_1,v_1)$, it has four neighbors, and so a
fifth strategy increases the maximum weight to at least 64.
This is why the root $(v_8,v_8)$, having eight neighbors, fares better.
On this basis one might expect the root $(v_4,v_4)$ to perform best.
In fact, it does even better for a different reason: it has eccentricity $4$.
It is amusing that the order of this graph is out of the range of
computing, while these strategies were found fairly easily by hand.

% ##########################################################################
%
%       GRAPH CLASSES
%
\subsection{Random Graphs}\label{Random}

Given our introductory comments and observations, we also tested some
modest sized random graphs.
Let $R_{15}$ denote the $15$-vertex graph given by the adjacency list
\begin{quote}
$([2,4,5,6,12,13], [1,3,4,8,11,12,14], [2,4,6,7], [1,2,3,5,7,9,14],
[1,4,6,8,11,15],$ $[1,3,5,9,13,14], [3,4,11,15], [2,5,10,13,14,15],
[4,6,10,11], [8,9,11],$\\ $[2,5,7,9,10,12,15], [1,2,11,13], [1,6,8,12],
[2,4,6,8], [5,7,8,11])$.
\end{quote}
We generated this graph with edge probability $.35$.
It has diameter three, with $10$ of the vertices (including $v_9$)
having distance two from all others.
We found that $R_{15}$ is at least nearly Class 0 as follows.
Our simple certificate of 3 basic strategies for the root $r=v_9$ is below.
\begin{eqnarray*}
2x_1+2x_2+2x_3+4x_4+2x_5+2x_7+x_8+x_{12}+x_{13}+2x_{14}+x_{15}&\le&20\\
2x_2+2x_7+2x_8+4x_{10}+4x_{11}+2x_{12}+x_{13}+2x_{15}&\le&19\\
2x_1+2x_3+2x_5+4x_6+x_8+x_{12}+2x_{13}+2x_{14}+x_{15}&\le&17
\end{eqnarray*}
We note that there are $20,422$ basic strategies of depth at most two rooted
at $r$.
These were generated by our java code in a few minutes on my
Core 2 Duo 2.66GhZ PC, with 2GB RAM and 250GB HD running Ubuntu linux 9.04.
CPLEX software solved the resulting linear optimization problem
instantly, delivering a dual certificate of a conic combination of $11$ 
of the strategies.
By hand, we grouped them according to their neighbor of $r$ and were able
to combine each into a single nonbasic strategy.
By trading weights between some of the strategies, we arrived at the
simplified certificate above.
The certificate might be more easily viewed in the table of coefficients,
below, in which the sum of the constraints appears below the line.
Thus $\pi(R_{15},v_9)\le 15$.

$$\left(
\begin{array}{rrrrrrrrrrrrrrr|r}
2&2&2&4&2&0&2&1&0&0&0&1&1&2&1&20\\
0&2&0&0&0&0&2&2&0&4&4&2&1&0&2&19\\
2&0&2&0&2&4&0&1&0&0&0&1&2&2&1&17\\
\hline
4&4&4&4&4&4&4&4&0&4&4&4&4&4&4&56\\
\end{array}\right)$$

The same upper bound of 15 was obtained for all other roots
(see the Appendix\footnote{We should point out that whenever there are
too many possible strategies to compute in a reasonable amount of time 
and space --- on random graphs with 30 vertices and diameter 4 it was 
typical to crash memory after a day of running --- one can generate a 
large collection of strategies at random rather than by exhaustion and
obtain optimal results far more quickly.  This was a useful speed up 
even for $R_{20}$ with roots $v_3$, $v_{10}$, $v_{12}$, $v_{17}$, 
and $v_{20}$.})
except for $v_{10}$, for which we found $\pi(R_{15},v_{10})\le 16$.
Note that $R_{15}$ contains the $4$-cycle $(v_4,v_3,v_6,v_9)$, with $v_9$
adjacent to the root $v_{10}$ --- the same dreaded configuration we
discussed in $Q^3$ obstructs us here: $(1,5,1,0)$.
More to the point, let $C$ be the size $15$ configuration 
$(1,$ $1,$ $5,$ $1,$ $1,$ $1,$ $1,$ $0,$ $0,$ $0,$ $0,$ $1,$ $1,$ $1,$ $1)$
on $R_{15}$.
Because $v_3$ is distance $3$ from the root $v_{10}$ and $C$ is
empty on the neighbors of $v_{10}$, the only way to move a pebble
to the root involves splitting the pebbles from $v_3$ as in the
dreaded configuration solution.
Hence no strategy can detect the $v_{10}$-solvability of $C$.
However, we can blend strategies with some case analysis and slight
amount of old fashioned analysis as follows.

\begin{thm}\label{R15class0}
The graph $R_{15}$ is Class 0.
\end{thm}

\proof
As we have already shown that $\pi(R_{15},r)=15$ for all $r$ except
for $v_{10}$, it suffices to prove that $\pi(R_{15},v_{10})$.
Let $C$ be an $r$-unsolvable configuration, where $r=v_{10}$.
We first consider the case that $C$ contains a pebble on a neighbor of $r$.
In this case we have $x_8+x_9+x_{11}\ge 1$, and the CPLEX certificate

{\tiny
$$\left(
\begin{array}{r|rrrrrrrrrrrrrrr|r}
333&0&0&0&0&0&0&0&-1&-1&0&-1&0&0&0&0&-1\\
80&2&0&2&3&0&4&0&0&8&0&0&0&2&2&0&23\\
37&2&0&2&4&2&4&0&0&8&0&0&0&1&2&0&25\\
50&2&2&2&4&2&3&2&0&8&0&0&0&0&2&0&27\\
3&2&4&2&2&0&1&0&8&0&0&0&2&4&4&0&29\\
1&2&4&2&2&0&2&0&8&0&0&0&2&4&0&3&29\\
136&2&1&2&2&4&0&4&0&0&0&8&4&2&0&2&31\\
25&2&3&2&2&4&2&4&0&0&0&8&4&0&0&2&33\\
4&2&4&2&2&3&0&0&0&0&0&8&4&2&0&4&31\\
37&2&4&2&2&1&2&0&8&0&0&0&2&0&4&4&31\\
127&2&4&2&0&1&2&2&8&0&0&0&2&4&4&4&35\\
1&2&4&2&2&4&2&1&0&0&0&8&4&2&2&4&37\\
\hline
&1002&1003&1002&1002&998&1003&999&1011&1003&0&995&1000&1003&1004&1001&14692\\
\end{array}\right)$$
}%
shows that $|C|\le 14$.
In this format, the left column holds the multipliers of each constraint, 
and the bottom row is the result of the linear combination of them.
Division by $995$ yields the result.  

Thus we may assume that $x_8+x_9+x_{11}=0$.
Next we consider the case that at most $5$ pebbles are distance $3$ from $r$;
that is, $x_1+x_3\le 5$.
Here we have the certificate

{\scriptsize
$$\left(
\begin{array}{r|rrrrrrrrrrrrrrr|r}
84&1&0&1&0&0&0&0&0&0&0&0&0&0&0&0&5\\
168&0&0&0&0&0&0&0&1&1&0&1&0&0&0&0&0\\
6&2&4&2&2&3&2&0&8&0&0&0&2&4&4&0&33\\
13&2&0&2&4&0&4&1&0&8&0&0&0&2&2&0&25\\
29&2&0&2&4&2&4&2&0&8&0&0&0&1&2&0&27\\
55&2&4&2&2&0&2&2&8&0&0&0&2&3&4&4&35\\
23&2&4&2&2&4&2&1&8&0&0&0&2&4&4&4&39\\
42&2&0&2&2&4&2&4&0&0&0&8&4&2&0&1&31\\
6&2&4&2&2&4&2&4&0&0&0&8&4&2&2&1&37\\
36&2&4&2&2&4&2&3&0&0&0&8&4&2&2&4&39\\
\hline
&504&504&504&504&504&504&504&840&504&0&840&504&504&504&504&7476\\
\end{array}\right)\ .$$
}%
On the other hand, at least $7$ pebbles at distance $3$ from $r$ yields
the following certificate.

{\footnotesize
$$\left(
\begin{array}{r|rrrrrrrrrrrrrrr|r}
20&-1&0&-1&0&0&0&0&0&0&0&0&0&0&0&0&-7\\
1&2&4&2&0&0&0&0&8&0&0&0&0&0&4&3&23\\
5&2&0&2&4&1&0&0&0&8&0&0&0&0&0&0&17\\
3&2&0&2&0&0&4&0&0&8&0&0&0&1&0&0&17\\
2&2&0&2&0&1&4&0&0&8&0&0&0&0&2&0&19\\
3&2&4&2&0&4&0&0&8&0&0&0&0&4&4&1&29\\
1&2&4&2&0&1&0&0&8&0&0&0&0&0&0&4&21\\
5&2&0&2&0&0&0&4&0&0&0&8&4&1&0&2&23\\
\hline
&20&20&20&20&20&20&20&40&80&0&40&20&20&20&20&280\\
\end{array}\right)$$
}%
Hence we may also assume that $x_1+x_3=6$, and so there are exactly $9$
pebbles on the $9$ vertices at distance two from $r$ if we assume that
$|C|=15$.

Now, $v_1$ and $v_3$ can both reach all neighbors of $r$, so if either
$x_1\ge 4$ or $x_3\ge 4$ then we can place a pebble on any neighbor of $r$.
This forces each vertex at distance two from $r$ to have at most, and hence
exactly $1$ pebble.
Thus the splitting of the 4 pebbles will place a pebble on $r$.
This contradiction means that $x_1=x_3=3$.

Similarly, then, the pair of $v_1$ and $v_3$ can place a pebble on any
neighbor of $r$ via a common neighbor, again forcing $1$ pebble on every
vertex at distance two from $r$.
However, this allows $v_1$ and $v_3$ to each place a pebble on the same
neighbor of $r$ via disjoint paths.

This final contradiction shows $|C|\le 14$ again, finishing the proof.
\pf

Next we considered the random graph $R_{20}$ having adjacency list
\begin{quote}
$([6,8,11,12,14,15,16,17]$, $[4,5,6,7,8,10,15,16,17,18,19,20]$,\\ 
$[4,6,8,12,14,20]$, $[2,3,5,6,8,9,12,15,18,19]$, $[2,4,7,12,14,15,16,18,20]$,\\ 
$[1,2,3,4,7,8,14,15,19]$, $[2,5,6,8,11,12,13,14,15,17,18]$,\\
$[1,2,3,4,6,7,10,11,14,15,17]$, $[4,10,11,13,14,17,19,20]$,\\
$[2,8,9,16,18,19,20]$, $[1,7,8,9,13,14,16,18,20]$, $[1,3,4,5,7,13,16]$,\\
$[7,9,11,12,19,20]$, $[1,3,5,6,7,8,9,11,18]$, $[1,2,4,5,6,7,8,19]$,\\
$[1,2,5,10,11,12,18,20]$, $[1,2,7,8,9,19]$, $[2,4,5,7,10,11,14,16,20]$,\\ 
$[2,4,6,9,10,13,15,17,20]$, $[2,3,5,9,10,11,13,16,18,19])$,
\end{quote}
generated with edge probability $.4$.
The graph has diameter two, does not have the form of the Class 1
characterization, and so must be Class 0.
Indeed, from the $9,371$ basic strategies of depth at most two for the 
root $r=v_1$ (other roots had between 8,000--28,000 such strategies), 
CPLEX delivered the following certificate involving 16 of them
(which we won't bother to simplify) in Figure \ref{R20a01}.

\begin{figure}
\centerline{
{\tiny
\rotatebox{90}{
$\left(
\begin{array}{r|rrrrrrrrrrrrrrrrrrrr|r}
1469&0&2&1&0&0&4&0&0&0&0&0&0&0&0&0&0&0&0&2&0&9\\
623&0&0&0&2&0&4&0&0&0&0&0&0&0&0&1&0&0&0&2&0&9\\
473&0&0&0&2&0&0&1&4&0&2&0&0&0&2&0&0&0&0&0&0&11\\
191&0&0&2&1&0&2&2&4&0&2&0&0&0&0&2&0&0&0&0&0&15\\
1308&0&0&2&2&0&0&1&4&0&2&0&0&0&0&0&0&2&0&0&0&13\\
527&0&0&0&0&0&0&0&4&0&2&0&0&0&0&1&0&2&0&0&0&9\\
488&0&0&0&0&0&0&0&0&0&0&4&0&2&2&0&0&0&1&0&2&11\\
2012&0&0&0&0&0&0&0&0&2&0&4&0&2&1&0&0&0&2&0&2&13\\
2500&0&0&2&2&2&0&1&0&0&0&0&4&2&0&0&0&0&0&0&0&13\\
1250&0&0&0&0&2&1&0&0&2&0&0&0&0&4&0&0&0&0&0&0&9\\
221&0&0&2&0&0&0&0&0&1&0&0&0&0&4&0&0&0&2&0&0&9\\
45&0&0&2&0&0&0&0&0&2&0&0&0&0&4&0&0&0&1&0&0&9\\
4234&0&0&0&0&0&0&1&0&0&0&0&0&0&0&2&0&0&0&1&0&4\\
2500&0&2&0&0&1&0&0&0&0&2&0&0&0&0&0&4&0&2&0&2&13\\
1031&0&2&0&0&0&0&0&0&2&0&0&0&0&0&0&0&4&0&1&0&9\\
551&0&0&0&0&0&0&2&0&2&0&0&0&0&0&0&0&4&0&1&0&9\\
\hline
0&0&10000&9999&9999&10000&10000&9999&9996&9999&9998&10000&10000&10000&9998&10000&10000&9998&9999&10000&10000&189985\\
\end{array}\right)$
}
}
}
\caption{Certificate for $\pi(R_{20},v_1)\le 20$.}\label{R20a01}
\end{figure}

\begin{thm}\label{R20class0}
The graph $R_{20}$ is Class 0.
\end{thm}

\proof
The certificate for $\pi(R_{20},v_1)=20$ is shown in Figure \ref{R20a01}.
The certificates for all remaining roots are shown in the Appendix.
\pf

% ##########################################################################
%
%       GRAPH CLASSES
%
\subsection{Graph Classes}\label{Classes}

Next we turn our attention to classes of graphs, and begin with 
an extremely simple proof of the pebbling numbers of cycles, first
proved in \cite{PaSnVo}.

\begin{thm}\label{cycles}
For $k\ge 1$ we have $\pi(C_{2k})=2^k$ and 
$\pi(C_{2k+1})=\lceil (2^{k+2}-1)/3\rceil$.
\end{thm}

\proof
For both results we use two basic strategies: one path in each direction.
For even cycles the paths of length $k$
will overlap in the vertex opposite the root.
This yields $\pi(C_{2k})\le 2(2^k-1)/2+1=2^k$.
For odd cycles with $k\ge 3$ the paths will be of length $k+3$, which gives
$\pi(C_{2k+1})\le\lfloor 2(2^{k+3}-1)/(2^2+2^3)\rfloor+1
=\lfloor (2^{k+2}-1/2)/3\rfloor+1=\lceil (2^{k+2}-1)/3\rceil$.
For $k\le 2$, paths of length $k+1$ suffice.
\pf

\begin{figure}
\centerline{
\includegraphics[height=2.0in,angle=90]{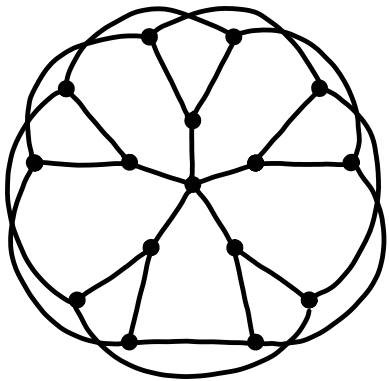}}
\caption{The generalized Petersen graph $P_{5,2}$.}\label{P52}
\end{figure}

Now we consider a generalization of the Petersen graph.
For $m\ge 3$ and $d\ge 2$ define $P_{m,d}$ to have vertices $u$ and
$v_{i,X}$, where $1\le i<m$ and $X$ is a binary $k$-tuple for $0\le k<d$.
Furthermore, $uv_{i,\mt}$ is an edge for all $i$, and $v_{i,X}v_{i,X^-}$
is an edge for all $i$ and nonempty $X$, where $X^-$ denotes the
truncation of $X$ obtained from dropping its final digit.
Finally, for every $i$ and length $d-1$ $X$, we include the edge
$v_{i,X}v_{i+1,X}$ (addition modulo m), with the exception that when
$i=m-1$ we use $v_{m-1,X}v_{0,X^+}$ instead, where $X^+$ denotes the 
$(d-1)$-tuple that satisfies $N(X^+)=N(X)+1$, and $N(X)$ is the
natural number represented by $X$ in binary.
Figure \ref{P52} shows the graph $P_{5,2}$; it is easy to check
that $P_{3,2}$ is the Petersen graph $P$.
Also, $P_{m,d}$ has diameter $2d$ when $m>3$ and $2d-1$ when $m=3$.

Interest regarding these graphs comes from two sources.
In \cite{PaSnVo} the problem is raised of finding the minimum number
$e(n)$ of edges in a Class 0 graph on $n$ vertices.
Because of the Petersen and Wheel graphs, we have $e(n)\le 2n-2$.
Blasiak, et al. \cite{BCFHS}, show that $e(n)\ge\lfloor 3n/2\rfloor$, and
conjecture that $e(n)=3n/2+o(1)$.
Furthermore, they believe that the graphs $P_{m,d}$ are all Class 0
which, if true, would prove the conjecture because it is $3$-regular
except for the central vertex $u$, having degree $m$ (with $n=m(2^d-1)+1$).
The graphs also appear in \cite{DanVol} as $n$-vertex graphs having the 
fewest edges among those of minimum degree 3 and radius $d$.
In fact, the intuition for the conjecture comes from this result.

Here, we are most interested in the graphs $P_{m,2}$ for $m>3$.
Note that $P_{m,2}$ has rotational symmetry: the exceptional edges
for $k=2$ can be ``rotated onward'' by swapping $v_{1,0}$ with
$v_{1,1}$ in the drawing.
This makes $P_{m,2}$ transitive on the set $\{v_{i,X}\}$ for fixed
length $X$.
Thus, when calculating $\pi(P_{m,2})$, one only need consider the
three vertices $u$, $v=v_{1,\mt}$, and $w=v_{1,0}$ as root.

\begin{figure}
\centerline{
\includegraphics[height=1.4in,angle=90]{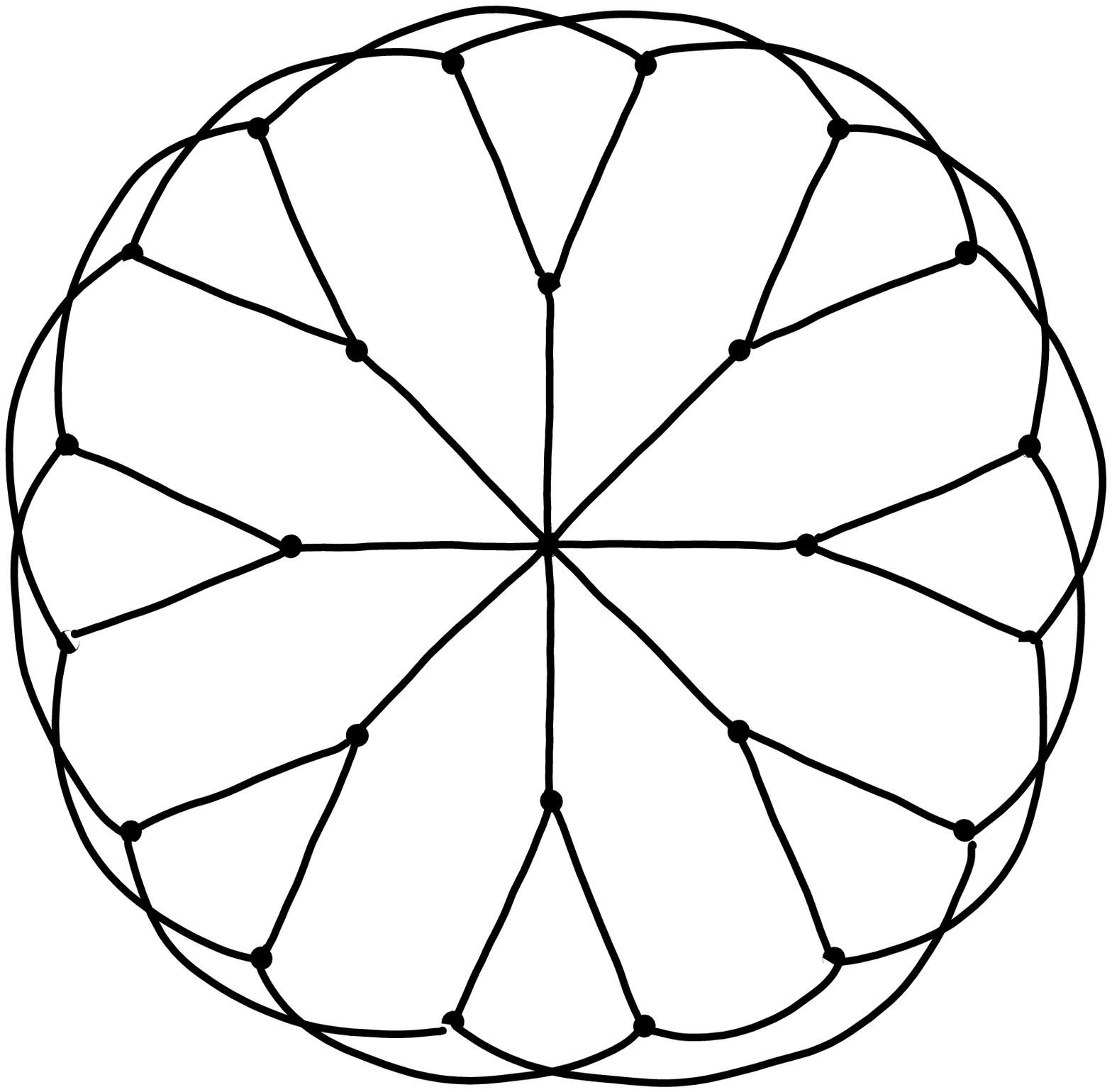}$\qquad$
\includegraphics[height=1.3in,angle=90]{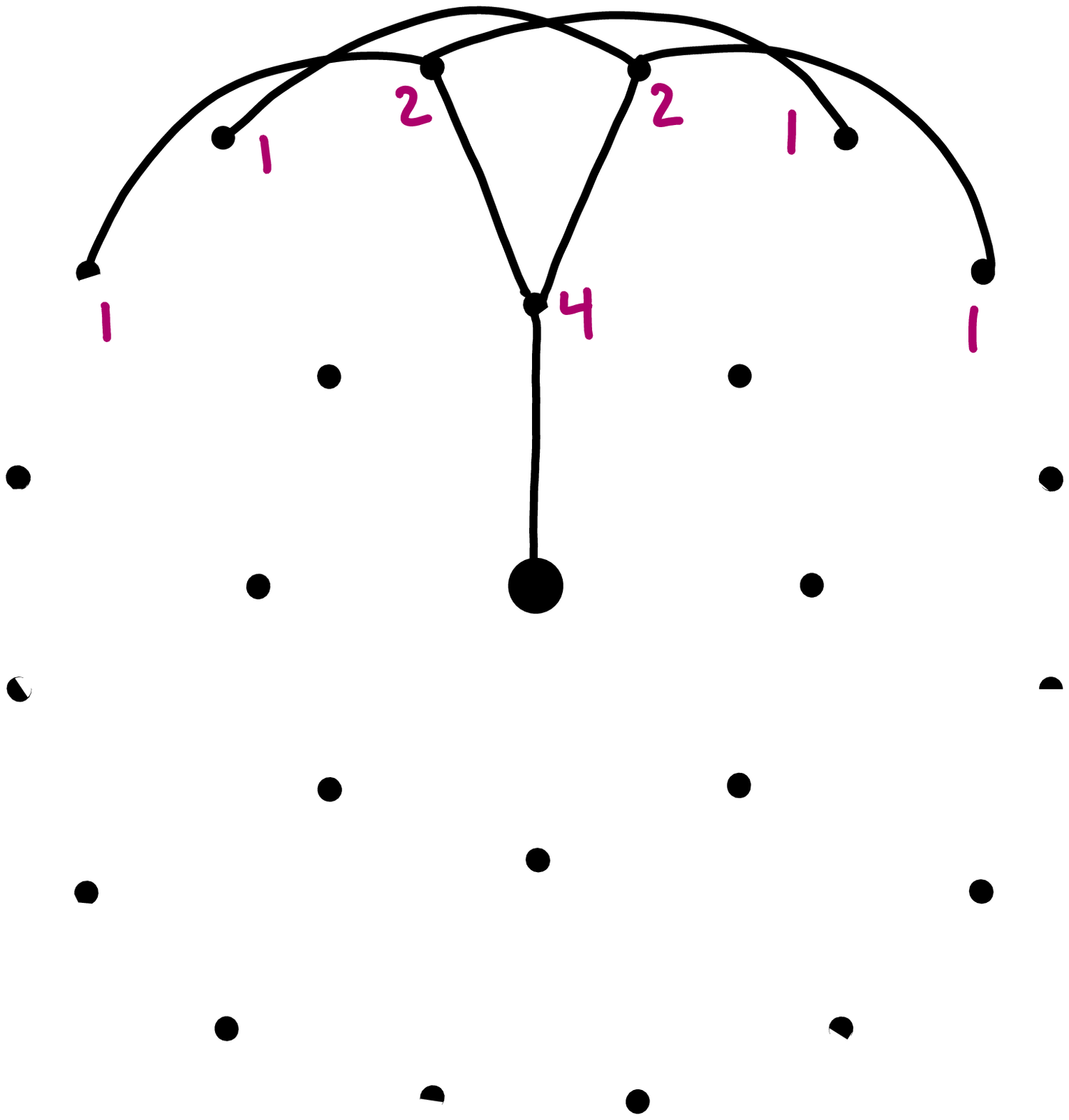}}
\caption{The graph $P_{8,2}$ and its main strategy $T$ for root $u$.}\label{Stratu}
\end{figure}

\begin{thm}\label{Pm2pi}
For all $m\ge 4$ we have the following, where $n=n(P_{m,2})=3m+1$:
\begin{enumerate}
\item
$\pi(P_{m,2},u)=n$,
\item
$\pi(P_{m,2},v)\le n+5$, and
\item
$\pi(P_{m,2},w)\le n+17$.
\end{enumerate}
\end{thm}

\proof
When the root is $u$, we use all rotations of the following basic strategy
$T$ (see Figure \ref{Stratu}).
Let $v_{1,\mt}$ have weight $4$, $v_{1,0}$ and $v_{1,1}$ have weight $2$,
and each of their two neighbors have weight $1$.
Clearly, of the sum of all rotations of $T$ has weight $4$ everywhere
but $u$, and so the Uniform Covering Lemma applies.

\begin{figure}
\centerline{
\includegraphics[height=1.3in,angle=90]{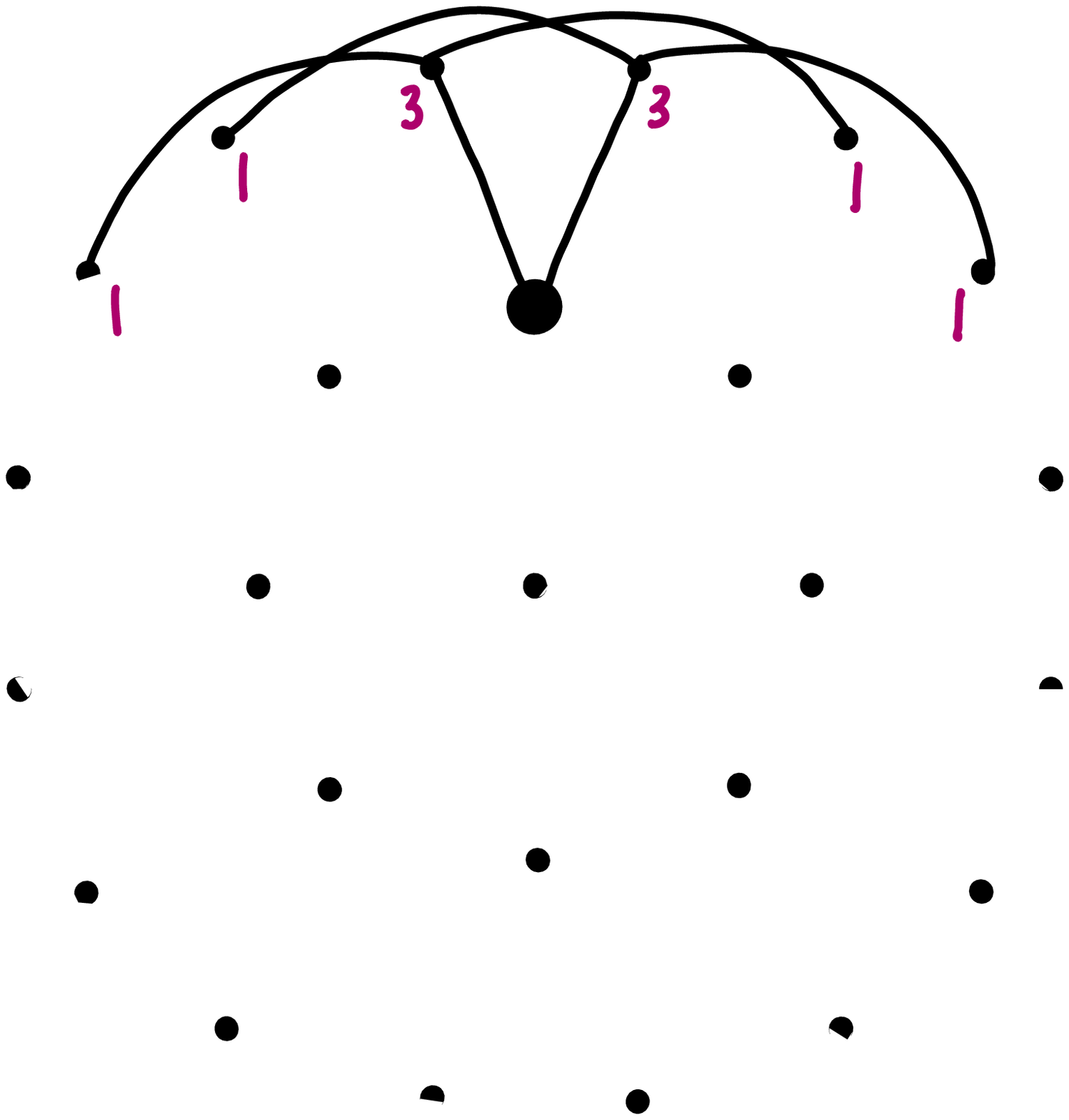}$\quad$
\includegraphics[height=1.4in,angle=90]{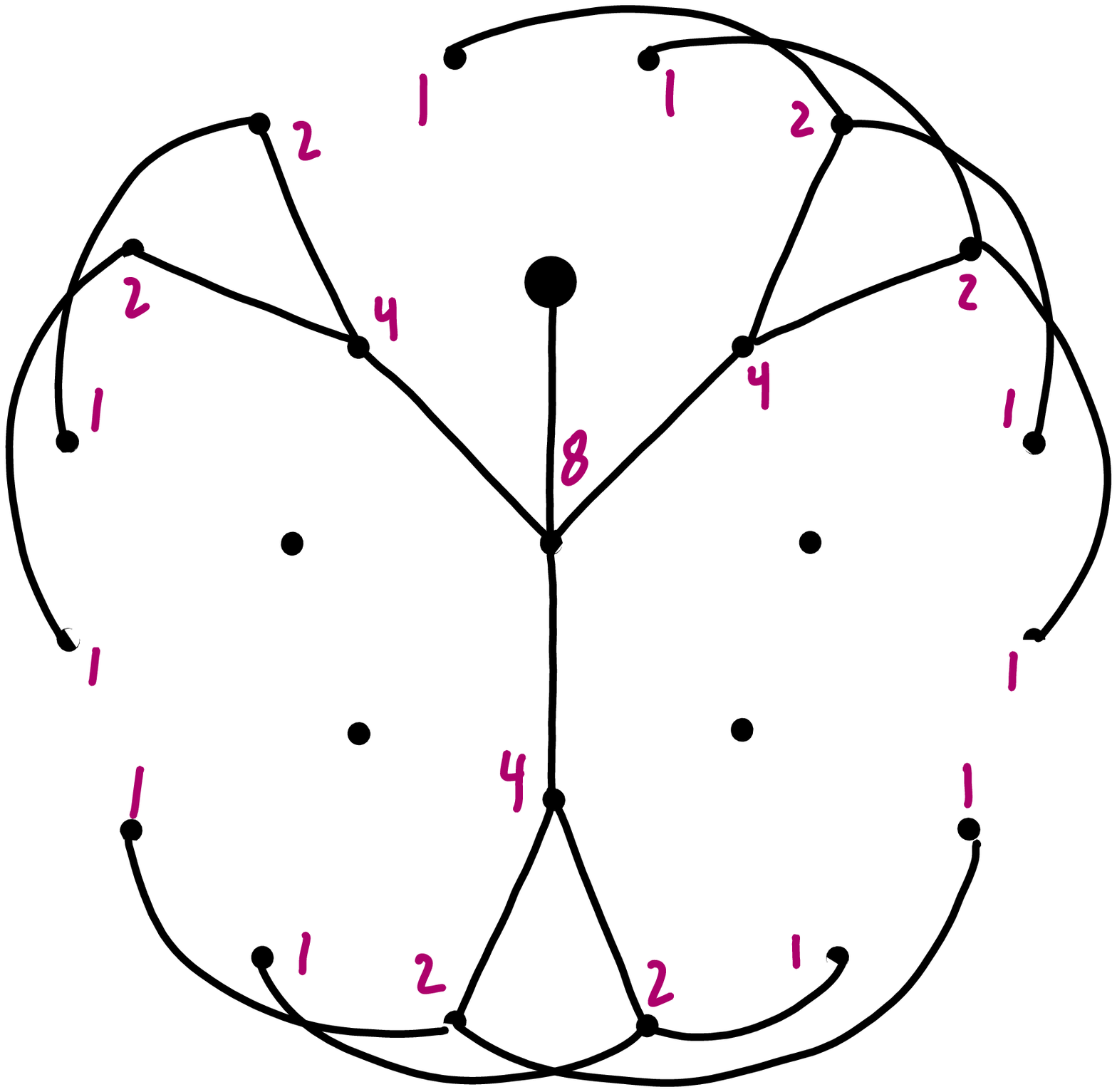}$\quad$
\includegraphics[height=1.4in,angle=90]{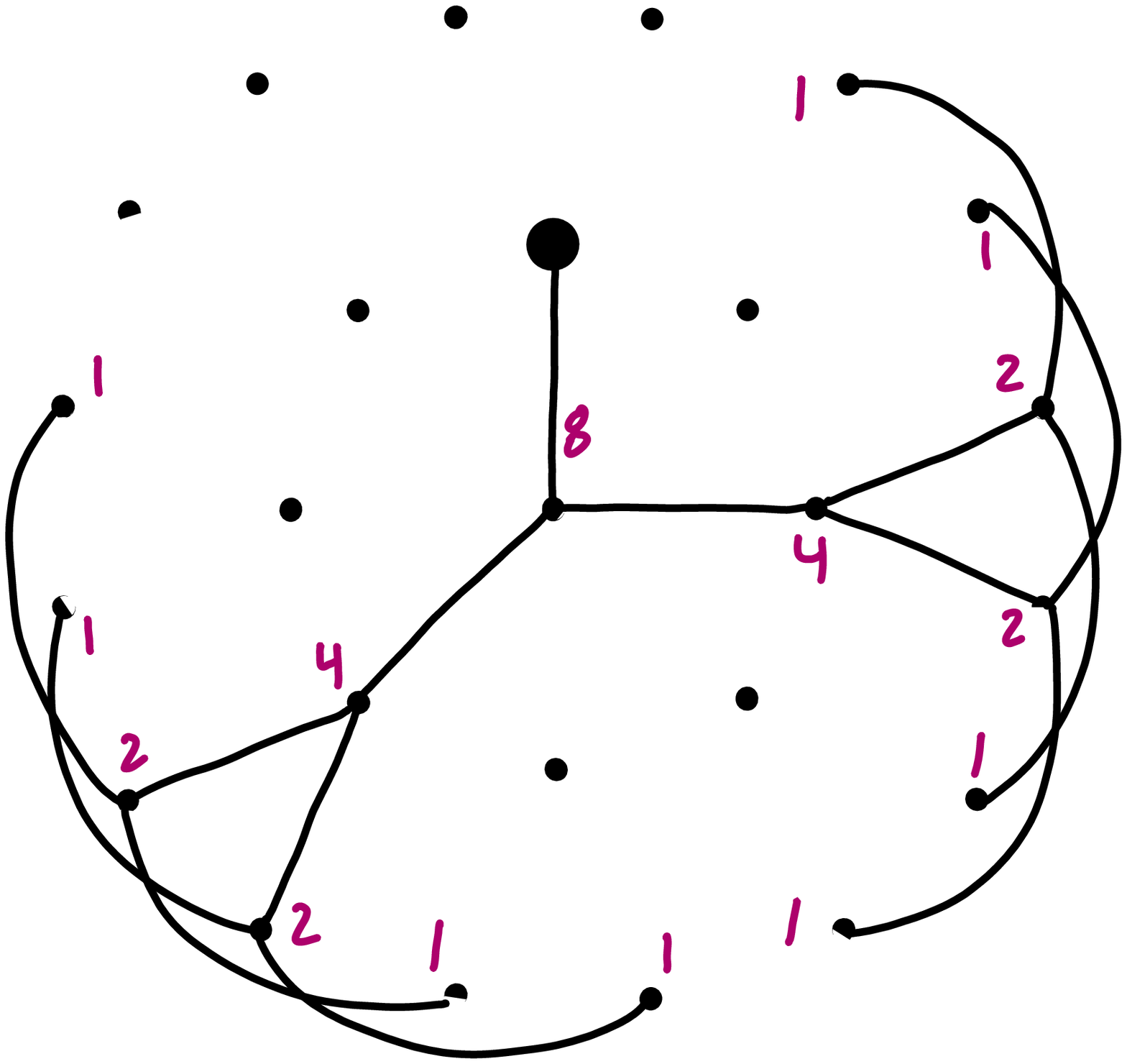}$\quad$
\includegraphics[height=1.4in,angle=90]{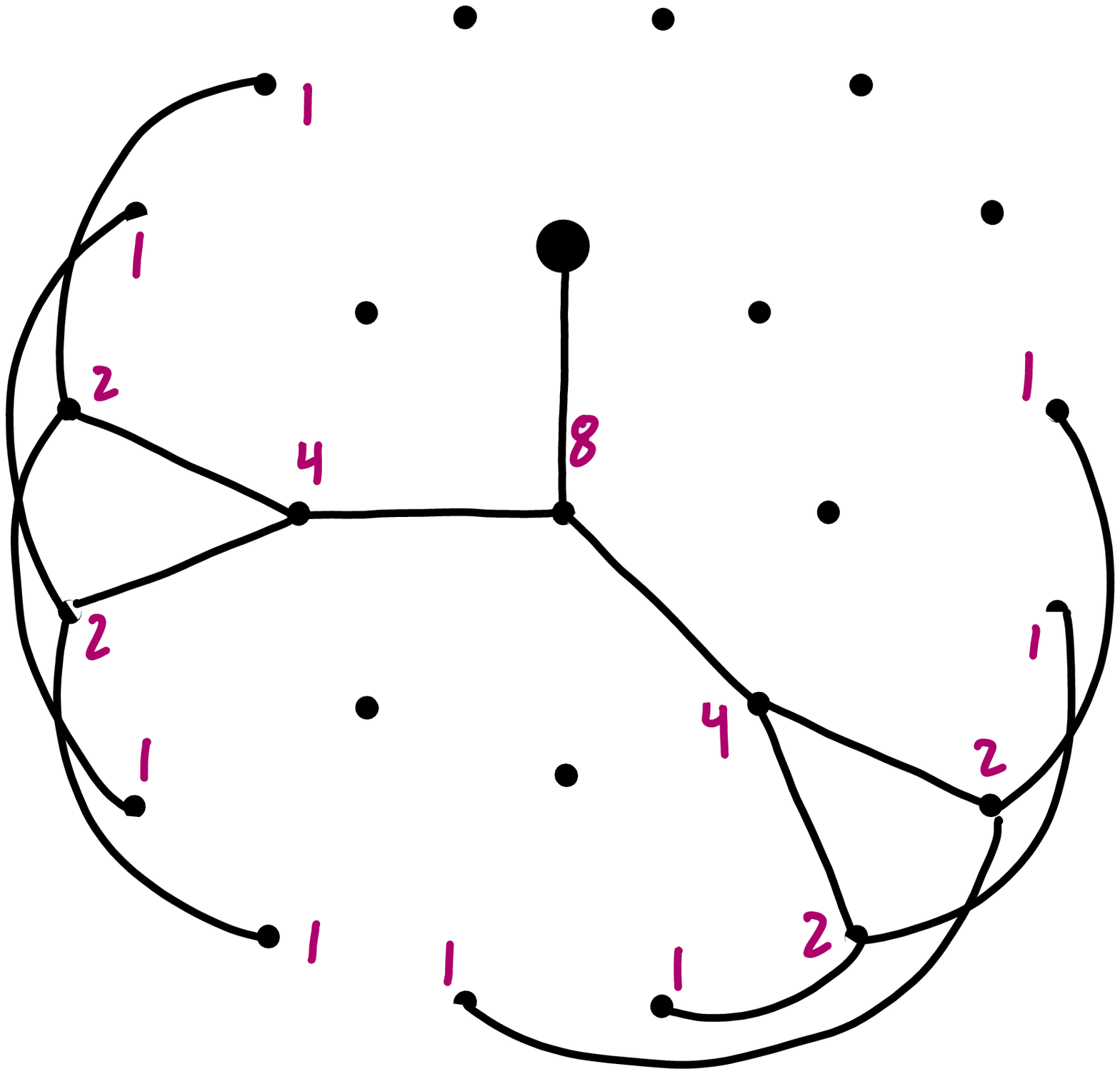}}
\caption{The four strategies for root $v$ in $P_{8,2}$.}\label{Stratv}
\end{figure}

When the root is $v$, we build slightly more complex strategies.
First we use the nonbasic strategy $S$, having weight $3$ on $v_{0,0}$ 
and $v_{0,1}$ and weight $1$ on each of their neighbors.
Next, write the rotations of $T$ as $T_0,\ldots,T_{m-1}$.
For $j\in\{0,1,2\}$, we build the basic strategies $S_j$ from these by 
combining all $T_i$ ($i\not=0$) with $i\equiv j\mod 3$, with the
exception that if $j=1$ and $i=m-1$ then we do not include the vertices
$v_{0,0}$ and $v_{0,1}$ twice.
Moreover, each $S_j$ includes weight $8$ on $u$ (see Figure \ref{Stratv}).
One quarter of the sum of these four strategies has weight $0$ on $v$,
weight $6$ on $u$, and weight $1$ everywhere else.
Now the Weight Function Lemma (actually Corollary \ref{Upper}) applies.

\begin{figure}
\centerline{
\includegraphics[height=1.3in,angle=90]{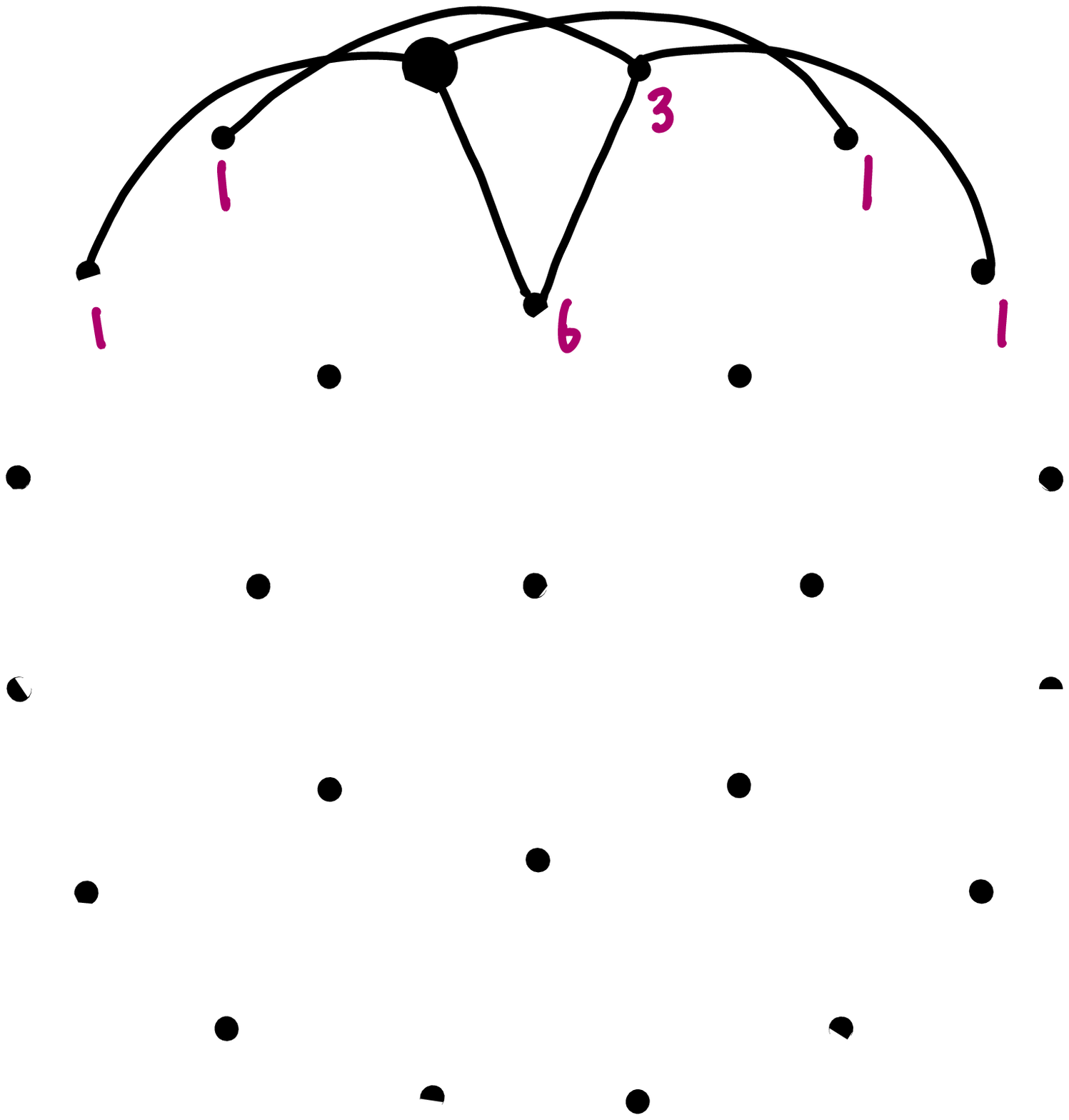}$\quad$
\includegraphics[height=1.4in,angle=90]{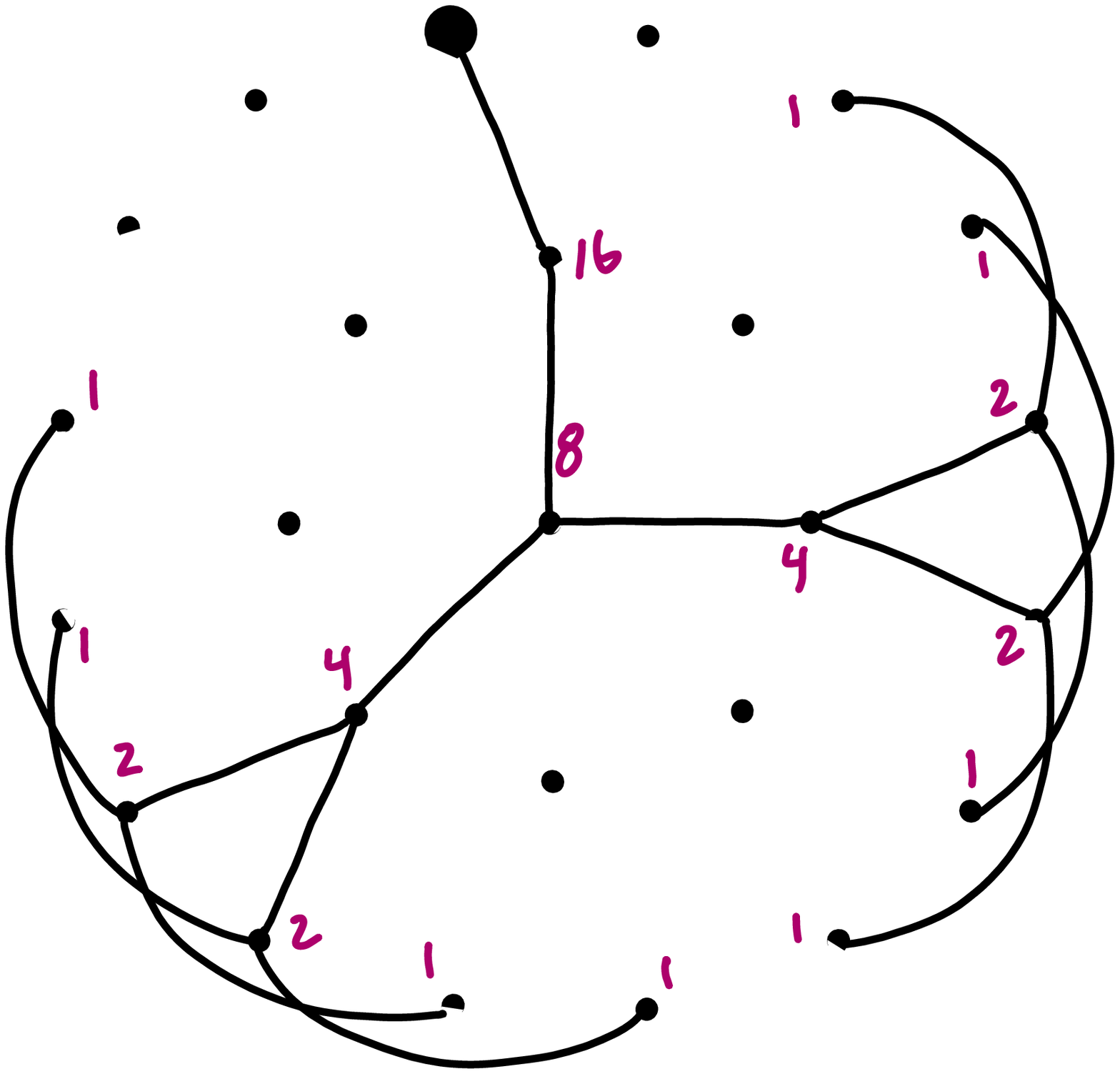}$\quad$
\includegraphics[height=1.4in,angle=90]{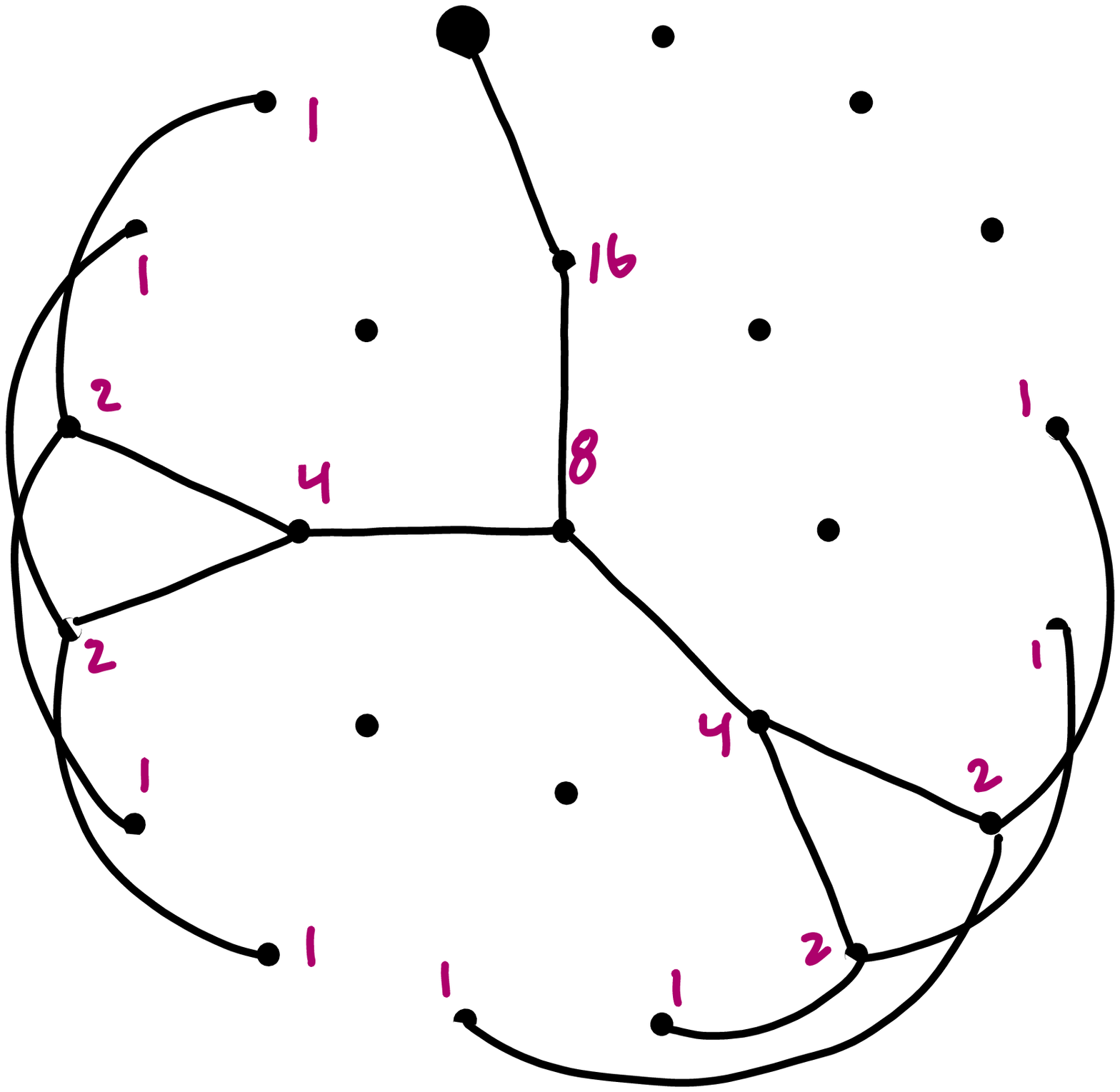}$\quad$
\includegraphics[height=1.4in,angle=90]{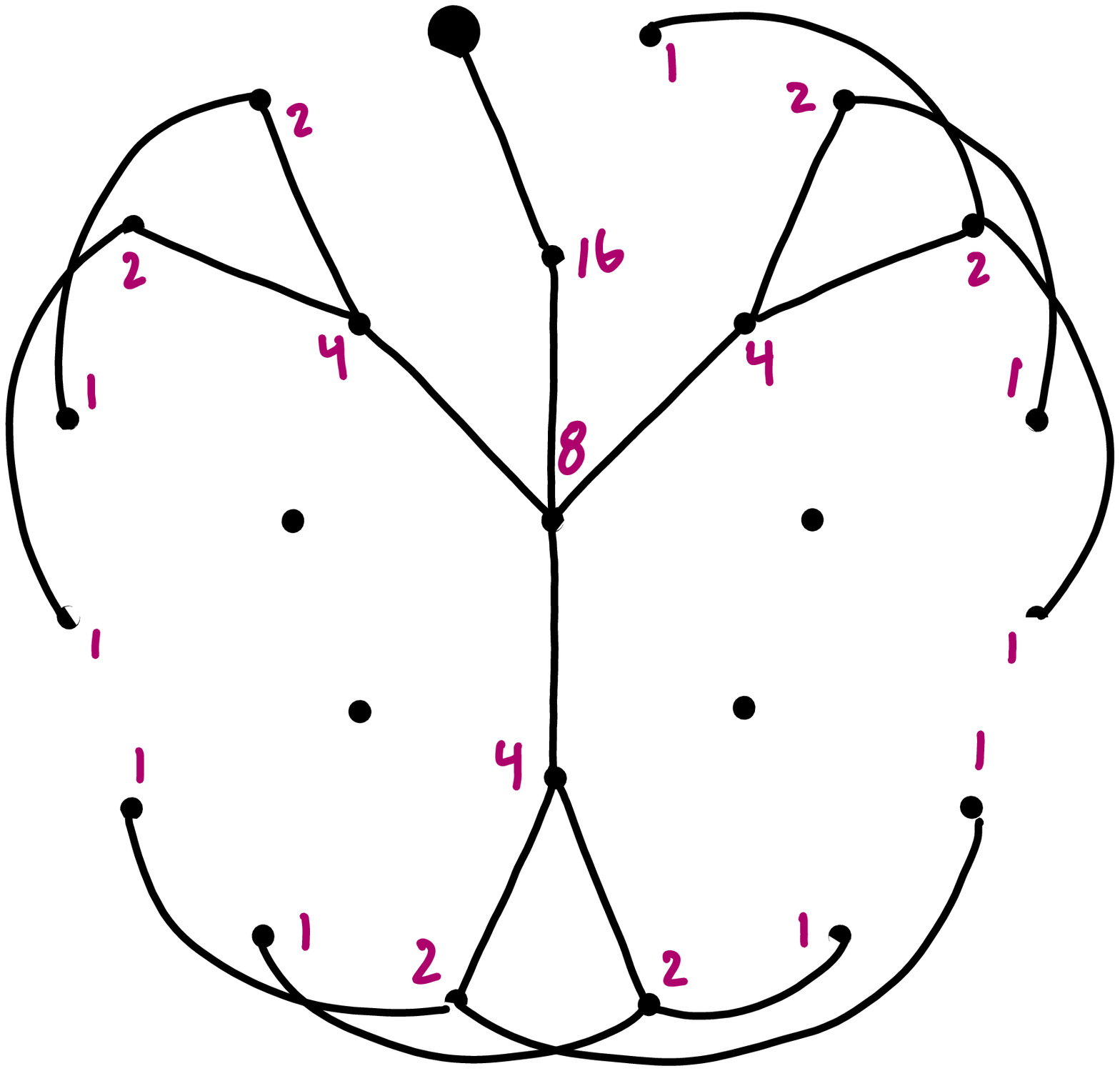}}
\caption{The four strategies for root $w$ in $P_{8,2}$.}\label{Stratw}
\end{figure}

The description for the strategies when the root is $w$ is almost
identical to that for root $v$ (see Figure \ref{Stratw}).
In this case, one quarter of the sum of the four strategies has weight
$0$ on $w$, weight $13.5$ on $v$, weight $6$ on $u$, and weight $1$
elsewhere.
\pf

Next we define the generalized Coxeter $p$-graph $C(p)$ for odd primes 
$p=2q+1$.
Like the graphs $P_{m,d}$, $C(p)$ has the potential (for $p\ge 5$)
to be a Class 0 graph with fewer than $2n-2$ edges.
Set $V_i=\{(i,j)\mid 0\le j<p\}$ and $V=\cup_{i=0}^qV_i$,
so that $C(p)$ has $n=p(q+1)=\binom{p+1}{2}$ vertices.
Define the $2pq=2n-2p$ edges $(i,j)(i,j+1)$ and $(0,j)(i,j)$ 
for each $1\le i\le q$ and $0\le j<p$.
The vertices in $V_0$ have degree $q$ and all others have degree 3,
making $C(7)$, the original Coxeter graph, 3-regular.
Here we show the following $n+O(\sqrt n)$ bound.

\begin{thm}\label{cox}
For every prime $p\ge 11$ we have $\pi(C(p))\le n+14q+14$.
\end{thm}

\proof
First we exhibit the symmetry of $G=C(p)$.
To do so, we note that the arithmetic that follows will be modular in $q$
in the first coordinate, with the speciality that we use the representative
$q$ in place of $0$, and modular in $p$ (as normal) in the second.
Let $\r$ be the automorphism that sends vertices $(i,j)$ to $(i,j+1)$
for all $0\le i\le q$ and $0\le j<p$.
This is a rotation of the vertices within each $V_i$.
Let $\a$ generate $Z_p^*$ and
define $\s$ to be the automorphism that sends $(0,j)$ to $(0,\a j)$ and 
$(i,j)$ to $(i+1,\a j)$ for all $1\le i\le q$ and $0\le j<p$.
While permuting $V_0$, the key aspect of $\s$ is that it also
permutes the sets of vertices from $V_i$ to $V_{i+1}$.
Together, $\r$ and $\s$ act transitively on $V_0$ and on $\overline{V_0}$.
This means that we only need calculate $\pi(G,(0,0))$ and $\pi(G,(1,0))$.
(It turns out that $C(7)$ is fully vertex transitive and so the root
$(0,0)$ suffices.)

Next we define the $q$ basic strategies $\bT_i$ ($1\le i\le q$) for the
root $r=(0,0)$, with each root child having weight 4.
The child of $r$ in $\bT_i$ has weight 4 and two children, $(i,i)$ and
$(i,-i)$.
We describe the descendents of $(i,i)$ only, as those of $(i,-i)$ are
identical but with their second coordinate negated.
The left subtree under $(i,i)$ is the path $(i,2i)$, $(i,3i)$, $\ldots$,
$(i,qi)$.
The right subtree under $(i,i)$ is a collection of $q-2$ paths,
starting with children $(k,i)$ for $1\le k\le q$ with 
$k\not\in\{i,(p\pm i)/2\}$.
The child $(k,i)=(k,sk)$, for $s=ik^{-1}$, and from it hangs the path
$(k,(s+1)k)$, $\ldots$, $(k,qk)$.

Now we look at the sum of the weights of vertices over all strategies,
and by the symmetries described above we need only consider the vertices
$(i,j)$ for $i\in\{0,1\}$ and $0\le j\le q$.
We note first that $(0,j)$ is the grandchild of $(j,0)$ in $\bT_j$, and
so has weight 1.
The left path under $(1,1)$ shows that vertex $(1,j)$ has weight $2^{2-j}$
in $\bT_i$.
For $j\ge 3$ we see that if $k\le j$ then vertex $(1,j)$ has weight
$2^{k-j-1}$.
Therefore, except for the $q$ weight 4 children of $r$ and their $2q$ weight 2 
children, every other nonroot vertex has weight 1.
This gives total weight $(n-1)+(4-1)q+(2-1)2q$, implying that
$\pi(G,r)\le n+5q$.

Consider the other root $r=(1,0)$.
Here we need many more basic strategies --- we will define $\bT_i$ for 
$i\in\{\pm 1\}\cup\{2,\ldots,q\}$ --- and the weights of the children of
$r$ equal to 8 instead of 4.
Similar to above, the strategy $\bT_{-1}$ will be identical to $\bT_i$
except that the second coordinates will be negated.
We will describe the first four levels of each strategy explicitly,
and their remaining levels implicitly.
The child of $r$ in $\bT_1$ is $(1,1)$, having children $(1,2)$ and $(0,1)$.
From $(1,2)$ hangs the path $(1,3)$ and $(1,4)$, while the children of
$(0,1)$ are $(2,1),\ldots (q,1)$.
The children of each $(i,1)$ are $(i,1\pm i)$, subject to them not having
already listed in some $\bT_j$, which we now describe.
The child of $r$ in $\bT_i$ is $(0,0)$, which has the single child $(i,0)$.
The two children of $(i,0)$ are $(i,\pm i)$, the two children of which are
$(i,\pm 2i)$ and $(0,\pm i)$, correspondingly, with the exception that 
$(i,\pm 2i)$ does not appear in $\bT_q$ because it already appears in 
$\bT_{\pm 1}$ (as $(q,\pm 1)$).
At this point, there are $q+1$ vertices of weight 8, $q+3$ vertices
of weight 4, and $4q-2$ vertices of weight 2, with all other vertices
having weight 1.

The remaining construction of the strategies proceeds in stages, with each
new vertex added with several fractional weights adding up to 1.
In the first stage, each vertex $(a,b)$ that is adjacent to a current leaf 
$(a,c)$ of some strategy is added to that strategy as a child of $(a,c)$, 
accounting for weight $1/2$.
The other $1/2$ weight comes from adding it as a child of $(0,b)$ as well.
In each subsequent stage, every vertex $(a,b)$ that is adjacent to a
current leaf $(a,c)$ is added as a child of $(a,c)$ in every strategy
that contains $(a,c)$, accounting for weight $1/2$.
The other $1/2$ weight comes from adding it as a child of $(0,b)$ as well.

Hence we obtain $\pi(G,r)\le (n-1)+7(q+1)+3(q+3)+1(4q-2)+1=n+14q+14$.
\pf

We note that the number of strategies for the root $(1,0)$ can be reduced
whenever two strategies (not including $\bT_{\pm 1}$) share no vertices
(other than $(0,0)$).
Thus one can define the intersection graph $G_p$ of the family of sets
$W_i$, $2\le i\le q$, where $W_i=V(\bT_i)-\{(0,0)\}$.
With chromatic number $\x_p=\x(G_p)$, we obtain the upper bound of
$n+7\x+7q+21$.
For example, $\x_p\le 1$ when $p\le 7$.
In those cases, further improvements can be made as well, and it is not
too hard to show that $\pi(C(5))\le n+6$ and $\pi(C(7))\le n+15$.

Finally, we discuss powers of cycles.
For a given graph $G$ and integer $k$, we denote by $G^{(k)}$ the graph on
the same vertex set as $G$, with edges $uv$ whenever the distance 
$\dist_G(u,v)\le k$ in $G$. 
For example, $G^{(n-1)}=K_n$ for every connected $G$, where $n=n(G)$.
Pachter, et al \cite{PaSnVo}, define the {\it pebbling exponent} of $G$ 
to be the minimum $e=e_\pi(G)$ for which $G^{(e)}$ is Class 0.
Consequently $e_\pi(G)\le n-1$ for all $G$.
The problem raised in \cite{PaSnVo} is to find $e_\pi(C_n)$.
Here we prove the following.

\begin{thm}\label{PebbExpo}
The pebbling exponent of the cycle $C_n$ satisfies 
$$\frac{n/2}{\lg n}\le e_\pi(C_n)\le\frac{n/2}{\lg n-\lg\lg n}\ .$$
\end{thm}

\proof
The lower bound follows from the general fact that 
$\pi(G)\ge 2^{\diam(G)}$ for all $G$, along with the observation that 
$\diam(C_n^{(e)})=\lceil n/2e\rceil$.
Therefore, a requirement for Class 0 is that $n\ge 2^{\lceil n/2e\rceil}$.

For the upper bound, we prove that $\pi(C_n^{(2^k)})=n$ for 
$n=(2k+1)2^k+3=f(k)$.
Then we show that $\pi(C_n^{(2^k)})=n$ for $f(k-1)<n<f(k)$.
Our method will be to split the cycle into two identical paths with endpoints
at the root $r$, and use identical strategies on each one.
We will invoke the Uniform Covering Lemma \ref{ucl} to obtain the result.

In more detail, let us give a useful labeling of the vertices $V$ of $C_n$
as follows.
We partition $V=\cup_{i=0}^{k+2}(U_i\cup W_i)$ so that
\begin{itemize}
\item
each $U_i$ and $W_i$ induces a path in $C_n$,
\item
$U_0=W_0=\{r\}$,
\item
$U_{i+1}$ follows $U_i$ when traversing $C_n$ clockwise from $r$,
\item
$W_{i+1}$ follows $W_i$ when traversing $C_n$ counterclockwise from $r$, and
\item
$U_{k+2}=W_{k+2}$.
\end{itemize}
The two identical paths mentioned above are seen to be $U=\cup_{i=0}^{k+2}U_i$
and $W=\cup_{i=0}^{k+2}W_i$, where now we see that the term `split' was a
white lie because of the overlap on $U_{k+2}=W_{k+2}$.
We will describe a family $\cT$ of strategies on $U$ with the property that,
for some $m$ (actually $2^{k+1}$) and all $u\in\cup_{i=1}^{k+1}U_i$, we have
$\sum_{T\in\cT}T(u)=m$, while for all $u\in U_{k+2}$ we have
$\sum_{T\in\cT}T(u)=m/2$.
Then we copy these strategies symmetrically onto $W$ and, because of the
overlap, the Uniform Covering Lemma \ref{ucl} applies.

Next we describe the vertices within each $U_i$.
First, $|U_i|=2^k$ for $i\in\{1,k+2\}$, and $|U_i|=2^k-2^{k-i+1}$
for $2\le i\le k+1$ (thus $n=(2k+1)2^k+3$).
We clockwise order the vertices of $U_i=\{v_{i,0},\ldots,v_{i,|U_i|-1}\}$
in their natural order by subscript, and will find it useful to identify
$v_{i,j}$ with the encoding $[i,b_j]$, where $b_j$ is the $k$-bit binary
representation of $j$ (so that leading zeros are not suppressed).
For example, $v_{3,6}$ is encoded as $[3,00110]$ in $C_{355}$ since
$k=5$ in that case.
Also, the root $r$ is encoded as $[0,0^k]$, where we write $x^j$ to denote
the concatenation $xx\cdots x$ of length $j$.
Furthermore, we use the notation $\bv_j$ to mean some binary word of length $j$,
subject to context.
That is, $[1,\bv_k]$ is always a vertex in $U_1$ but $[2,\bv_k]$ is not
always in $U_2$ --- every vertex in $U_2$ looks like $[2,0\bv_{k-1}]$.
Moreover, for $2\le i\le k+1$, vertices in $U_i$ look like everything
but $[i,1^{i-1}\bv_{k-i+1}]$.

Now we describe a partially ordered set $P$ that we will use to define
our strategies.
The elements of $P$ are the vertices $U$, and the covering relations
are given by
\begin{enumerate}
\item
$[0,0^k]>[1,\bv_k]$ for all $\bv_k$,
\item
$[1,\bv_{k-1}x]>[2,0\bv_{k-1}]$ for all $\bv_{k-1}$ and all $x\in\{0,1\}$,
\item\label{stepi1}
$[i,\bv_{i-1}\bv_{k-i}x]>[i+1,\bv_{i-1}b\bv_{k-i}]$ for all
$2\le i\le k$, all $\bv_{i-1}\not=1^{i-1}$, all $\bv_{k-i}$, 
and each $x,b\in\{0,1\}$, 
\item\label{stepi2}
$[i,1^{i-2}0\bv_{k-i}x]>[i+1,1^{i-1}0\bv_{k-i}]$ for all
$2\le i\le k$, all $\bv_{k-i}$, and all $x\in\{0,1\}$,
\item
$[k+1,\bv_k]>[k+2,\bv_k]$ for all $\bv_k\not=1^k$, and
\item
$[k+1,1^{k-1}0]>[k+2,1^d]$.
\end{enumerate}
All other relations of $P$ are determined by transitivity.

For $z\in P$, define its {\it downset} $D(z)=\{y\in P\mid y<z\}$.
Notice that each $D=D([1,\bv_k])$ forms a tree in $P$.
Indeed, simple induction shows that if $\bv_k=v_1\cdots v_k$ then
every vertex of $U_{i+1}$ in $D$ has the form 
$[i+1,\bx_iv_1\cdots v_{k-i}]$, where $\bx_i$ is anything but $1^i$.
This means that $[i+1,x_1\cdots x_iv_1\cdots v_{k-i}]$ has exactly 
one element from $D$ that covers it, namely either
$[i,x_1\cdots x_{i-1}v_1\cdots v_{k-i}v_{k-i+1}]$ or
$[i,x_2\cdots x_iv_1\cdots v_{k-i}v_{k-i+1}]$.
So no cycles exist in $D$.  

For each $\bv_k$, then, define the basic strategy $T=T(\bv_k)$ to have
vertices $\{r\}\cup D([1,\bv_k])$, with edge $zy$ whenever $z>y$
is a covering relation in $P$.
In order that $T$ is a strategy in $C_n^{(2^k)}$ we must verify that
the distance $d=\dist(z,y)$ between $z$ and $y$ in $C_n$ is at most $2^k$.
When $z=r$, we have $y=v_{1,j}$ for some $0\le j\le 2^k-1$, so $d=j+1$.
When $z=v_{1,j}$, we have $y=v_{2,\lfloor j/2\rfloor}$, so 
$d=2^k-\lceil j/2\rceil$.
When $z=v_{i,j}$ and $2\le i\le k+1$, then we can write $j=h2^{k-i+1}+l$
for some $0\le l<2^{k-i+1}$.
From relation \ref{stepi1} or \ref{stepi2} we have $y=v_{2,j^\pr}$, where 
$j^\pr=h2^{k-i+1}+b2^{k-i}+\lfloor l/2\rfloor$ and $b\in\{0,1,2\}$.
Thus the greatest distance comes from relation \ref{stepi2},
in which case $d=|U_i|-j+j^\pr=2^k-\lceil l/2\rceil$.  

Note that the characterization of elements in $D([1,\bv_k])$
implies that each vertex of $U_i$ is in $2^{i-1}$ basic strategies
when $1\le i\le k+1$, and in $2^k$ when $i=k+2$.
Because each $T$ is basic, this means that $\sum_TT(v)=2^{k+1}$
for all $v\in U_i$ (resp. $W_i$), $1\le i\le k+1$, and $2^k$
for $v\in U_{k+2}$ (resp. $W_{k+2}$).
Now the overlap from $U_{k+2}=W_{k+2}$ gives sum $m=2^{k+1}$
for all $v\not=r$, and the Uniform Covering Lemma \ref{ucl} applies.

Finally, whenever $n$ is smaller than $f(k)$ we simply erase sufficiently
many vertices of $U_{k+2}=W_{k+2}$ (but don't renumber any indicies/encodings),
as they are leaves in all strategies and so don't destroy the uniform
covering.
When such vertices are exhausted continue erasing vertices of
$U_{k+1}\cup W_{k+1}$ with similar results.
Since $f(k)-f(k-1)=2^{k+1}$, no more considerations are necessary.  
\pf

%\bigskip

% ##########################################################################
%
%       REMARKS
%
\section{Remarks}\label{Remarks}

In this paper we have shown several different strengths of the
Weight Function Lemma in combination with linear optimization,
highlighting its versatility.
It has been used to compute upper bounds on and exact values of
the pebbling number of small graphs.
It has also been successful in calculating the pebbling numbers of much
larger graphs than previous algorithms.
For such graphs having too many strategies than time allows to construct,
the technique of creating a smaller set of them at random seems to 
perform just as well.
This is most likely due to the property that nonoptimal solutions
(derived from having fewer constraints) seem, in most instances, to be 
near optimal (have the same floor function).
In fact, by restricting strategies to be breadth-first search, one
obtains upper bounds on the greedy pebbling numbers of graphs (which
requires pebbling steps to move toward the root).
The method also yields results for many families of graphs, in many
cases by hand, with much simpler and remarkably shorter proofs than 
given in previously existing arguments.
This is especially so with highly symmetric graphs.
We note also that the technique can be used in conjunction with more
traditional arguments, as in Theorem \ref{R15class0}, and it has
delivered an array of upper bounds, such as $n$, $n+c$, and $n+o(n)$,
most of which are the best known and might possibly be best possible.
It's two main shortcomings are the inability to overcome the kind of
splitting structure found in cubes, for example, in which solutions
to some configurations require nontree solutions, and the difficulty
in dealing with large diameter, although success has been found with
cycles and their graph powers, in addition to Petersen and Coxeter
generalizations.
When the technique gives upper bounds, it would be of great use to
know how good the bound might be.
That is, does the Weight Function Lemma yield an approximation
algorithm for graph pebbling?

\begin{qst}\label{ApproxAlgo}
Is there a constant $c$ such that, for all graphs $G$ and every
root $r\in V(G)$, $z_{G,r}\le c\pi(G,r)$?
\end{qst}

For example, is $c=2$?
The cube $Q^d$ shows that it couldn't be any smaller.

\begin{thm}\label{Qapprox}
For all $d\ge 1$ and all $r\in V(Q^d)$ we have $z_{Q^d,r}<2\pi(Q^d,r)$.
\end{thm}

\proof
We exploit the symmetry of $Q^d$ as follows.
Because $Q^d$ is vertex transitive we need only consider one root $r$.
We identify $V(Q^d)$ with the power set of $\{1,\ldots,d\}$ and
take $r=\mt$.
We define a single strategy $\bT$ for $r$ and then apply every 
permutation of $\{1,\ldots,d\}$ to obtain other strategies.
Finally we average over this collection of strategies.
The result will be that $\pi(Q^d)$ will be at most one more than the
sum of the weights in this average.

For each $1\le k\le d$ such that $2^k\le\binom{d}{k}$ we define 
$a_k=\lfloor\binom{d}{k}/2^k\rfloor$ and $b_k=\binom{d}{k}\mod 2^k$.
We form $\bT$ by first taking a neighbor of $r$.
This neighbor will be a set of size 1, to which we assign the weight $2^{d-1}$
--- this is step $i=1$, with $k=d-i$.
For future steps $i$, while $2^k>\binom{d}{k}$, we continue adding a
single neighbor, a set of size $i$, to the current leaf of $\bT$, and
assign the weight $2^k$.
If $2^k\le\binom{d}{k}$, however, we add $a_k+1$ vertices to the current
leaves: $a_k$ of these will have weight $2^k$ and one will have weight $b_k$.
All of these will be connected to leaves of weight $2^{k+1}$ and none
will be connected to leaves of weight $b_{k+1}$.
This is possible because the degree of each vertex of size $i-1$ is $k+1$,
and $a_k<(k+1)a_{k+1}$ whenever $k+1<d-1$, so there are enough potential
neighbors to accomplish this.

Hence we obtain weight $1$ on average for vertices of size $i$ for which
$2^k\le\binom{d}{k}$, and weight $2^k/\binom{d}{k}$ otherwise.
This yields the bound $z_{Q^d,r}\le 1+\sum_{k=0}^{d-1}\max\{\binom{d}{k},2^k\}
<\sum_{k=0}^{d}\binom{d}{k}+\sum_{k=0}^{d-1}2^k<2^{d+1}=2\pi(Q^d)$.
\pf

If one restricts their attention to only polynomially many strategies, 
this linear optimization technique becomes a polynomial algrorithm.
It would be useful to investigate how good the approximation can be
under these circumstances.

%\bigskip

% ##########################################################################
%
%       ACKNOWLEDGEMENTS
%
\section{Acknowledgements}\label{Acknowledgements}

We thank Andrzej Czygrinow for converting the author's Maple code for
generating all the tree strategies of a rooted graph into java.

%\bigskip

% ##########################################################################
% ##########################################################################
%
%       BIBLIOGRAPHY
%

% ##########################################################################
% ##########################################################################
%
%       APPENDIX
%
\section{Appendix}\label{Append}

\subsection{Certificates for $R_{15}$}\label{R15}

We list the missing certificates in order, with the all-zeros
column signifying the root.
The format is the same as that for $R_{15}$ at $v_{10}$ in Section \ref{Random}.
No attempt was made to find the simplest set of strategies.

{\tiny
$$\left(
\begin{array}{r|rrrrrrrrrrrrrrr|r}
577&0&8&4&0&0&0&2&4&0&2&0&0&2&1&2&25\\
40&0&8&0&0&0&2&0&3&0&2&4&0&0&4&2&25\\
1250&0&0&4&8&0&0&4&0&4&2&0&0&0&3&2&27\\
542&0&0&0&0&8&0&2&4&0&2&4&1&2&2&2&27\\
190&0&0&0&0&8&4&0&4&0&2&1&0&0&0&4&23\\
217&0&0&0&0&8&4&1&4&0&2&4&0&2&0&4&29\\
62&0&0&4&0&0&8&2&0&4&2&1&0&0&0&0&21\\
646&0&0&0&0&0&8&0&2&4&2&1&0&4&4&0&25\\
94&0&4&2&0&0&0&2&2&1&2&4&8&0&2&2&29\\
91&0&8&4&0&0&1&0&4&2&2&4&2&2&4&2&35\\
300&0&0&2&0&8&4&2&1&2&2&4&2&2&2&4&35\\
344&0&4&0&0&0&2&2&2&0&2&4&8&4&1&2&31\\
175&0&4&2&0&0&1&0&0&2&2&4&8&4&2&2&31\\
472&0&4&2&0&0&1&2&2&2&2&4&8&4&0&2&33\\
\hline
&0&10004&10002&10000&9992&9998&9999&10000&10002&10000&9998&10004&10002&10001&9998&14000\\
\end{array}\right)$$

$$\left(
\begin{array}{r|rrrrrrrrrrrrrrr|r}
4&0&0&0&0&0&0&0&0&0&0&0&1&0&0&0&1\\
4&0&0&0&0&0&0&0&0&0&0&0&0&0&1&0&1\\
2&4&0&0&0&2&2&0&0&0&0&0&1&2&0&0&11\\
2&0&0&4&0&0&2&1&0&0&0&0&0&0&0&0&7\\
1&0&0&0&4&0&0&0&0&2&0&0&0&0&1&0&7\\
1&0&0&0&4&0&0&2&0&2&0&0&0&0&1&0&9\\
2&0&0&0&0&0&0&0&4&0&2&0&0&2&1&2&11\\
2&0&0&0&0&2&0&2&0&2&2&4&1&0&0&2&15\\
\hline
&8&0&8&8&8&8&8&8&8&8&8&8&8&8&8&112\\
\end{array}\right)$$

$$\left(
\begin{array}{r|rrrrrrrrrrrrrrr|r}
2&0&8&0&0&2&0&0&4&0&2&1&4&2&0&0&23\\
6&2&8&0&0&0&0&0&4&0&2&1&4&2&4&2&29\\
34&2&8&0&0&0&0&0&4&0&2&4&4&2&3&2&31\\
23&4&0&0&8&4&0&0&0&4&2&0&1&2&4&0&29\\
1&0&0&0&0&0&8&0&0&4&2&1&2&4&0&0&21\\
42&0&0&0&0&0&0&8&2&0&2&4&1&0&0&4&21\\
16&4&0&0&8&4&0&0&2&4&2&1&2&2&4&0&33\\
3&4&0&0&8&4&0&0&0&4&2&0&2&2&3&2&31\\
19&0&0&0&0&4&8&0&2&4&2&0&1&4&0&2&27\\
7&4&0&0&0&4&8&0&2&4&2&1&2&4&0&2&33\\
15&4&0&0&0&4&8&0&0&4&2&0&2&4&3&2&33\\
\hline
&336&336&0&336&336&336&336&336&336&336&336&336&336&336&336&4704\\
\end{array}\right)$$

$$\left(
\begin{array}{r|rrrrrrrrrrrrrrr|r}
135&0&0&2&0&0&0&4&0&0&0&1&0&0&0&2&9\\
31&0&0&0&0&0&4&0&0&8&4&3&0&2&0&0&21\\
2&0&8&0&0&0&0&2&0&1&2&4&4&2&0&2&25\\
1&2&4&8&0&0&4&0&0&1&0&2&2&2&0&0&25\\
38&4&0&0&0&8&0&2&3&0&2&4&2&2&0&4&31\\
10&2&0&0&0&0&4&0&0&8&4&4&2&0&1&2&27\\
1&0&0&2&0&0&4&0&0&8&4&4&2&2&2&1&29\\
103&0&4&2&0&1&0&0&4&0&2&0&2&2&8&2&27\\
39&8&0&0&0&4&4&0&0&0&0&0&2&4&1&2&25\\
30&8&2&1&0&0&4&0&2&0&0&0&4&4&0&0&25\\
49&2&8&4&0&0&0&2&4&0&2&0&4&2&0&1&29\\
4&2&8&4&0&0&1&2&0&0&2&4&4&0&4&0&31\\
54&0&0&2&0&8&4&2&1&0&2&4&2&2&2&4&33\\
82&2&1&2&0&0&4&2&2&8&4&4&2&2&0&0&33\\
\hline
&996&998&1000&0&995&996&998&1000&995&996&994&996&998&999&996&13957\\
\end{array}\right)$$

$$\left(
\begin{array}{r|rrrrrrrrrrrrrrr|r}
10&0&0&0&0&0&0&0&0&0&0&0&0&0&0&1&1\\
6&4&0&0&0&0&0&0&0&0&0&0&2&1&0&0&7\\
2&0&2&2&4&0&0&1&0&0&0&0&0&0&0&0&9\\
1&0&0&2&4&0&0&2&0&0&0&0&0&0&1&0&9\\
3&0&0&2&4&0&0&2&0&2&0&0&0&0&1&0&11\\
6&0&0&2&0&0&4&0&0&1&0&0&0&2&2&0&11\\
2&0&0&0&0&0&0&0&4&0&2&0&0&1&0&0&7\\
2&0&2&0&0&0&0&0&4&0&2&0&0&0&2&1&11\\
2&0&2&0&0&0&0&0&4&0&2&0&0&2&2&1&13\\
6&0&2&0&0&0&0&2&0&2&2&4&2&0&0&1&15\\
2&0&0&0&0&0&0&1&0&0&0&0&0&0&0&2&3\\
\hline
&24&24&24&24&0&24&24&24&24&24&24&24&24&24&24&336\\
\end{array}\right)$$

$$\left(
\begin{array}{r|rrrrrrrrrrrrrrr|r}
7&8&0&0&4&0&0&2&0&1&0&0&4&4&0&0&23\\
5&8&0&0&0&4&0&0&0&0&0&0&4&3&0&2&21\\
79&0&0&0&0&0&0&2&0&8&2&4&2&0&0&1&19\\
1&0&4&2&3&0&0&0&0&0&0&0&2&0&8&0&19\\
14&8&0&2&4&0&0&1&0&0&0&2&4&4&2&0&27\\
120&2&4&8&3&0&0&4&0&0&0&0&2&0&0&2&25\\
31&4&2&0&0&8&0&2&0&1&2&4&0&2&0&4&29\\
91&0&0&0&4&8&0&2&1&0&2&4&2&0&0&4&27\\
4&0&2&0&3&0&0&0&2&8&4&4&2&0&0&2&27\\
20&4&0&0&0&0&0&0&4&0&2&0&1&8&0&0&19\\
85&4&0&0&0&0&0&0&4&0&2&0&2&8&0&1&21\\
104&0&4&0&0&0&0&0&4&0&2&0&1&0&8&0&19\\
37&0&0&0&4&0&0&2&1&8&4&4&0&0&2&2&27\\
7&2&4&1&4&0&0&2&4&0&2&0&2&0&8&2&31\\
\hline
&1006&998&997&999&996&0&998&1000&998&998&996&1002&1001&998&998&13985\\
\end{array}\right)$$

$$\left(
\begin{array}{r|rrrrrrrrrrrrrrr|r}
99&2&4&8&0&0&4&0&1&0&0&0&2&2&2&0&25\\
26&2&4&8&0&0&4&0&2&2&0&0&2&2&1&0&27\\
2&4&0&0&8&4&1&0&0&4&2&0&2&2&4&0&31\\
59&0&0&0&0&0&0&0&1&2&2&4&2&1&0&0&12\\
51&2&4&0&0&0&0&0&0&3&4&8&4&2&2&0&29\\
31&0&0&0&0&4&0&0&1&4&4&8&4&2&0&0&27\\
12&2&0&0&0&4&0&0&4&0&2&0&0&2&1&8&23\\
65&0&0&0&0&4&2&0&4&0&2&0&0&2&1&8&23\\
48&2&1&0&0&4&2&0&4&0&2&0&0&2&2&8&27\\
123&4&2&0&8&3&2&0&2&4&2&0&2&2&4&0&35\\
13&2&0&0&0&0&2&0&1&4&4&8&4&2&0&0&27\\
\hline
&998&998&1000&1000&1001&1000&0&1000&999&998&996&998&999&999&1000&13986\\
\end{array}\right)$$

$$\left(
\begin{array}{r|rrrrrrrrrrrrrrr|r}
203&0&0&0&0&0&0&0&0&2&4&1&0&0&0&0&7\\
208&2&0&0&0&0&1&0&0&0&0&0&2&4&0&0&9\\
7&0&0&2&4&8&0&0&0&2&0&3&0&0&0&0&19\\
66&0&0&0&4&0&4&1&0&2&0&0&0&0&8&0&19\\
37&2&3&2&4&0&4&2&0&0&0&0&0&0&8&0&25\\
116&0&0&2&1&0&0&4&0&0&0&0&0&0&0&8&15\\
38&4&8&4&0&0&2&2&0&2&1&4&2&2&4&0&35\\
30&4&0&2&0&8&4&0&0&2&0&0&2&0&1&0&23\\
5&2&0&2&4&8&4&2&0&0&0&0&0&1&0&0&23\\
77&2&0&2&4&8&0&2&0&2&1&4&2&0&0&0&27\\
6&0&0&2&0&8&4&2&0&2&0&4&0&0&0&1&23\\
37&0&8&4&0&0&2&2&0&2&2&4&4&2&0&1&31\\
6&2&8&4&0&0&1&2&0&2&0&4&4&2&0&0&29\\
31&2&8&4&4&0&2&2&0&2&0&4&4&0&0&1&33\\
\hline
&1000&1007&1004&1008&1000&1002&1004&0&1002&1001&1004&1002&999&1006&1002&14041\\
\end{array}\right)$$

$$\left(\begin{array}{r|rrrrrrrrrrrrrrr|r}
82&2&0&2&4&0&0&8&0&1&0&0&0&0&2&0&19\\
15&0&0&2&4&0&0&8&0&2&0&0&0&0&2&3&21\\
98&0&0&0&0&0&2&0&1&0&0&0&8&4&0&0&15\\
3&0&0&0&0&0&0&0&4&0&2&0&0&0&1&8&15\\
12&0&0&0&0&0&0&0&4&0&2&0&0&1&2&8&17\\
4&0&8&4&4&0&1&0&0&2&0&0&4&0&4&0&27\\
38&4&8&4&0&0&0&1&4&0&0&0&4&2&4&2&33\\
11&0&8&4&4&0&2&0&4&2&0&0&3&0&4&0&31\\
67&2&0&2&0&8&4&0&1&2&0&0&0&2&2&4&27\\
58&4&0&2&0&8&4&1&0&2&0&0&0&0&2&4&27\\
14&2&0&4&4&0&1&8&0&0&0&0&0&0&2&0&21\\
120&0&0&0&2&0&1&0&4&4&8&0&0&2&0&2&23\\
64&4&8&4&4&0&2&0&1&2&0&0&0&2&4&0&31\\
8&4&8&4&0&0&2&2&4&0&1&0&2&2&4&2&35\\
\hline
&998&1000&1000&1000&1000&1000&1000&997&1000&998&0&1001&998&999&997&13988\\
\end{array}\right)$$

$$\left(
\begin{array}{r|rrrrrrrrrrrrrrr|r}
28&8&0&2&4&3&0&2&0&2&0&0&0&0&2&0&23\\
4&8&0&0&4&2&3&2&0&2&0&0&0&0&2&0&23\\
3&0&0&0&0&4&0&4&0&1&4&8&0&0&0&0&21\\
4&0&0&2&0&3&0&4&0&0&4&8&0&0&0&0&21\\
27&8&0&2&4&4&2&2&0&1&0&0&0&0&2&2&27\\
33&0&8&4&4&1&2&2&4&0&2&0&0&0&4&2&33\\
52&0&0&0&0&2&1&4&0&4&4&8&0&0&0&4&27\\
59&0&0&2&0&1&4&0&4&2&2&0&0&8&2&2&27\\
26&0&8&4&4&2&2&2&4&2&2&0&0&0&4&1&35\\
\hline
&472&472&472&472&472&472&472&472&472&472&472&0&472&472&472&6608\\
\end{array}\right)$$

$$\left(
\begin{array}{r|rrrrrrrrrrrrrrr|r}
8&8&0&2&4&0&0&2&0&1&0&0&0&0&2&0&19\\
10&0&0&0&0&0&0&2&0&1&2&4&8&0&0&2&19\\
15&8&0&2&4&4&0&2&0&0&0&0&0&0&2&1&23\\
14&8&0&0&4&4&0&2&0&2&0&0&0&0&2&1&23\\
5&0&2&4&0&3&8&2&0&4&2&2&0&0&4&0&31\\
13&0&4&0&2&1&0&2&8&0&4&2&0&0&4&4&31\\
27&0&4&2&1&0&0&2&0&2&2&4&8&0&2&2&29\\
2&0&2&4&0&4&8&2&0&4&2&2&0&0&0&1&29\\
17&0&0&4&2&4&8&2&0&4&2&2&0&0&0&1&29\\
13&0&2&4&1&4&8&2&0&4&2&2&0&0&0&2&31\\
24&0&4&2&2&1&0&2&8&2&4&2&0&0&4&4&35\\
\hline
&296&296&296&296&296&296&296&296&296&296&296&296&0&296&296&4144\\
\end{array}\right)$$

$$\left(
\begin{array}{r|rrrrrrrrrrrrrrr|r}
19&2&4&2&0&4&0&0&8&0&4&2&2&0&0&3&31\\
109&0&8&4&0&0&0&2&0&0&1&4&4&2&0&2&27\\
3&0&0&2&8&3&0&4&0&4&2&2&0&0&0&0&25\\
22&0&0&0&0&4&1&2&8&0&0&2&2&4&0&4&27\\
27&2&0&0&0&1&0&0&8&0&4&0&0&4&0&4&23\\
30&0&0&0&0&2&1&0&8&0&4&0&2&4&0&4&25\\
16&0&0&0&0&0&0&2&8&2&4&0&0&3&0&4&23\\
3&2&8&4&0&2&1&2&4&0&2&0&4&0&0&0&29\\
7&4&0&4&8&4&0&1&2&4&2&2&0&0&0&2&33\\
116&4&0&0&8&4&0&4&0&4&2&2&2&1&0&2&33\\
54&4&0&4&0&2&8&2&0&4&2&2&2&1&0&0&31\\
28&2&0&4&0&0&8&1&0&4&2&2&2&4&0&0&29\\
35&4&0&4&0&4&8&2&0&4&2&2&0&3&0&2&35\\
8&0&4&2&0&0&1&2&8&0&4&0&2&4&0&4&31\\
\hline
&1002&1004&1004&1008&1006&999&1005&1002&1004&1001&1004&1002&1001&0&1003&14045\\
\end{array}\right)$$

$$\left(
\begin{array}{r|rrrrrrrrrrrrrrr|r}
14&4&2&2&4&8&2&0&0&1&0&0&2&2&0&0&27\\
6&4&2&2&1&8&4&0&0&2&0&0&2&2&0&0&27\\
1&4&2&0&2&8&4&0&0&1&0&0&2&2&2&0&27\\
5&4&0&2&2&8&4&0&0&1&0&0&2&2&2&0&27\\
2&2&1&4&4&0&2&8&0&0&0&0&0&0&0&0&21\\
14&2&2&4&4&0&1&8&0&0&0&0&0&0&2&0&23\\
10&2&0&4&4&0&1&8&0&2&0&0&0&0&2&0&23\\
2&2&4&1&2&0&2&0&8&2&4&0&2&4&0&0&31\\
24&2&1&0&0&0&2&0&8&2&4&0&2&4&4&0&29\\
26&0&4&2&1&0&2&0&0&4&4&8&4&2&2&0&33\\
\hline
&208&208&208&208&208&208&208&208&208&208&208&208&208&208&0&2912\\
\end{array}\right)$$
}

\subsection{Certificates for $R_{20}$}\label{R20}
We list the missing certificates in order, as above.

%\newpage

\tikzstyle{mybox} = []
\begin{tikzpicture}[transform shape, rotate=90, baseline=6.8in]
\node [mybox] (box) {%
    \begin{minipage}[t!]{0.5\textwidth}

{\tiny
$\left(
\begin{array}{r|rrrrrrrrrrrrrrrrrrrr|r}
354&0&0&2&4&0&0&0&0&1&0&0&2&0&0&0&0&0&0&0&0&9\\
917&0&0&2&4&2&1&0&0&2&0&0&2&0&0&0&0&0&0&0&0&13\\
1046&0&0&2&4&0&0&0&0&2&0&0&0&0&0&0&0&0&0&1&0&9\\
1390&0&0&0&0&4&0&0&0&0&0&0&2&0&0&0&2&0&2&0&1&11\\
3783&1&0&1&0&0&2&0&0&0&0&0&0&0&1&1&0&0&0&1&0&7\\
2019&0&0&0&0&0&0&4&0&0&0&0&2&2&0&1&0&2&0&0&0&11\\
320&0&0&0&0&2&2&4&0&0&0&2&2&2&0&2&0&1&0&0&0&17\\
214&2&0&1&0&0&0&0&4&0&0&2&0&0&0&0&0&2&0&0&0&11\\
2136&2&0&0&0&0&0&0&4&0&1&2&0&0&2&0&0&2&0&0&0&13\\
71&2&0&2&0&0&2&0&4&0&2&1&0&0&2&0&0&2&0&0&0&17\\
321&2&0&0&0&0&0&2&1&0&0&0&0&0&0&4&0&0&0&0&0&9\\
368&2&0&0&2&2&2&0&0&0&0&0&0&0&0&4&0&0&0&1&0&13\\
394&0&0&0&0&0&0&0&0&0&0&2&0&0&1&0&2&0&4&0&0&9\\
1411&0&0&0&0&0&0&0&0&0&2&2&0&0&1&0&2&0&4&0&0&11\\
401&0&0&0&0&0&0&0&0&2&2&0&0&1&0&2&0&2&0&4&0&13\\
410&0&0&0&0&0&0&0&0&2&0&0&0&2&0&0&0&0&0&4&1&9\\
490&0&0&0&0&0&0&0&0&2&2&2&0&2&0&0&1&0&0&0&4&13\\
945&0&0&0&0&0&0&0&0&2&2&0&0&2&0&0&2&0&0&1&4&13\\
615&0&0&2&0&2&0&0&0&2&2&0&0&2&0&0&2&0&0&1&4&17\\
\hline
0&10003&0&10003&10004&10000&10001&9998&10005&10002&10002&10001&10000&9999&10002&10000&10000&10002&10000&10001&10000&190023\\
\end{array}\right)$
}
\vspace{0.2 in}

{\tiny
$\left(
\begin{array}{r|rrrrrrrrrrrrrrrrrrrr|r}
354&0&0&2&4&0&0&0&0&1&0&0&2&0&0&0&0&0&0&0&0&9\\
917&0&0&2&4&2&1&0&0&2&0&0&2&0&0&0&0&0&0&0&0&13\\
1046&0&0&2&4&0&0&0&0&2&0&0&0&0&0&0&0&0&0&1&0&9\\
1390&0&0&0&0&4&0&0&0&0&0&0&2&0&0&0&2&0&2&0&1&11\\
3783&1&0&1&0&0&2&0&0&0&0&0&0&0&1&1&0&0&0&1&0&7\\
2019&0&0&0&0&0&0&4&0&0&0&0&2&2&0&1&0&2&0&0&0&11\\
320&0&0&0&0&2&2&4&0&0&0&2&2&2&0&2&0&1&0&0&0&17\\
214&2&0&1&0&0&0&0&4&0&0&2&0&0&0&0&0&2&0&0&0&11\\
2136&2&0&0&0&0&0&0&4&0&1&2&0&0&2&0&0&2&0&0&0&13\\
71&2&0&2&0&0&2&0&4&0&2&1&0&0&2&0&0&2&0&0&0&17\\
321&2&0&0&0&0&0&2&1&0&0&0&0&0&0&4&0&0&0&0&0&9\\
368&2&0&0&2&2&2&0&0&0&0&0&0&0&0&4&0&0&0&1&0&13\\
394&0&0&0&0&0&0&0&0&0&0&2&0&0&1&0&2&0&4&0&0&9\\
1411&0&0&0&0&0&0&0&0&0&2&2&0&0&1&0&2&0&4&0&0&11\\
401&0&0&0&0&0&0&0&0&2&2&0&0&1&0&2&0&2&0&4&0&13\\
410&0&0&0&0&0&0&0&0&2&0&0&0&2&0&0&0&0&0&4&1&9\\
490&0&0&0&0&0&0&0&0&2&2&2&0&2&0&0&1&0&0&0&4&13\\
945&0&0&0&0&0&0&0&0&2&2&0&0&2&0&0&2&0&0&1&4&13\\
615&0&0&2&0&2&0&0&0&2&2&0&0&2&0&0&2&0&0&1&4&17\\
\hline
0&10003&0&10003&10004&10000&10001&9998&10005&10002&10002&10001&10000&9999&10002&10000&10000&10002&10000&10001&10000&190023\\
\end{array}\right)$
}
\vspace{0.2 in}

{\tiny
$\left(
\begin{array}{r|rrrrrrrrrrrrrrrrrrrr|r}
283&0&0&0&1&4&0&2&0&4&2&4&0&2&8&2&2&1&0&2&1&35\\
1023&0&0&0&2&0&0&0&0&0&0&0&0&0&0&1&0&0&1&0&0&4\\
449&0&0&0&8&0&1&0&0&4&2&0&0&0&0&0&2&0&4&2&0&23\\
346&0&0&0&8&1&0&1&0&1&0&1&0&2&0&2&0&2&0&4&0&22\\
7&0&1&0&2&16&0&8&0&4&4&4&4&16&0&0&8&2&16&8&32&125\\
1829&0&2&0&0&2&0&0&0&1&2&2&0&0&0&0&2&0&2&1&4&18\\
789&0&2&0&0&2&0&4&8&2&4&2&0&2&4&4&0&4&0&2&0&40\\
595&0&2&0&0&4&0&1&0&0&0&1&8&4&0&0&4&0&2&0&0&26\\
133&0&4&0&2&0&8&4&1&0&2&0&0&2&0&2&0&2&2&4&0&33\\
5&2&0&0&0&0&1&8&0&32&16&16&0&4&64&4&0&0&32&0&0&179\\
63&4&4&0&0&0&0&1&2&16&8&8&0&16&0&0&16&2&4&16&32&129\\
226&4&4&0&0&1&8&4&2&1&2&1&0&2&2&4&2&2&0&4&0&43\\
299&8&0&0&0&0&0&1&2&4&0&4&16&8&0&0&0&4&0&0&0&47\\
230&8&0&0&0&0&0&2&0&0&0&0&0&1&16&0&0&4&0&2&0&33\\
148&8&0&0&0&0&0&4&16&0&0&4&0&2&0&0&0&8&2&0&1&45\\
117&8&0&0&0&0&16&0&0&2&2&0&0&0&0&8&0&4&0&4&0&44\\
65&8&0&0&0&0&16&8&0&0&2&1&0&0&2&8&0&4&4&4&0&57\\
12&8&1&0&8&8&2&0&1&0&2&0&16&0&0&4&4&4&4&4&0&66\\
233&8&8&0&4&2&16&8&0&2&0&0&1&0&0&8&4&4&4&0&0&69\\
\hline
0&9998&9997&0&9997&9994&9990&9993&10001&9999&9996&10000&9997&9992&10002&9995&9998&9997&9993&9991&9987&189917\\
\end{array}\right)$
}

    \end{minipage}
    };
\end{tikzpicture}

\tikzstyle{mybox} = []
\begin{tikzpicture}[transform shape, rotate=90, baseline=6.8in]
\node [mybox] (box) {%
    \begin{minipage}[t!]{0.5\textwidth}

{\tiny
$\left(
\begin{array}{r|rrrrrrrrrrrrrrrrrrrr|r}
470&0&4&0&0&2&2&2&1&0&0&0&0&0&0&0&0&2&0&0&0&13\\
831&0&4&0&0&0&0&2&0&0&2&0&0&0&0&2&2&1&2&0&0&15\\
1118&0&4&0&0&0&0&2&0&0&2&0&0&0&0&1&2&2&0&0&2&15\\
680&0&0&4&0&0&2&0&0&0&0&0&2&0&2&0&0&0&0&0&1&11\\
10&0&2&0&0&4&0&2&0&0&0&0&2&0&1&0&0&0&0&0&0&11\\
2255&0&0&0&0&4&0&2&0&0&0&0&1&0&0&0&2&0&0&0&0&9\\
221&1&1&1&0&0&2&1&1&0&0&0&0&0&1&0&0&0&0&0&0&8\\
1074&2&0&0&0&0&4&0&0&0&0&0&0&0&0&1&0&0&0&0&0&7\\
42&0&2&0&0&0&1&2&4&0&0&2&0&0&0&0&0&2&0&0&0&13\\
1651&2&0&2&0&0&1&0&4&0&2&2&0&0&2&0&0&2&0&0&0&17\\
61&0&0&0&0&0&2&1&4&0&2&2&0&0&2&0&0&2&0&0&0&15\\
573&2&0&1&0&0&2&0&4&0&0&0&0&0&2&2&0&2&0&0&0&15\\
1161&0&0&0&0&0&0&0&0&4&0&2&0&2&1&0&0&0&0&0&0&9\\
1339&0&0&0&0&0&0&0&0&4&2&0&0&2&2&0&0&1&0&0&0&11\\
1591&2&0&2&0&0&0&0&0&0&0&0&4&0&0&0&1&0&0&0&0&9\\
1951&0&0&0&0&0&0&0&0&0&0&2&0&0&0&0&0&0&4&0&1&7\\
133&0&0&0&0&0&0&2&0&0&0&2&0&0&0&0&0&0&4&0&1&9\\
5000&0&0&0&0&0&0&0&0&0&0&0&0&1&0&1&0&0&0&2&1&5\\
\hline
0&9999&10001&9998&0&10000&9999&10000&9999&10000&10000&9998&9999&10000&10000&10000&9999&10000&9998&10000&10000&189990\\
\end{array}\right)$
}
\vspace{0.2 in}

{\tiny
$\left(
\begin{array}{r|rrrrrrrrrrrrrrrrrrrr|r}
1033&0&4&0&0&0&0&0&2&0&2&0&0&0&0&0&0&2&0&1&0&11\\
1467&0&4&0&0&0&0&0&0&0&2&0&0&0&0&0&1&2&0&2&0&11\\
113&0&0&2&4&0&2&0&0&0&0&0&2&0&0&0&0&0&1&0&0&11\\
192&0&0&2&4&0&2&0&0&2&0&0&1&0&0&0&0&0&2&0&0&13\\
1262&0&0&0&4&0&2&0&2&0&0&0&2&0&0&1&0&0&0&2&0&13\\
933&0&0&2&4&0&2&0&2&2&0&0&2&0&0&1&0&0&0&2&0&17\\
1882&0&0&0&0&0&0&4&0&0&0&0&2&2&0&0&0&2&1&0&0&11\\
618&0&0&0&0&0&0&4&2&0&0&0&2&2&1&2&0&2&2&0&0&17\\
1083&2&0&0&0&0&0&0&2&1&0&2&0&0&4&0&0&0&0&0&0&11\\
264&0&0&2&0&0&0&0&0&2&0&0&0&0&4&0&0&0&1&0&0&9\\
858&0&0&2&0&0&2&0&0&1&0&0&0&0&4&0&0&0&2&0&0&11\\
141&2&0&2&0&0&0&0&1&2&0&0&0&0&4&0&0&0&2&0&0&13\\
1462&2&0&0&0&0&2&0&0&0&0&0&0&0&0&4&0&0&0&1&0&9\\
2037&2&0&0&0&0&0&0&0&0&1&2&0&0&0&0&4&0&0&0&0&9\\
96&2&0&0&0&0&0&0&0&0&1&2&2&0&0&0&4&0&2&0&0&13\\
1067&0&0&2&0&0&0&0&0&2&0&2&0&2&0&0&0&0&1&0&4&13\\
1433&0&0&2&0&0&0&0&0&2&2&1&0&2&0&0&0&0&2&0&4&15\\
\hline
0&9998&10000&10002&10000&0&10000&10000&9999&10001&9999&9999&10000&10000&10002&9999&9999&10000&10002&9999&10000&189999\\
\end{array}\right)$
}
\vspace{0.2 in}

{\tiny
$\left(
\begin{array}{r|rrrrrrrrrrrrrrrrrrrr|r}
5000&2&0&0&0&0&0&0&0&0&0&1&0&0&0&0&1&0&0&0&0&4\\
294&0&4&0&2&1&0&0&0&0&2&0&0&0&0&2&2&2&2&0&0&17\\
414&0&4&0&0&0&0&0&0&0&2&0&0&0&0&0&2&1&0&0&2&11\\
1792&0&4&0&0&2&0&0&0&0&2&0&0&0&0&1&2&2&0&0&2&15\\
1994&0&0&4&2&0&0&0&0&0&0&0&2&0&1&0&0&0&0&0&2&11\\
1012&0&0&2&4&1&0&0&0&2&0&0&2&0&0&2&0&0&2&0&0&15\\
349&0&0&0&0&0&0&4&2&0&0&0&2&2&2&1&0&2&0&0&0&15\\
1012&0&0&0&0&2&0&4&0&0&0&0&1&2&0&2&0&0&2&0&0&13\\
1138&0&0&0&0&0&0&4&0&0&0&2&2&2&1&0&0&2&2&0&0&15\\
1018&0&0&0&0&0&0&0&4&0&0&2&0&0&0&1&0&0&0&0&0&7\\
687&0&0&0&2&0&0&0&4&0&2&1&0&0&0&2&0&0&0&0&0&11\\
235&0&0&0&0&0&0&0&4&0&2&0&0&0&0&1&0&2&0&0&0&9\\
1542&0&0&0&0&2&0&0&1&2&0&0&0&0&4&0&0&0&2&0&0&11\\
1375&0&0&0&0&0&0&0&0&2&2&0&0&2&0&0&0&1&0&4&0&11\\
297&0&0&0&0&0&0&0&0&2&1&0&0&2&0&2&0&2&0&4&0&13\\
1547&0&0&0&0&0&0&0&0&1&0&0&0&1&0&0&0&0&0&2&1&5\\
55&0&0&0&0&0&0&0&0&0&2&0&0&2&0&0&0&0&0&4&1&9\\
\hline
0&10000&10000&10000&9998&9998&0&9996&10000&9999&10001&9999&9998&9999&9998&9998&10000&9999&9996&10002&10002&189983\\
\end{array}\right)$
}

    \end{minipage}
    };
\end{tikzpicture}

\tikzstyle{mybox} = []
\begin{tikzpicture}[transform shape, rotate=90, baseline=6.8in]
\node [mybox] (box) {%
    \begin{minipage}[t!]{0.5\textwidth}

{\tiny
$\left(
\begin{array}{r|rrrrrrrrrrrrrrrrrrrr|r}
898&0&4&0&0&0&0&0&0&0&2&0&0&0&0&0&2&2&1&2&0&13\\
533&0&4&0&2&1&0&0&0&0&0&0&0&0&0&0&2&2&0&2&2&15\\
1068&0&4&0&0&2&0&0&0&0&0&0&0&0&0&0&2&2&1&2&2&15\\
1013&2&0&0&0&0&4&0&0&0&0&0&0&0&1&2&0&0&0&0&0&9\\
629&2&0&2&1&0&0&0&4&0&2&0&0&0&0&2&0&2&0&0&0&15\\
42&2&0&0&2&0&2&0&4&0&1&0&0&0&0&2&0&2&0&0&0&15\\
1443&0&0&2&1&0&2&0&4&0&2&0&0&0&2&2&0&2&0&0&0&17\\
386&2&0&0&0&0&1&0&4&0&0&2&0&0&2&2&0&2&0&0&0&15\\
2061&2&0&0&0&0&0&0&0&1&0&4&0&0&0&0&0&0&0&0&2&9\\
173&0&0&0&0&0&0&0&0&2&0&4&0&0&0&0&2&0&0&0&1&9\\
2500&0&0&2&2&2&0&0&0&0&0&0&4&0&0&0&1&0&0&0&0&11\\
2500&0&0&0&0&0&0&0&0&2&0&0&0&4&0&0&0&0&0&2&1&9\\
869&2&0&0&0&1&2&0&0&2&0&0&0&0&4&0&0&0&0&0&0&11\\
427&0&0&2&0&1&2&0&0&2&0&0&0&0&4&0&0&0&0&0&0&11\\
743&0&0&0&2&1&0&0&0&0&0&0&0&0&0&4&0&0&0&0&0&7\\
1863&0&0&0&0&0&0&0&0&0&2&0&0&0&0&0&1&0&4&0&0&7\\
145&0&0&0&2&2&0&0&0&0&2&2&0&0&1&0&2&0&4&0&0&15\\
\hline
0&10000&9996&9998&9998&9998&10000&0&10000&9999&9998&9998&10000&10000&10000&9998&9997&9998&9998&9998&9997&189971\\
\end{array}\right)$
}
\vspace{0.2 in}

{\tiny
$\left(
\begin{array}{r|rrrrrrrrrrrrrrrrrrrr|r}
1551&4&0&0&0&0&0&0&0&0&0&0&2&0&2&0&1&2&0&0&0&11\\
959&4&0&0&0&0&2&0&0&0&0&0&2&0&0&1&2&2&0&0&0&13\\
105&0&4&0&0&0&0&0&0&0&2&0&0&0&0&2&0&2&1&0&0&11\\
919&0&4&0&0&2&0&0&0&0&0&0&0&0&0&0&2&0&0&0&1&9\\
290&0&4&0&0&2&2&0&0&0&0&0&0&0&0&0&0&1&2&0&2&13\\
737&0&4&0&0&2&0&0&0&0&1&0&0&0&0&2&0&2&2&2&2&17\\
2170&0&0&4&2&0&0&0&0&0&0&0&0&0&0&0&0&0&0&0&1&7\\
660&0&0&2&4&2&2&0&0&1&0&0&2&0&0&0&0&0&2&2&0&17\\
755&0&0&0&4&2&0&0&0&0&0&0&0&0&0&2&0&0&2&1&0&11\\
962&0&0&0&0&0&4&0&0&0&0&0&0&0&2&1&0&0&0&2&0&9\\
55&0&0&0&0&2&0&4&0&0&0&0&2&2&1&2&0&0&0&0&0&13\\
633&0&0&0&0&2&1&4&0&0&0&0&0&2&0&2&0&2&0&0&0&13\\
898&0&2&0&0&2&0&4&0&0&0&0&2&2&0&2&0&1&0&0&0&15\\
861&0&0&0&0&0&2&4&0&0&0&0&2&2&0&2&0&1&0&0&0&13\\
52&0&0&0&0&2&0&4&0&0&0&0&1&2&2&0&0&0&2&0&0&13\\
2263&0&0&0&0&0&0&0&0&2&4&0&0&0&0&0&2&0&0&2&1&11\\
2407&0&0&0&0&0&0&0&0&2&0&4&0&2&2&0&0&0&2&0&1&13\\
93&0&0&0&0&0&0&0&0&0&0&4&0&2&0&0&2&0&1&0&2&11\\
\hline
0&10000&10000&10000&10000&9998&10001&9996&0&10000&9999&10000&10000&9998&9999&9999&9999&9999&10000&9999&9999&189986\\
\end{array}\right)$
}
\vspace{0.2 in}

{\tiny
$\left(
\begin{array}{r|rrrrrrrrrrrrrrrrrrrr|r}
308&0&2&2&4&2&0&0&2&0&0&0&2&0&0&2&0&0&1&0&0&17\\
188&0&2&1&4&0&2&0&2&0&0&0&2&0&0&2&0&0&2&0&0&17\\
48&0&2&0&0&0&0&0&1&0&4&0&0&0&0&0&0&0&0&0&0&7\\
708&0&0&0&0&0&0&0&0&0&2&0&0&0&0&0&1&0&1&0&0&4\\
100&2&0&0&0&0&0&0&1&0&0&4&0&0&0&0&2&0&0&0&2&11\\
373&2&0&0&0&0&0&2&2&0&0&4&0&0&0&0&2&0&0&0&1&13\\
992&0&0&0&0&0&0&0&0&0&0&0&1&2&0&0&0&0&0&0&0&3\\
69&2&0&0&0&2&2&2&0&0&0&1&0&0&4&0&0&0&0&0&0&13\\
402&2&0&2&0&2&2&2&0&0&0&0&0&0&4&0&0&0&1&0&0&15\\
25&0&0&0&0&2&0&0&2&0&0&0&0&0&4&0&0&0&1&0&0&9\\
200&0&2&0&0&0&0&1&0&0&0&0&0&0&0&0&0&4&0&0&0&7\\
48&2&2&0&0&0&0&2&1&0&0&0&0&0&0&0&0&4&0&0&0&11\\
326&0&1&0&0&0&0&0&0&0&0&0&0&0&0&1&0&1&0&2&1&6\\
319&0&0&0&0&0&2&0&0&0&0&0&0&0&0&2&0&2&0&4&1&11\\
14&0&2&0&0&0&2&0&0&0&0&0&0&0&0&2&0&2&0&4&1&13\\
23&0&2&2&0&2&0&0&0&0&2&1&0&0&0&0&0&0&0&0&4&13\\
165&0&0&2&0&2&0&0&0&0&2&0&0&0&0&0&2&0&1&0&4&13\\
\hline
0&1984&1984&1984&1984&1984&1984&1984&1984&0&1984&1984&1984&1984&1984&1984&1984&1984&1984&1984&1984&37696\\
\end{array}\right)$
}

    \end{minipage}
    };
\end{tikzpicture}

\tikzstyle{mybox} = []
\begin{tikzpicture}[transform shape, rotate=90, baseline=6.8in]
\node [mybox] (box) {%
    \begin{minipage}[t!]{0.5\textwidth}

{\tiny
$\left(
\begin{array}{r|rrrrrrrrrrrrrrrrrrrr|r}
255&0&0&0&0&2&0&0&0&0&0&2&0&0&1&0&0&0&0&0&4&9\\
1444&0&0&2&0&0&1&0&0&0&0&0&1&2&0&0&0&0&0&0&4&10\\
367&0&0&2&0&0&4&2&0&0&0&0&1&4&0&2&0&0&0&8&0&23\\
906&0&2&0&0&1&0&0&0&0&0&0&0&0&0&0&0&1&0&0&0&4\\
186&0&2&2&0&1&4&8&16&0&0&8&2&4&8&0&1&0&0&2&0&58\\
1038&0&4&0&1&0&0&0&0&0&0&0&0&0&0&2&0&2&0&0&0&9\\
896&1&0&0&0&0&0&0&0&0&0&2&0&1&0&0&0&0&4&0&0&8\\
130&1&4&0&8&8&4&4&0&0&0&0&2&0&8&0&0&0&16&0&0&55\\
806&2&0&0&4&0&0&1&0&2&0&2&2&4&0&4&0&4&0&8&2&35\\
1128&2&0&1&0&1&1&0&0&0&0&0&2&0&0&0&4&0&0&0&0&11\\
118&2&0&2&0&4&0&0&0&0&0&0&2&1&4&0&0&0&8&0&0&23\\
81&2&0&2&4&0&0&0&0&8&0&4&1&0&0&0&0&4&2&0&0&27\\
570&2&2&4&1&1&4&4&8&0&0&1&0&0&2&4&0&4&0&0&2&39\\
29&4&1&0&0&0&0&0&0&16&0&0&0&0&8&2&0&8&0&0&0&39\\
392&4&1&2&8&4&4&2&1&16&0&0&4&0&8&2&0&0&2&0&0&58\\
94&4&1&8&2&0&0&4&16&4&0&8&0&0&8&0&0&2&0&0&0&57\\
91&4&16&0&1&8&8&8&0&2&0&2&0&0&4&4&0&4&4&2&0&67\\
83&8&0&8&4&32&0&16&4&4&0&32&16&8&2&0&64&2&2&0&4&206\\
30&16&1&0&2&8&4&32&8&4&0&4&16&0&32&16&0&8&64&2&4&221\\
\hline
0&10000&9997&10000&10003&10004&10000&10004&10004&10006&0&10006&10004&10002&10005&10000&10010&10000&10004&9998&10000&190047\\
\end{array}\right)$
}
\vspace{0.2 in}

{\tiny
$\left(
\begin{array}{r|rrrrrrrrrrrrrrrrrrrr|r}
641&4&0&0&0&0&2&0&0&0&0&0&1&0&0&0&2&0&0&0&0&9\\
1859&4&0&0&0&0&2&0&0&0&0&0&2&0&0&1&2&0&0&0&0&11\\
476&0&0&0&0&2&0&4&0&0&0&0&2&0&1&2&0&2&0&0&0&13\\
605&0&0&0&0&2&2&4&0&0&0&0&2&0&2&2&0&1&0&0&0&15\\
540&0&2&0&0&2&1&4&0&0&0&0&0&0&0&2&0&0&2&0&0&13\\
18&0&2&2&0&0&2&2&4&0&1&0&0&0&2&0&0&0&0&0&0&15\\
1855&0&2&2&2&0&0&0&4&0&1&0&0&0&0&2&0&2&0&0&0&15\\
563&0&0&2&1&0&2&0&4&0&2&0&0&0&0&2&0&2&0&0&0&15\\
65&0&2&2&2&0&2&0&4&0&2&0&0&0&0&1&0&2&0&0&0&17\\
1740&0&0&0&2&0&0&0&0&4&2&0&0&1&2&0&0&2&0&2&2&17\\
1740&0&0&0&0&0&0&2&0&1&0&0&2&4&0&0&0&0&0&2&0&11\\
980&0&0&2&0&1&2&0&0&0&0&0&0&0&4&0&0&0&0&0&0&9\\
408&0&0&0&2&1&0&0&0&0&0&0&0&0&0&0&2&0&4&0&0&9\\
171&0&2&0&0&0&0&0&0&0&1&0&0&0&0&0&2&0&4&0&0&9\\
619&0&2&0&0&2&0&0&0&0&1&0&0&0&0&0&2&0&4&0&0&11\\
212&0&2&0&2&1&0&0&0&0&2&0&0&0&0&0&2&0&4&0&0&13\\
440&0&0&0&2&2&0&0&0&0&2&0&0&0&2&0&2&0&4&0&1&15\\
870&0&2&2&0&2&0&0&0&0&0&0&0&0&0&0&0&0&1&2&4&13\\
650&0&2&2&0&2&0&0&0&2&2&0&0&2&0&0&2&0&1&2&4&21\\
\hline
0&10000&10000&10002&10003&10000&10002&10000&10004&10000&10003&0&10001&10000&10002&10002&10000&10003&10000&10000&10000&190022\\
\end{array}\right)$
}
\vspace{0.2 in}

{\tiny
$\left(
\begin{array}{r|rrrrrrrrrrrrrrrrrrrr|r}
190&0&0&0&0&0&0&0&0&1&0&2&0&0&0&0&4&0&1&0&2&10\\
270&0&0&0&0&2&0&0&0&4&8&16&0&32&2&0&1&0&4&16&16&101\\
86&0&0&0&0&2&0&1&1&4&2&0&0&8&0&0&2&2&2&0&4&28\\
215&0&0&32&0&0&16&0&0&8&4&0&0&0&16&0&0&0&0&2&1&79\\
76&0&0&32&2&0&0&4&16&0&8&2&0&0&0&2&4&8&8&4&16&106\\
314&0&4&1&0&0&0&2&2&0&0&1&0&0&2&0&8&0&4&0&2&26\\
391&0&4&1&0&8&2&0&2&0&2&2&0&0&0&4&4&2&4&2&0&37\\
274&0&8&0&16&0&4&1&8&8&4&1&0&0&0&8&0&0&2&0&0&60\\
297&0&8&0&16&4&8&0&8&2&4&0&0&1&4&2&0&4&8&8&0&77\\
44&1&0&0&2&0&16&32&8&4&2&16&0&8&8&16&1&2&0&4&4&124\\
137&1&0&0&2&2&2&8&4&2&2&0&0&0&1&4&0&4&0&2&0&34\\
116&1&8&0&0&16&0&2&4&4&0&2&0&0&8&0&0&4&0&0&8&57\\
365&4&4&0&0&4&0&16&2&4&1&2&0&0&1&8&0&8&2&0&0&56\\
154&8&0&0&0&0&1&0&0&0&8&0&0&0&0&0&16&4&0&2&8&47\\
219&8&0&0&0&0&4&0&0&0&0&0&0&0&0&4&0&1&0&2&0&19\\
3&8&0&0&0&1&0&0&4&4&8&2&0&0&2&2&16&0&8&4&4&63\\
111&8&2&0&2&4&2&1&4&4&8&8&0&0&0&1&16&2&8&4&0&74\\
5&32&0&0&0&128&16&4&8&1&0&32&0&2&64&64&16&16&8&0&4&395\\
130&32&0&0&1&2&0&0&1&8&2&8&0&0&16&0&0&16&4&1&4&95\\
\hline
0&9973&9998&10017&10002&9965&10004&9999&10000&9995&9997&9982&0&9987&9984&9987&9994&9987&9996&9994&9991&189852\\
\end{array}\right)$
}

    \end{minipage}
    };
\end{tikzpicture}

\tikzstyle{mybox} = []
\begin{tikzpicture}[transform shape, rotate=90, baseline=6.8in]
\node [mybox] (box) {%
    \begin{minipage}[t!]{0.5\textwidth}

{\tiny
$\left(
\begin{array}{r|rrrrrrrrrrrrrrrrrrrr|r}
68&0&0&0&0&0&1&2&1&0&0&0&0&0&1&1&0&1&1&0&0&8\\
61&0&2&0&0&2&2&4&2&0&0&0&0&0&1&2&0&2&2&0&0&19\\
40&0&0&0&2&0&0&0&0&4&1&0&0&0&2&0&0&0&0&0&0&9\\
19&2&0&0&0&0&0&0&2&0&0&4&0&0&1&0&0&0&0&0&0&9\\
76&2&0&0&0&0&0&0&2&0&0&4&0&0&2&0&2&0&1&0&0&13\\
9&2&0&2&1&2&0&0&0&0&0&0&4&0&0&0&0&0&0&0&0&11\\
86&2&0&2&2&1&0&0&0&0&0&0&4&0&0&0&2&0&0&0&0&13\\
46&0&0&0&1&0&1&0&0&0&1&0&0&0&0&1&0&1&0&2&0&7\\
1&0&2&0&2&0&2&0&0&0&1&0&0&0&0&2&0&2&0&4&0&15\\
71&0&2&0&1&0&2&0&0&2&2&0&0&0&0&2&0&2&0&4&0&17\\
38&0&0&2&0&2&0&0&0&0&2&0&0&0&0&0&1&0&0&0&4&11\\
39&0&2&2&0&2&0&0&0&2&1&0&0&0&0&0&0&0&2&0&4&15\\
18&0&2&2&0&0&0&0&0&0&2&0&0&0&0&0&1&0&2&0&4&13\\
\hline
0&380&380&380&380&380&380&380&380&380&380&380&380&0&380&380&380&380&380&380&380&7220\\
\end{array}\right)$
}
\vspace{0.2 in}

{\tiny
$\left(
\begin{array}{rrrrrrrrrrrrrrrrrrrrrr}
260&4&0&0&0&0&0&0&0&0&0&0&0&0&0&0&2&1&0&0&0&7\\
1298&4&0&0&0&0&0&0&0&0&0&0&2&0&0&0&2&1&0&0&0&9\\
866&0&0&4&2&0&0&0&0&0&0&0&1&0&0&0&0&0&0&0&0&7\\
344&0&0&0&0&4&0&0&0&0&0&0&2&0&0&0&2&0&1&0&0&9\\
1289&0&2&0&0&4&0&0&0&0&0&0&1&0&0&2&0&0&0&0&2&11\\
867&0&2&0&2&4&0&0&0&0&0&0&2&0&0&2&0&0&0&0&1&13\\
1929&0&0&0&1&0&2&0&0&0&0&0&0&0&0&0&0&0&0&1&0&4\\
1535&2&0&2&0&0&4&0&2&0&0&0&0&0&0&1&0&0&0&2&0&13\\
1024&0&0&0&0&0&0&4&0&0&0&0&2&2&0&1&0&0&0&0&0&9\\
778&0&2&0&0&0&0&4&0&0&0&0&1&2&0&0&0&2&0&0&0&11\\
698&0&0&0&0&0&0&4&0&0&0&2&0&1&0&2&0&2&2&0&0&13\\
1732&0&0&2&2&0&0&0&4&0&2&0&0&0&0&1&0&2&0&0&0&13\\
2974&0&0&0&0&0&0&0&0&2&1&0&0&1&0&0&0&0&0&1&1&6\\
1013&0&0&0&0&0&0&0&0&4&0&0&0&2&0&0&0&2&0&2&1&11\\
349&2&0&0&0&0&0&0&0&0&0&4&0&2&0&0&1&0&0&0&0&9\\
1428&0&0&0&0&0&0&0&0&0&0&4&0&0&0&0&2&0&0&0&1&7\\
1496&0&2&0&0&0&0&0&0&0&2&1&0&0&0&0&2&0&4&0&0&11\\
570&0&2&0&2&0&0&0&0&0&1&0&0&0&0&0&0&0&4&0&2&11\\
0&10000&10000&9998&9999&10000&9998&10000&9998&10000&10000&10000&9999&10000&0&9999&10001&10000&10004&9999&10000&189995\\
\end{array}\right)$
}
\vspace{0.2 in}

{\tiny
$\left(
\begin{array}{r|rrrrrrrrrrrrrrrrrrrr|r}
132&4&0&0&0&0&0&0&0&0&0&2&0&0&1&0&2&2&0&0&0&11\\
54&0&4&0&0&0&2&0&0&0&0&0&0&0&0&0&0&1&0&0&2&9\\
58&0&4&0&0&0&0&0&0&0&0&0&0&0&0&0&2&0&1&0&2&9\\
52&0&4&0&0&0&2&0&0&0&1&0&0&0&0&0&2&2&2&0&2&15\\
164&0&0&2&4&0&0&0&0&2&0&0&1&0&0&0&0&0&2&0&0&11\\
48&0&0&0&0&4&0&0&0&0&0&0&2&0&1&0&0&0&0&0&2&9\\
13&0&0&0&0&4&0&0&0&0&0&0&0&0&2&0&0&0&1&0&2&9\\
69&0&0&0&0&4&0&0&0&0&0&0&0&0&0&0&2&0&1&0&2&9\\
34&0&0&0&0&4&0&0&0&0&0&0&2&0&0&0&1&0&2&0&2&11\\
148&0&0&0&0&0&2&4&0&0&0&1&2&2&2&0&0&0&0&0&0&13\\
16&0&0&0&0&0&0&4&0&0&0&2&2&2&0&0&0&2&1&0&0&13\\
26&0&0&2&0&0&0&0&4&0&0&0&0&0&1&0&0&0&0&0&0&7\\
74&0&0&2&0&0&2&0&4&0&2&2&0&0&0&0&0&1&0&0&0&13\\
64&2&0&2&0&0&0&0&4&0&2&1&0&0&2&0&0&2&0&0&0&15\\
328&0&0&0&0&0&0&0&0&1&1&0&0&1&0&0&0&0&0&2&0&5\\
\hline
0&656&656&656&656&656&656&656&656&656&656&656&656&656&656&0&656&656&656&656&656&12464\\
\end{array}\right)$
}

    \end{minipage}
    };
\end{tikzpicture}

\tikzstyle{mybox} = []
\begin{tikzpicture}[transform shape, rotate=90, baseline=6.8in]
\node [mybox] (box) {%
    \begin{minipage}[t!]{0.5\textwidth}

{\tiny
$\left(
\begin{array}{r|rrrrrrrrrrrrrrrrrrrr|r}
80&4&0&0&0&0&2&0&1&0&0&0&0&0&0&2&0&2&0&0&0&11\\
6&0&4&0&0&0&2&1&0&0&0&0&0&0&0&2&0&2&2&0&0&13\\
40&0&4&0&0&0&2&0&2&0&1&0&0&0&0&0&0&2&2&2&0&15\\
26&0&4&0&0&0&2&1&0&0&0&0&0&0&0&2&0&2&2&2&0&15\\
8&0&4&0&0&0&2&0&2&0&0&0&0&0&0&2&0&2&2&1&0&15\\
80&0&0&0&2&4&0&0&0&0&0&0&0&0&2&1&0&0&0&0&0&9\\
50&0&0&0&0&0&0&0&2&1&4&0&0&0&0&0&0&0&0&2&0&9\\
22&0&0&0&0&0&0&2&1&2&0&4&0&2&0&0&0&0&0&0&0&11\\
3&0&0&0&0&0&0&0&0&0&0&4&0&2&1&0&0&0&0&0&0&7\\
22&0&0&0&0&0&0&2&1&0&0&4&0&2&2&0&0&0&0&0&0&11\\
33&0&0&0&0&0&0&0&0&2&0&4&0&2&1&0&0&0&0&0&0&9\\
80&0&0&2&1&0&0&2&0&0&0&0&4&0&0&0&0&0&0&0&0&9\\
40&0&0&0&2&0&0&1&0&0&2&0&0&0&2&0&0&0&4&0&0&11\\
80&0&0&2&0&0&0&0&0&2&0&0&0&2&0&0&0&0&0&1&4&11\\
\hline
0&320&320&320&320&320&320&320&320&320&320&320&320&320&320&320&0&320&320&320&320&6080\\
\end{array}\right)$
}
\vspace{0.2 in}

{\tiny
$\left(
\begin{array}{r|rrrrrrrrrrrrrrrrrrrr|r}
14&0&0&0&0&0&8&0&16&0&0&0&0&0&0&0&1&0&0&4&2&31\\
26&0&0&0&0&1&0&0&0&4&16&8&8&16&0&0&2&0&4&32&8&99\\
1551&0&0&0&2&1&0&0&0&4&2&0&0&2&0&0&1&0&1&0&2&15\\
550&0&0&4&0&4&2&16&1&4&2&4&8&8&8&8&4&0&4&0&2&79\\
85&0&0&4&2&4&4&0&0&16&0&8&2&1&0&1&4&0&0&8&8&62\\
479&0&0&4&4&0&2&0&8&0&0&4&0&0&2&1&0&0&0&0&0&25\\
37&0&0&4&8&0&8&0&0&0&8&1&4&2&0&0&0&0&0&16&8&59\\
430&0&1&4&4&0&4&0&8&0&4&0&0&0&4&2&0&0&0&2&1&34\\
44&0&1&16&8&0&32&8&2&0&0&1&8&32&8&4&0&0&16&64&32&232\\
194&0&2&0&8&4&0&1&4&0&4&2&2&0&2&8&0&0&4&16&8&65\\
99&0&16&0&1&8&0&0&0&0&8&2&0&0&0&8&4&0&2&0&0&49\\
13&0&16&0&2&0&0&0&1&0&0&0&0&0&0&0&8&0&4&0&8&39\\
109&0&16&0&2&4&0&0&1&0&4&2&0&0&1&0&8&0&4&0&8&50\\
130&0&32&4&4&16&16&4&4&1&4&2&0&4&8&8&0&0&16&8&0&131\\
715&8&2&2&0&2&0&0&0&0&1&4&4&0&0&0&4&0&2&0&0&29\\
112&16&0&4&0&0&8&0&4&0&1&8&2&0&8&0&4&0&4&0&2&61\\
130&16&0&4&0&2&8&1&0&0&0&0&8&0&0&4&8&0&0&0&0&51\\
14&16&0&4&1&4&4&0&0&0&0&0&8&0&8&0&8&0&0&2&0&55\\
6&32&0&0&2&8&1&0&0&0&0&16&16&0&4&16&0&0&4&0&0&99\\
\hline
0&10008&9988&10002&9997&9995&10012&9996&10000&9998&9985&10001&9998&10005&9999&10000&9989&0&10003&10008&10004&189988\\
\end{array}\right)$
}
\vspace{0.2 in}

{\tiny
$\left(
\begin{array}{r|rrrrrrrrrrrrrrrrrrrr|r}
1598&0&4&0&0&0&2&0&0&0&2&0&0&0&0&2&0&2&0&1&0&13\\
710&0&4&0&0&2&0&0&2&0&1&0&0&0&0&2&0&2&0&2&0&15\\
192&0&4&0&0&0&0&0&0&0&2&0&0&0&0&0&1&2&0&0&2&11\\
305&0&0&2&4&0&0&0&2&0&0&0&2&0&0&1&0&0&0&0&0&11\\
2047&0&0&2&4&0&0&0&0&2&0&0&0&0&0&1&0&0&0&2&0&11\\
1516&0&0&0&0&4&0&0&0&0&0&0&2&0&0&2&1&0&0&0&0&9\\
295&0&0&0&2&4&0&0&0&0&0&0&2&0&0&0&0&0&0&0&1&9\\
2500&0&0&0&0&0&2&4&2&0&0&0&2&2&1&0&0&2&0&0&0&15\\
1201&0&0&0&0&0&0&0&2&0&4&0&0&0&0&0&1&0&0&2&0&9\\
2216&2&0&0&0&0&0&0&0&2&0&4&0&1&0&0&0&0&0&0&0&9\\
284&0&0&0&0&0&0&0&2&2&0&4&0&0&1&0&0&0&0&0&0&9\\
1804&2&0&2&0&0&1&0&0&0&0&0&0&0&4&0&0&0&0&0&0&9\\
596&2&0&0&0&1&0&0&0&0&0&0&0&0&0&0&4&0&0&0&0&7\\
384&2&0&0&0&1&0&0&0&0&0&0&2&0&0&0&4&0&0&0&0&9\\
453&0&0&0&0&0&0&0&0&2&2&0&0&0&0&0&1&0&0&0&4&9\\
1034&0&0&0&0&0&0&0&0&0&0&0&0&2&0&0&1&0&0&0&4&7\\
358&0&0&2&0&1&0&0&0&0&0&0&0&2&0&0&2&0&0&0&4&11\\
486&0&0&2&0&0&0&0&0&0&0&0&0&0&0&0&2&0&0&1&4&9\\
\hline
0&10000&10000&10000&9998&10002&10000&10000&10000&10000&10000&10000&10000&10000&10000&10000&10004&10000&0&10000&10003&190007\\
\end{array}\right)$
}

    \end{minipage}
    };
\end{tikzpicture}

\tikzstyle{mybox} = []
\begin{tikzpicture}[transform shape, rotate=90, baseline=6.8in]
\node [mybox] (box) {%
    \begin{minipage}[t!]{0.5\textwidth}

{\tiny
$\left(
\begin{array}{r|rrrrrrrrrrrrrrrrrrrr|r}
627&0&4&0&0&2&0&0&0&0&0&0&0&0&0&0&2&0&1&0&0&9\\
899&0&4&0&0&0&0&0&2&0&2&0&0&0&0&1&2&0&2&0&0&13\\
494&0&4&0&0&2&0&0&2&0&2&0&0&0&0&0&1&0&2&0&2&15\\
141&0&4&0&0&2&0&2&2&0&2&0&0&0&0&0&0&1&2&0&2&17\\
2500&0&0&2&4&1&0&0&2&0&0&0&2&0&0&0&0&0&2&0&0&13\\
1558&2&0&0&0&0&4&1&0&0&0&0&0&0&2&0&0&0&0&0&0&9\\
942&2&0&2&0&0&4&0&0&0&0&0&0&0&2&1&0&0&0&0&0&11\\
2500&0&0&0&0&0&0&0&0&4&0&2&0&0&2&0&0&1&0&0&0&9\\
553&0&2&0&0&0&0&0&1&0&4&0&0&0&0&0&2&0&0&0&0&9\\
1054&0&0&0&0&0&0&0&0&0&4&0&0&0&0&0&2&0&1&0&0&7\\
126&0&2&0&0&0&0&0&0&0&4&0&0&0&0&0&1&0&2&0&0&9\\
2500&0&0&0&0&0&0&2&0&0&0&2&2&4&0&0&0&0&0&0&1&11\\
660&2&0&0&0&1&0&2&0&0&0&0&0&0&0&4&0&0&0&0&0&9\\
1379&0&0&0&0&2&0&0&1&0&0&0&0&0&0&4&0&0&0&0&0&7\\
1840&2&0&0&0&0&0&1&0&0&0&0&0&0&0&0&0&4&0&0&0&7\\
1558&0&0&2&0&1&0&0&0&0&0&0&0&0&0&0&2&0&0&0&4&9\\
\hline
0&10000&10002&10000&10000&10000&10000&10000&10000&10000&10000&10000&10000&10000&10000&9997&10002&10001&10001&0&10002&190005\\
\end{array}\right)$
}
\vspace{0.2 in}

{\tiny
$\left(
\begin{array}{r|rrrrrrrrrrrrrrrrrrrr|r}
822&0&0&0&0&0&0&0&0&0&0&2&2&0&1&0&4&0&0&0&0&9\\
442&0&0&0&0&0&4&8&0&0&0&0&0&16&0&1&0&0&0&2&0&31\\
130&0&0&0&2&0&0&0&0&4&0&2&0&0&0&0&0&0&0&1&0&9\\
200&0&1&16&0&2&2&4&8&0&0&0&8&0&8&4&0&0&0&0&0&53\\
813&0&4&0&0&0&0&2&2&0&2&0&0&0&0&1&0&2&0&0&0&13\\
3&0&8&0&0&4&8&1&8&32&16&16&4&16&16&2&8&16&8&8&0&171\\
124&0&16&0&0&4&0&2&4&0&2&0&2&1&0&8&1&8&8&4&0&60\\
17&0&16&0&1&0&8&8&2&0&0&0&0&0&0&4&0&0&0&0&0&39\\
343&1&0&0&0&1&0&0&0&4&0&0&4&8&0&0&2&0&0&0&0&20\\
742&1&0&0&0&2&0&0&0&0&8&0&1&0&0&1&2&0&4&0&0&19\\
175&1&4&0&8&4&2&8&4&0&8&2&0&0&4&4&2&2&16&1&0&70\\
10&2&16&0&4&4&0&2&0&0&0&0&0&1&0&8&32&4&0&8&0&81\\
239&4&0&16&4&4&8&0&0&0&0&0&8&0&0&0&4&0&2&4&0&54\\
453&4&4&4&8&0&8&0&2&8&1&0&0&0&1&8&0&8&2&16&0&74\\
110&4&4&8&2&32&8&8&4&0&0&2&1&0&16&4&2&2&0&0&0&97\\
192&8&0&0&2&4&0&4&4&0&0&16&1&0&8&2&0&0&0&0&0&49\\
225&8&0&0&4&0&4&0&8&0&0&16&0&0&2&4&0&4&1&0&0&51\\
128&8&0&0&16&1&0&0&8&32&0&4&8&0&16&0&2&16&8&0&0&119\\
144&8&8&2&1&8&0&4&4&2&2&2&8&0&4&0&16&1&4&0&0&74\\
\hline
0&10000&9996&10004&9993&9999&9994&9993&9994&9996&9999&9994&10008&9998&9993&9991&10012&9992&9993&9993&0&189942\\
\end{array}\right)$
}

    \end{minipage}
    };
\end{tikzpicture}

\end{document}